\documentclass[12pt]{amsart}

\usepackage{amsmath,amscd,amssymb,amsfonts}
\usepackage{epsfig}

\newcommand{\pf}{\begin{proof}}
\newcommand{\epf}{\end{proof}}

\newtheorem{prop}{Proposition}[section]
\newtheorem{cor}[prop]{Corollary}
\newtheorem{thm}[prop]{Theorem}
\newtheorem{lem}[prop]{Lemma}
\newtheorem{conj}[prop]{Conjecture}

\theoremstyle{definition}
\newtheorem{defn}[prop]{Definition}
\newtheorem{exmp}[prop]{Example}

\theoremstyle{remark}
\newtheorem{rem}[prop]{Remark}

\numberwithin{equation}{section}

\newenvironment{enu}%
        {\begin{list}{\textup{(\arabic{enumi})}}%
		{\usecounter{enumi}
		\setlength{\labelwidth}{1em}
		\setlength{\labelsep}{1em}
		\setlength{\itemindent}{2em}%{\labelwidth+\labelsep}
		\setlength{\itemsep}{6pt}
		\setlength{\rightmargin}{0pt}
		\setlength{\leftmargin}{0pt}}}%
	{\end{list}}

\newenvironment{enu2}%
        {\begin{list}{(\alph{enumii})}%
		{\usecounter{enumii}
		\setlength{\labelwidth}{1em}
		\setlength{\labelsep}{1em}
		\setlength{\itemindent}{-4pt}
		\setlength{\itemsep}{2pt}
		\setlength{\rightmargin}{0pt}
		\setlength{\leftmargin}{4em}}}%
	{\end{list}}

\newenvironment{ite}%
        {\begin{list}{$\bullet$}%
		{%\usecounter{enumi}
		\setlength{\labelwidth}{.5em}
		\setlength{\labelsep}{.5em}
		\setlength{\itemindent}{1em}%{\labelwidth+\labelsep}
		\setlength{\itemsep}{6pt}
		\setlength{\rightmargin}{0pt}
		\setlength{\leftmargin}{1em}}}%
	{\end{list}}

\newcommand{\card}{\operatorname{card}}
\newcommand{\id}{\operatorname{id}}
\newcommand{\GL}{\operatorname{GL}}
\newcommand{\Aut}{\operatorname{Aut}}
\newcommand{\Out}{\operatorname{Out}}
\newcommand{\Int}{\operatorname{Int}}
\newcommand{\Hom}{\operatorname{Hom}}
\newcommand{\End}{\operatorname{End}}
\newcommand{\Ind}{\operatorname{Ind}}
\newcommand{\Inn}{\operatorname{Inn}}
\newcommand{\im}{\operatorname{Im}}
\newcommand{\ord}{\operatorname{ord}}

\def\eps{\varepsilon}
\def\toba{{{\mathfrak B}}}
\def\wtoba{{\widehat{\mathfrak B}}}

\newcommand{\sgn}{\operatorname{sgn}}
\newcommand{\fun}{\operatorname{Fun}}
\newcommand{\rdu}{{\mathbb G}_{\infty}}
\newcommand{\kk}{\mathbb C}
\newcommand{\ydskg}{{}^G_G\mathcal{YD}}

\newcommand{\trid}{\triangleright}
\newcommand{\Sim}{\mathbb S}
\newcommand{\Tre}{\mathbb B}
\newcommand{\Die}{\mathbb D}
\newcommand{\N}{\mathbb N}
\newcommand{\Z}{\mathbb Z}
\newcommand{\R}{\mathbb R}
\newcommand{\fp}{{\mathbb F}_p}
\newcommand{\fq}{{\mathbb F}_q}
\newcommand{\fpt}{\mathbb{F}_{p^t}}
\newcommand{\ral}[1]{\Z\{#1\}}
\newcommand{\qal}[1]{\Z(#1)}

\newcommand{\wkappa}{\widetilde\kappa}

\newcommand{\rs}{\mathfrak{R}}
\newcommand{\rsx}{\rs\vert_{X}}
\newcommand{\qs}{\mathfrak{Q}}
\newcommand{\qg}{\mathfrak{q}}
\newcommand{\qsx}{\qs\vert_{X}}
\newcommand{\cs}{\mathfrak{CrS}}

\newcommand{\fllad}[1]{\displaystyle\mathop{\longrightarrow}^{#1}}%flecha der ley. alta

\def\al{\alpha}
\def\be{\beta}
\def\gm{\gamma}
\def\dt{\delta}
\def\sgm{\sigma}
\def\ep{\epsilon}
\def\gx{\Aut_{\trid}(X)}
\def\gxiota{\Aut_{\trid^{\iota}}(X)}
\newcommand{\gxd}[1]{\Aut_{\trid}(#1)}
\def\g{{\mathcal G}}
\def\gax{\Inn_{\trid} (X)}
\def\gay{\Inn_{\trid} (Y)}
\def\gaxiota{\Inn_{\trid^{\iota}} (X)}
\def\orb{{\mathcal O}}
\def\forget{\mathfrak F}
\newcommand{\mdpd}[4]{\left(\begin{smallmatrix}#1 & #2 \\ #3 & #4\end{smallmatrix}\right)}

\begin{document}
\title{From racks to pointed Hopf algebras}
\author{Nicol\'as Andruskiewitsch}
\address{N.A.: \newline\indent
Facultad de Matem\'atica, Astronom\'\i a y F\'\i sica, \newline\indent
Universidad Nacional de C\'ordoba, \newline\indent
CIEM -- CONICET, \newline\indent
(5000) Ciudad Universitaria, \newline\indent
C\'ordoba, Argentina}
\email{andrus@mate.uncor.edu}
\urladdr{www.mate.uncor.edu/andrus}
\author{Mat\'{\i}as Gra\~na}
\address{M.A.: \newline\indent
FCEyN -- UBA, \newline\indent
Pab. I -- Ciudad Universitaria, \newline\indent
(1428) Buenos Aires, Argentina \newline\indent
\textup{Current address}: \newline\indent
MIT, Mathematics Department, \newline\indent
77 Mass. Ave., \newline\indent
02139 Cambridge, MA - USA}
\email{matiasg@math.mit.edu}
\thanks{This work was partially supported by  ANPCyT, CONICET, Agencia C\'ordoba Ciencia,
Fundaci\'on Antorchas, Secyt (UBA) and Secyt (UNC)}
%\date{ \today }
\begin{abstract} 
A fundamental step in the classification of 
finite-dimensional complex pointed Hopf algebras is the
determination of all finite-dimensional
Nichols algebras of braided vector spaces arising from groups.
The most important class of braided vector spaces arising from groups
is the class of braided vector spaces  $(\kk X, c^q)$, where
$X$ is a rack and $q$ is a $2$-cocycle on $X$ with values in $\kk^\times$.
Racks and cohomology of racks appeared also in the work of topologists.
This leads us to the study of the structure of racks, their cohomology groups and the
corresponding Nichols algebras. We will show advances in these three directions.
We classify simple racks in group-theoretical terms; we describe projections of racks
in terms of general cocycles; we introduce a general cohomology theory of racks
containing properly the existing ones. We introduce a ``Fourier transform" on
racks of certain type; finally, we compute some new examples of finite-dimensional
Nichols algebras.
\end{abstract}
\maketitle

\setcounter{tocdepth}{1}
\tableofcontents

%%%%%%%%%%%%%%%%%%%%%%%%%%%%%%%%%%%%%%%%%%%%%%%%%%%%%%%%%%%%%%%
%%%%%%%%%%%%%%%% Seccion  Introduction %%%%%%%%%%%%%%%%%%%%%%%%
%%%%%%%%%%%%%%%%%%%%%%%%%%%%%%%%%%%%%%%%%%%%%%%%%%%%%%%%%%%%%%%
\section*{Introduction}

\emph{1.}
This paper is about braided vector spaces arising from pointed Hopf algebras,
and their Nichols algebras. 
Our general reference for pointed Hopf algebras is \cite{AS2}. 
We shall work over the field $\kk$ of complex numbers;
many results below are valid over more general fields.
We denote by $\rdu$ the group of roots of unity of $\kk$.

\bigskip
\emph{2.}
The determination of all complex finite dimensional pointed Hopf algebras $H$
with group of group-likes $G(H)$ isomorphic to a fixed finite group $\Gamma$ is
still widely open. Even the existence of such Hopf algebras $H$
(apart from the group algebra $\kk\Gamma$) is unknown for many
finite groups $\Gamma$. If $\Gamma$ is abelian, substantial advances
were done via the theory of quantum groups at roots of unit \cite{AS-adv}.
The results can be adapted to a non-necessarily abelian group $\Gamma$:
if $(a_{ij})_{1 \le i,j \le \theta}$ is a finite Cartan matrix, 
$g_1, \dots, g_\theta$ are \emph{central} elements in 
$\Gamma$, and $\chi_1, \dots, \chi_\theta$ are 
multiplicative characters  of $\Gamma$, such that 
$\chi_i(g_j)\chi_j(g_i) = \chi_i(g_i)^{a_{ij}}$ 
(plus some technical hypotheses on the orders of $\chi_i(g_j)$),
then a finite-dimensional
pointed Hopf algebra $H$ with group $G(H) \simeq \Gamma$ can be constructed 
from this datum. Besides these, only a small number of examples
have appeared in print; in these examples, $\Gamma$ is $\mathbb S_{3}$, $\mathbb S_{4}$,
$\mathbb D_{4}$ (see \cite{MiS}), $\mathbb S_{5}$ (using \cite{FK}), 
$\mathbb A_{4} \times \mathbb Z/2$ (see \cite{gr}).

\bigskip
\emph{3.} An important invariant of a pointed Hopf algebra $H$
is its \emph{infinitesimal braiding}; this is a braided vector space,
that is, a pair $(V, c)$, where $V$ is a vector space
and $c\in \Aut (V\otimes V)$ is a solution of the braid equation: 
$(c\otimes \id)(\id \otimes c)(c\otimes \id) = (\id \otimes c)(c\otimes \id)
(\id \otimes c)$. Let $\toba(V)$ be the Nichols algebra of $(V,c)$,
see \cite{AS2}.
If $\dim H$ is finite, then $\dim \toba(V)$ is finite and divides
$\dim H$. Conversely, given a braided vector space $(V,c)$,
where $V$ is a Yetter-Drinfeld module over the group algebra
$\kk\Gamma$, then the Radford's biproduct, or bosonization,
$H = \toba(V) \# \kk\Gamma$ is a pointed Hopf algebra
with $G(H) \simeq \Gamma$. 
Thus, a fundamental problem is to determine 
the dimension of the Nichols algebra of a finite-dimensional
braided vector space.
We remark that the same braided vector space can arise
as the infinitesimal braiding of  pointed Hopf algebras
$H$ with very different $G(H)$, 
as in the examples with Cartan matrices above.
Analogously, there are  infinitely many 
finite-dimensional pointed Hopf algebras $H$,
with non-isomorphic groups $G(H)$ and the same
infinitesimal braiding as, respectively, 
the examples above related to the
groups $\mathbb S_{3}$, $\mathbb S_{4}$,
$\mathbb D_{4}$, $\mathbb S_{5}$, $\mathbb A_{4} \times \mathbb Z/2$;
and such that they are link-indecomposable in the sense of \cite{MiS}, \cite{m}.

\bigskip
\emph{4.}  It is then more convenient to study
in a first stage Nichols algebras of braided vector spaces
of group-type \cite[Def. 1.4.10]{gr}. However, there are strong constraints 
on a braided vector space of group-type to have a finite-dimensional
Nichols algebra  \cite[Lemma 3.1]{gr}.
Briefly, we shall consider in this paper
the class of braided vector spaces of the form $(\kk X, c^q)$; where $(X, \trid)$
is a finite rack, $q$ is a $2$-cocycle with values on the multiplicative
group of invertible elements of $\kk$, and $c^q$ is given by
$c^q(i\otimes j) = q_{ij} \, i\trid j \otimes i$ (see the precise definitions
in the main part of the text). 
See also the discussion in \cite[Chapter 5]{bariloche}.

We are faced with the following questions:
to determine the general structure of finite racks;
to compute their second cohomology groups; and to decide whether
the dimensions of the corresponding Nichols algebras are finite.

\bigskip
\emph{5.} We now describe the contents of this paper.
The notions of ``rack" and ``quandle",
sets provided with a binary operation like the conjugation of a group,
have been considered in the literature, mainly as
a way to produce knot invariants (cf. \cite{Bk,jo,ma,k,fr,quandle1,de}).
Section \ref{sn:1} is a short survey of the theory of racks and quandles,
addressed to non-specialists on these
structures, including a variety of examples relevant for this paper.

The determination of all finite racks is a very hard task. 
There are two successive approximations to this problem. 
First, any finite rack is a union of indecomposable components.
However, indecomposable racks can be put together in many different ways,
and the description of all possible ways is again very difficult.
In other words, even the determination of all indecomposable finite racks
would not solve the general question.

In Section \ref{sn:2},  we describe epimorphisms of racks and quandles
by general cocycles. We  then introduce modules over a quandle, resp. a rack, 
$X$. We show that modules over $X$ are in one-to-one correspondence
with the abelian group objects
in the category of arrows over $X$, if $X$ is indecomposable.
Our definition of modules over $X$ generalizes those in \cite{ces,cens}.

We say that a non-trivial rack is \emph{simple} if it has no proper quotients.
Then any indecomposable finite rack with cardinality $>1$ is an extension of
a simple rack. We study simple racks in Section \ref{sn:3}. 
One of our main results is the explicit classification of all
finite simple racks, see Theorems \ref{clasipotprimo} and \ref{clasipqr}. 
The proof is based on a group-theoretical result
kindly communicated to us by R. Guralnick.

After acceptance of this paper, it became to our attention Joyce's article
\cite{jo2}, where results similar to those in Subsection \ref{ssn:sr} are
obtained. However, notice that our classification in
Theorem~\ref{clasipotprimo} includes more quandles than
that of \cite[Thm. 7 (2)]{jo2}. This is one of the reasons why we decided to
leave our results. The other reason is that, by using a result in \cite{egs}
(which depends upon the classification of simple groups), we can
split the simple quandles into two classes regarding their cardinality.
These classes appear naturally in \cite{jo2} also, but the fact that
they are split by cardinality was impossible to prove in 1982.

It is natural to define homology and cohomology
theories of racks and quandles as standard homology and cohomology
theories for abelian group objects in the category of arrows over $X$
\cite{q}. We propose in Section \ref{sn:4} a complex that, conjecturally,
would be suitable to compute these homology and cohomology theories.
We show that the homology and cohomology theories known so far
(see \cite{fr,quandle1,gr,cens}) are special cases
of ours. We discuss as well nonabelian cohomology theories.

A braided vector space of the form $(\kk Y, c^q)$ \emph{does not}
determine the rack $Y$ and the $2$-cocycle $q$. We provide a
general way of constructing two braided vector spaces 
$(\kk Y, c^q)$ and $(\kk \widetilde Y, c^{\widetilde q})$ where
the racks $Y$ and $\widetilde Y$ are not isomorphic in general,
but such that the corresponding Nichols algebras have the same dimension. 
In our construction, $Y$ is an extension $X\times_{\alpha} A$,
where $A$ is an $X$-module; $\widetilde Y$ is an extension 
$X\times_{\beta} \widehat A$,
where $\widehat A$, the group of characters of $A$, is also an $X$-module.
The construction can be thought as a Fourier transform. We show how to use
this construction to obtain new examples of pointed Hopf algebras
with non-abelian group of group-like elements.
This is the content of Section \ref{sn:5}.

In Section \ref{sn:6}, we present several new examples of finite dimensional
Nichols algebras $\toba(V)$ over finite groups.
First, we show that some Nichols algebras can be computed
by reduction to Nichols algebras of diagonal type, via Fourier transform. 
Next we use Fourier transform again to compute a Nichols algebra
related to the faces of the cube, starting from a Nichols algebra
related to the transpositions of $\mathbb S_4$ computed in \cite{MiS}. 
Finally we establish some relations that hold in Nichols algebras
related to affine racks, and use them to compute Nichols algebras
related to the vertices of the tetrahedron (a result announced in \cite{gr})
and the affine rack $(\mathbb Z/5, \trid^2)$. 
Support to our proofs is given by Theorem \ref{strate}
which gives criteria to insure that a finite dimensional braided
Hopf algebra is a Nichols algebra.

In most of the paper, we shall only consider finite racks,
or quandles, or crossed sets, and
omit the word ``finite" when designing them, unless explicitly
stated.

\bigskip
\emph{6.}
In conclusion, we remark that the next natural step in the classification
of finite dimensional pointed Hopf algebras is to deal with
Nichols algebras of braided vector spaces arising from \emph{simple} racks.

\vfill\eject
%%%%%%%%%%%%%%%%%%%%%%%%%%%%%%%%%%%%%%%%%%%%%%%%%%%%%%%%%%%%%%%
%%%%%%%%%%%%%%%% Seccion  Prelimin. %%%%%%%%%%%%%%%%%%%%%%%%%%%
%%%%%%%%%%%%%%%%%%%%%%%%%%%%%%%%%%%%%%%%%%%%%%%%%%%%%%%%%%%%%%%
\section{Preliminaries}\label{sn:1}

%%%%%%%%%%%%%%%%%%%%%%%%%%%%%%%%%%%%%%
%%%%%% SubSeccion  R/Q/CS     %%%%%%%%
%%%%%%%%%%%%%%%%%%%%%%%%%%%%%%%%%%%%%%
\subsection{Racks, quandles and crossed sets}\label{ss:rqcs}

\begin{defn}\label{df:rack}
A \emph{rack} is a pair $(X,\trid)$ where $X$ is a non-empty set and
$\trid:X\times X\to X$ is a function, such that
\begin{align}
\label{ccc1}& \phi_i:X\to X, \quad \phi_i (j) = i\trid j, \quad
	\text{is a bijection for all $i\in X$}, \\
\label{ccc4} & i\trid(j\trid k)=(i\trid j)\trid(i\trid k)\ \forall i,j,k\in X. \\
\intertext{A \emph{quandle} is a rack $(X,\trid)$ which further satisfies}
\label{ccc2} &i\trid i=i, \text{ for all } i\in X. \\
\intertext{A \emph{crossed set} is a quandle $(X,\trid)$ which further satisfies}
\label{ccc3}& j\trid i=i \text{ whenever } i\trid j=j.
\end{align}
\end{defn}

A morphism of racks is defined in the obvious way:
$\psi: (X,\trid) \to (Y,\trid)$ is a morphism of racks if
$\psi(i\trid j) = \psi(i) \trid \psi(j)$, for all $i, j \in X$.
Morphisms of quandles (resp. crossed sets) are morphisms of
racks between quandles (resp. crossed sets).

In particular, a subrack of a rack $(X,\trid)$ is a non-empty
subset $Y$ such that $Y \trid Y = Y$. If $X$ is a crossed set and
$Y$ is a subrack, then, clearly, it is a crossed subset; same for
quandles.

\begin{defn}
If $\Gamma$ is a group, any non-empty subset $X\subseteq\Gamma$
stable under conjugation by $\Gamma$ (i.e, a union of conjugacy
classes) is a crossed set with the structure $i\trid j=iji^{-1}$.
A crossed set isomorphic to one of these shall be called \emph{standard}.
\end{defn}
A primary goal is to compare arbitrary crossed sets with standard ones.

\begin{defn}\label{de:ddp}
Let $(X,\trid)$ be a rack and let $\Sim_X$ denote the group of symmetries of $X$.
By \eqref{ccc1}, we have a map $\phi: X \to \Sim_X$. Let $(X,\trid)$ be a rack. We set
\begin{align*}
\gx &:= \{g\in  \Sim_X: g(i\trid j) = g(i) \trid g(j)\},\\
\gax &:= \text{ the subgroup of } \Sim_X \text{ generated by } \phi(X).
\end{align*} \end{defn}

By \eqref{ccc4}, $\gax$ is a subgroup of $\gx$. On the other hand,
it is easy to see that
\begin{equation}\label{conjugada}
g\phi_x g^{-1} = \phi_{g(x)}, \qquad \forall g\in \gx, x\in X.
\end{equation}

Therefore, $\phi(X) \subset \gx$ is a standard crossed subset,
$\phi:X\to \gx$ is a morphism of racks, and $\gax$ is a normal
subgroup of $\gx$. It is not true in general that $\gax = \gx$.

\begin{exmp} Let $X = \{\pm i, \pm j\}$, a standard subset of the group
of units of the quaternions. Then $\gax \neq \gx$. \end{exmp}
\pf
It is easy to compute $\gx$ and $\gax$. One sees that
$\gax=\;<\phi_i,\,\phi_j>\,\simeq C_2\times C_2$ and
$\gx=\;<\phi_i,\,\phi_j,\,\sigma_0>$ has order $8$, where
$\sigma_0(\pm i)=\pm j$ and $\sigma_0(\pm j)=\pm i$.
\epf

Another basic group attached to $X$ is the following one:
\begin{defn}\label{gx}\cite{jo,Bk,fr}
Let $(X,\trid)$ be a rack. We define the \emph{enveloping group} of $X$ as
$$G_X = F(X) / \langle  xyx^{-1}=x\trid y, \quad x,y\in X\rangle,$$
where $F(X)$ denotes the free group generated by $X$. The assignment
$x\mapsto \phi_x$ extends to a group homomorphism $\pi_X: G_X \to \gax$;
the kernel of $\pi_X$ is called the \emph{defect group} of $X$ in the
literature and coincides with the subgroup $\Gamma$ considered in \cite[Th. 2.6]{s}.
\end{defn}
The name ``enveloping group" is justified by the following fact, contained essentially
in \cite{jo}:
\begin{lem}\label{lm:neej}%nombre enveloping esta justificado
The functor $X\mapsto G_X$ is left adjoint to the forgetful functor $H\mapsto \forget H$
from the category of groups to that of racks. That is,
$$\Hom_{\mbox{groups}}(G_X,H) \simeq \Hom_{\mbox{racks}}(X,\forget H)$$
by natural isomorphisms.\hfill\qed
\end{lem}

The definition of rack was proposed a long time ago by Conway and Wraith,
see the historical account in \cite{fr}.
Quandles were introduced independently in \cite{jo} and \cite{ma}
and studied later in \cite{Bk} and other articles.
They are being extensively studied nowadays in relation with knot
invariants, see \cite{cs} and references therein.
In \cite{gr}, it was proposed to consider crossed sets with the
normalizing conditions \eqref{ccc2} and \eqref{ccc3}; the conditions also appear in
\cite{s}. It is worth noting that in most of these articles the quandle
structure is the opposite to the one here (i.e., $x*y$ for our $\phi_y^{-1}(x)$).
Racks, quandles and crossed sets are related as follows.

\subsubsection*{\bf From racks to quandles.}

We follow  \cite{Bk}. Let $X$ be a rack and let $C_{\trid}(X)$ be
the centralizer in $\gx$ of $\gax$. For $\psi\in C_{\trid}(X)$,
define $\trid^{\psi}$ by
\begin{equation}\label{eq:frtq}
a\trid^{\psi} b = a\trid \psi (b) =
\psi (a\trid b) = \psi (a) \trid \psi (b).
\end{equation}
Then $(X, \trid^{\psi})$ is again a rack; we say that it is
\emph{conjugated to $(X, \trid)$ via $\psi$}.

Let now $\iota: X \to X$ be given by $a\trid \iota(a) =
a$, which is well defined by \eqref{ccc1}. Then, by \eqref{ccc4},
$$ a\trid b = a\trid (b \trid \iota(b)) = (a\trid b) \trid (a
\trid \iota(b)); $$ hence $a\trid \iota (b) = \iota (a\trid b)$.
In particular, $a= \iota (a\trid a)$, and $\iota$ is surjective.
Also, $$ a\trid (\iota(a) \trid \iota(b)) = (a\trid \iota(a))
\trid (a \trid \iota(b)) = a \trid (a \trid \iota(b)), $$ so that,
by \eqref{ccc1}, $\iota(a) \trid \iota(b) = a \trid \iota(b)=
\iota (a\trid b)$. Suppose now that $\iota(a) = \iota(b)$. Then $$
a = a\trid \iota(a) = \iota(a)\trid \iota(a) = \iota(b)\trid
\iota(b) = b; $$ That is, $\iota$ is injective, and belongs to
$C_{\trid}(X)$. We can then consider $(X, \trid^{\iota})$, which
is a quandle.

In conclusion, any rack $(X, \trid )$ is conjugated to a unique quandle,
called \emph{the quandle associated to $(X, \trid )$}. If $F: X\to Y$ is a
morphism of racks, then $F\iota = \iota F$; hence $F$ is a morphism of
the associated quandles. It follows that $\gx \simeq \gxiota$.
(But $\gax$ and $\gaxiota$ may be very different).

\subsubsection*{\bf From quandles to  crossed sets.}

We exhibit a functor $Q$ from the category of finite quandles to
that of crossed sets, which assigns to a quandle $X$ a quotient
crossed set $Q(X)$ with the expected universal property: any
morphism of quandles $X\to Y$ is uniquely factorized through $
Q(X)$ whenever $Y$ is a crossed set. To see this, take on $X$
$\sim$ as the equivalence relation generated by
\begin{equation}\label{eq:qacs}%quandles a crossed sets
x\sim x'\mbox{ if }\exists y\mbox{ s.t. }x\trid y=y\mbox{ and
}y\trid x=x'.
\end{equation}
Then $\sim$ coincides with the identity relation if and only if
$X$ is a crossed set. Take $X_1:=X/\sim$. We must see that $X_1$
inherits the structure of a quandle. First, suppose $x,x',y$ are
as in \eqref{eq:qacs}. Then $x'\trid y=(y\trid x)\trid (y\trid
y)=y\trid(x\trid y)=y$. Next, for $z\in X$, we have
$$y\trid(x'\trid z)=(x'\trid y)\trid(x'\trid z)=x'\trid(y\trid z)
    =(y\trid x)\trid(y\trid z)=y\trid(x\trid z),$$
whence we see that $\phi_x=\phi_{x'}$, and then
$\phi_x=\phi_{x''}$ for any $x''\sim x$. Thus, it makes sense to
consider $\trid:(X/\sim)\times X\to X$. Finally, we have for $z\in
X$,
\begin{align*}
&(z\trid x)\trid(z\trid y)=z\trid(x\trid y)=(z\trid y)\mbox{ and}
\\ &(z\trid y)\trid(z\trid x)=z\trid(y\trid x)=(z\trid x').
\end{align*}
Then $(z\trid x)\sim(z\trid x')$, and it makes sense to consider
$\trid:(X/\sim)\times(X/\sim)\to(X/\sim)$. If $X_1$ is not a crossed set
then take $X_2=X_1/\sim$, and so on. Since $X$ is finite, we must
eventually arrive to a crossed set. The functoriality and
universal property are clear.

%%%%%%%%%%%%%%%%%%%%%%%%%%%%%%%%%%%%%%
%%%%%% SubSeccion  Basic Def. %%%%%%%%
%%%%%%%%%%%%%%%%%%%%%%%%%%%%%%%%%%%%%%
\subsection{Basic definitions}

We collect now a number of definitions and results; many of them
appear already in previous papers on racks or quandles, see
\cite{jo,Bk,quandle1,fr}.

We shall say that $X$ is \emph{trivial} if $i\trid j=j$ for all $i, j\in X$.

\begin{lem}\label{primerlema}
Let $(X,\trid)$ be a rack, $H$ a group and $\varphi: X \to H$ an injective
morphism of racks such that the image is invariant under conjugation in $H$
(thus $(X,\trid)$ is actually a crossed set).
Then the map $\widetilde{\varphi}: H \to \Sim_X$, given by
$\widetilde{\varphi}_h(x) = \varphi^{-1}(h\varphi(x)h^{-1})$, is a
group homomorphism and its image is contained in $\gx$.
\end{lem}
\pf Left to the reader.
\epf

The assignment $X\mapsto\gax$ is not functorial in general; just
take $i\in X$ with $\phi_i \neq \id$; the inclusion $\{i\} \subset
X$ does not extend to a morphism $\Inn_{\trid} (\{i\}) \to \gax$
commuting with $\phi$. But we have:

\begin{lem}\label{le:gammafuntor}
If $\pi: X \to Y$ is a surjective morphism of racks, then it
extends to a group homomorphism $\Inn_{\trid}  (\pi): \gax \to
\gay$.
\end{lem}

\pf Let $\Inn_{\trid}  (\pi)$ be defined on $\phi(X)$ by
$\Inn_{\trid}  (\pi)(\phi(x))=\phi(\pi(x))$. It is well defined
and it extends to a morphism of groups since $\pi$ is surjective.
\epf

We determine now the structure of $\gax$ when $X$ is standard.

\begin{lem}\label{le:gammastandard}
Let $H$ be a group and let $X\subset H$ be a standard subset.
\begin{enumerate}
\item\label{gms1}
$\gax\simeq C/Z(C)$, where $C=\langle X\rangle$ is the
	subgroup of $H$ generated by $X$, which is clearly normal.
\item\label{gms2}
If $X$ generates $H$ then $\gax\simeq H/Z(H)$.
\item\label{gms3}
If $H$ is simple non-abelian and $X\neq\{1\}$, then $H=\gax$.
\end{enumerate}
\end{lem}

\pf We prove (\ref{gms2}); (\ref{gms1}) and (\ref{gms3}) will follow.
By Lemma \ref{primerlema}, we have a morphism\linebreak
$\psi:H\to\gx$, whose image is $\gax$. Now,
$$h\in\ker(\psi)\iff hxh^{-1}=x\ \forall x\in X\iff h\in Z(H).$$
\epf

\begin{cor}\label{co:gammastandard}
Let $X$ be a rack, let $H$ be a group and let $\psi: X\to H$ be a morphism.
If $Z(\langle \psi(X) \rangle)$ is trivial, then $\psi$ extends to a morphism
$\Psi: \gax \to H$.
\end{cor}
\pf
If $Y = \psi(X)$, then $\Psi: \gax \fllad{\Inn_{\trid}  (\psi)}
\gay \simeq \langle Y \rangle \hookrightarrow H $.
\epf

The map $\phi:X\to\gax$ is not injective, in general.

\begin{defn}We shall say that the rack $(X, \trid)$ is \emph{faithful} when
the corresponding $\phi$ is injective. Observe that in this case $X$ is
a crossed set, since it is standard.
\end{defn}

\begin{rem}\label{rm:centrotrivial}
If $(X, \trid)$ is faithful then the center of $\gax$ is trivial.
More generally, if $z\in Z(\gax)$, then $\phi_{z(i)}=\phi_i$, for all $i\in X$.
\end{rem}

\begin{defn} A \emph{decomposition} of a rack $(X,\trid)$ is a disjoint union 
$X = Y \cup Z$ such
that $Y$ and $Z$  are both  subracks of $X$.
(In particular, both $Y$ and $Z$ are non-empty).
$X$ is \emph{decomposable} if it admits a decomposition,
and \emph{indecomposable} otherwise.
\end{defn}

The image of an indecomposable rack under a morphism is
again indecomposable. 

We shall occasionally denote $i^n \trid j :=
\phi_i^n(j)$, $n\in \Z$.
The \emph{orbit} of an element $x\in X$ is the subset
$$\orb_x := \{i_s^{\pm 1}\trid(i_{s-1}^{\pm 1} \trid( \dots
	(i_1^{\pm 1} \trid x) \dots ))\ |\ i_1, \dots, i_s\in X\}.$$
That is, $\orb_x$ is the orbit of $x$ under the natural action of
the group $\gax$. (If $X$ is finite, then
$\orb_x=\{i_s\trid(i_{s-1} \trid(\dots (i_1 \trid x) \dots ))
	\ |\  i_1, \dots, i_s \in X  \}$).

\begin{lem}\label{indecomposable} 
Let $(X,\trid)$ be a rack, $Y \neq X$ a non-empty subset and $Z=X-Y$.
Then the following are equivalent:
\begin{enumerate}
\item $X = Y \cup Z$ is a decomposition of $X$.
\item $Y\trid Z \subseteq Z$ and $Z\trid Y \subseteq Y$.
\item $X\trid Y \subseteq Y$. \hfill\qed
\end{enumerate}
\end{lem}

\begin{lem} Let $(X,\trid)$ be a rack. Then the following are equivalent:
\begin{enumerate}
\item $X$ is indecomposable.
\item $X = \orb_x$ for all (for some) $x \in X$.
\hfill\qed
\end{enumerate}
\end{lem}

Note that a standard crossed set $X\subset H$ need not be
indecomposable, even if it consists of only one $H$-orbit.
However, it is so when $H$ is simple by Lemma
\ref{le:gammastandard}.

\begin{exmp}
If $X\subset H$ is a conjugacy class with two elements, then it is
trivial as crossed set. As another example, take $A$ an abelian group
and $G$ the group of automorphisms of $A$; let $H=A\rtimes G$, and let
$X\subset A$ be any orbit for the action of $G$. Let $1_G\in G$ be the
unit and consider $X\subset H$ as $X\rtimes 1_G$. Then $X$ is trivial.
\end{exmp}

\begin{prop}
Any rack $X$ is the disjoint union of maximal indecomposable subracks.
\end{prop}
\pf
Given $Y\subset X$ a subset, consider
$$Y'=Y\cup (Y\trid Y)\cup (Y^{-1}\trid Y)
	=Y\cup \{y\trid z\ |\ y,z\in Y\}\cup\{y^{-1}\trid z\ |\ y,z\in Y\}.$$
Then $Y'\supset Y$ and any subrack of $X$ containing $Y$ contains $Y'$.
The \emph{subrack generated} by $Y$ is thus $\bigcup_{n\in\N}Y^n$, where
$Y^{n+1}=(Y^n)'$ and $Y^1=Y$. This is the smallest subrack of $X$ containing $Y$.

For $Y\subset X$, we say that it is \emph{connectable} if
for any two elements $y_1,y_2\in Y$ there exist
$u_1^{\ep_1},\ldots,u_n^{\ep_n}$, where $u_i\in Y$ and
$\ep_i\in\{\pm 1\}$ $\forall i$, such that
$y_2=u_1^{\ep_1}\trid(u_2^{\ep_2}\trid\cdots (u_n^{\ep_n}\trid y_1))$
(here the intermediate elemments
$u_i^{\ep_i}\trid(u_{i+1}^{\ep_{i+1}}\trid\cdots (u_n^{\ep_n}\trid y_1))$
may not belong to $Y$).
Then it is easy to see that if $Y$ is connectable then so is $Y'$.
Hence, for $Y$ connectable, the
subrack generated by $Y$ is connectable and, being a subrack, it
is indecomposable. Also, since the union of intersecting indecomposable
subracks is connectable, we see that they generate an indecomposable
subrack. Hence the \emph{indecomposable component} of $x\in X$, the
union of all indecomposable subracks containing $x$, is an indecomposable
subrack. Now, $X$ is the disjoint union of such components.
%Given $x\in X$, $\{x^n\trid x,\ n\in\Z\}$ is an 
%indecomposable subrack containing $x$. Indeed,
%$$(x^n\trid x) \trid(x^m\trid y)= (\phi_x^n x) \trid(\phi_x^m y)
%= (\phi_x^n( x \trid(\phi_x^{m - n} x))) = \phi_x^{m + 1} x.$$
%Define  the ``indecomposable component" of $x$ by
%$$\bigcup_{\substack{x\in Y\subset X \\ Y \text{ indecomposable}}} Y.$$
%Then $X$ is the disjoint union of the indecomposable components
%of all elements, which are maximal indecomposable subracks.
\epf

Unlike the situation of Lemma \ref{indecomposable}, the indecomposable components
may not be stable under the action of $X$. The case of two
components is more satisfactory because  we can describe how to glue two racks.

\begin{lem}\label{segundolema}\
\begin{enu}
\item Let $Y,Z$ be two racks and $X=Y\sqcup Z$ be their disjoint union.
The following are equivalent:
\begin{enu2}
\item Structures of rack on $X$ such that $X = Y \cup Z$ is a
decomposition.

\item Pairs  $(\sigma, \tau)$ of morphisms of racks
$\sigma: Y \to \Aut_{\trid}(Z)$, $\tau: Z \to \Aut_{\trid}(Y)$
such that
\begin{flalign}
\label{decomp2} & y\trid\tau_{z}(u)=\tau_{\sigma_{y}(z)}(y\trid u),\quad
	\forall y,u\in Y,\ z\in Z,&&
	\text{i.e., } \phi_y\tau_z=\tau_{\sigma_y(z)}\phi_y; \\
\label{decomp3} & z\trid\sigma_{y}(w)=\sigma_{\tau_{z}(y)}(z\trid w),\quad
	\forall y\in Y,\ z,w \in Z,&&
	\text{i.e., } \phi_z\sigma_y=\sigma_{\tau_z(y)}\phi_z.
\end{flalign}
\end{enu2}
\item Assume that $Y$ and $Z$ are crossed sets and \eqref{decomp2}, 
\eqref{decomp3} hold.
Then $X$ is a crossed set exactly when 
\begin{flalign}\label{decomp1}
&\sigma_{y}(z)=z \text{ if and only if } \tau_{z}(y)=y, \quad
	\forall y\in Y, \quad z \in Z.&
\end{flalign}
\end{enu}
\end{lem}

\pf Left to the reader. \epf

If the conditions of the Lemma are satisfied, we shall say that $X$ is the
\emph{amalgamated sum} of $Y$ and $Z$. If $\sigma$ and $\tau$ are trivial,
we say that $X$ is the \emph{disjoint sum} of $Y$ and $Z$. Clearly, one can
define the disjoint sum of any family of racks (resp. quandles, crossed sets).

For example, let $X$ be a rack (resp. quandle, crossed set)
and set $X \times 2=2X=X\times \{1, 2\}$; this is a rack (resp. quandle,
crossed set) with $(x, i) \trid (y, j) = (x\trid y, j)$, and
$2X = X_1 \cup X_2$ is a decomposition, where $X_i = X\times \{ i \}$.
Note that $\phi_{(x, i)} =  \phi_{(x, j)}$; $\phi$ is not injective.
In an analogous way, we define the crossed set $nX$, for any positive
integer $n$. More generally, we have

\begin{exmp}\label{prodcs}
Let $X$, $Y$ be two racks (resp. quandles, crossed sets). Then $X\times Y$ is a
rack (resp. quandle, crossed set), with $(x,y)\trid(u,v)=(x\trid u,y\trid v)$;
this is the direct product of $X$ and $Y$ in the category of racks
(resp. quandles, crossed sets).
\end{exmp}

\begin{lem} Let $X$, $Y$ be two racks (resp. quandles, crossed sets).
If $X\times Y$ is indecomposable then $X$ and $Y$ are indecomposable.
The converse is true if $X$ or $Y$ is a quandle.
\end{lem}
\pf
Since the canonical projections $X\times Y\to X$ and
$X\times Y\to Y$ are rack homomorphisms, the first statement is immediate.
For the second one, let $(x_1,y_1),(x_2,y_2)\in X\times Y$. As $X,Y$ are
indecomposable, there exist $u_1^{\ep_1},\ldots,u_n^{\ep_n}$,
$v_1^{\ep'_1},\ldots,v_m^{\ep'_m}$ where $u_i\in X$, $v_j\in Y$
and $\ep_i,\ep'_j\in\{\pm 1\}$ $\forall i,j$, such that
$x_2=u_1^{\ep_1}\trid(u_2^{\ep_2}\trid\cdots (u_n^{\ep_n}\trid x_1))$ and
$y_2=v_1^{\ep'_1}\trid(v_2^{\ep'_2}\trid\cdots (v_m^{\ep'_m}\trid y_1))$.
Suppose $X$ is a quandle. Adding if necessary at the end of the sequence
in $Y$ pairs $y_1,y_1^{-1}$, we may suppose that $m\ge n$. Adding if necessary
at the end of the sequence in $X$ elements $x_1$, we may suppose that
$m=n$. Then,
$(x_2,y_2)=(u_1^{\ep_1},v_1^{\ep'_1})\trid((u_2^{\ep_2},v_2^{\ep'_2})
	\trid\cdots((u_n^{\ep_n},v_n^{\ep'_n})\trid (x_1,y_1)))$.
%suppose that $X$ is a quandle and let $Z\subseteq X\times Y$
%be the orbit of a point $(x,y)$. It is clear that $\pi_X(Z)=X$, since $X$
%is indecomposable. Now, as $X$ is a quandle and $Y$ is indecomposable,
%if $(x',y')\in Z$ then $x'\times Y\subset Z$, and we are done.
\epf

Let $X$ be a rack, take $\phi:X\to\gax$ as usual. Let us abbreviate
$F_y := \phi^{-1}(y)$, a fiber of $\phi$. If $X$ is a quandle, any
fiber $F_y$ is a trivial subquandle of $X$.

\begin{lem}
\begin{enumerate}
\item\label{yuti1}
For $x, y\in \phi(X)$ the fibers $F_y$ and $F_{x\trid y}$  have the
same cardinality.
\item\label{yuti2}
If $X$ is an indecomposable crossed set, then the fibers of
$\phi$ all have the same cardinality.
\end{enumerate}
\end{lem}

\pf
We claim that $i\trid F_y \subseteq F_{\phi_i\trid y}$.
Indeed, if $j\in F_y$, then
$\phi_{i\trid j}=\phi_i\phi_j\phi_i^{-1}=\phi_iy\phi_i^{-1}=\phi_i\trid y $;
the claim follows. Similarly, $i^{-1}\trid F_y \subseteq
F_{\phi_i^{-1}\trid y}$, hence (\ref{yuti1}).
Now (\ref{yuti2}) follows from (\ref{yuti1}).
\epf

We give finally some definitions of special classes of crossed
sets, following \cite{jo}.

\begin{defn} Let $(X, \trid)$ be a quandle.  We shall say that $X$
is \emph{involutory} if $\phi_x^2=\id$ for all $x\in X$. That is,
if $x\trid (x\trid y) = y$ for all $x, y\in X$.

We shall say that $X$ is \emph{abelian} if $(x\trid w) \trid
(y\trid z) = (x\trid y) \trid (w\trid z)$ for all $x, y, w,z\in
X$.
\end{defn}

%%%%%%%%%%%%%%%%%%%%%%%%%%%%%%%%%%%%%%
%%%%%% SubSeccion  Examples   %%%%%%%%
%%%%%%%%%%%%%%%%%%%%%%%%%%%%%%%%%%%%%%
\subsection{Examples}\label{ss:examples}

\subsubsection*{\bf A rack which is not a quandle.}
Take $X$ any set and $f\in\Sim_X$ any function. Let $x\trid y=f(y)$.
This is a rack, and it is not a quandle if $f\neq\id_X$. This rack is
called \emph{permutation rack}.
\subsubsection*{\bf A quandle which is not a crossed set.}
Take $X=\{x,+,-\}$, $x\trid\pm=\mp$, $\phi_{\pm}=\id_X$.
\subsubsection*{\bf Amalgamated sums.}
Let $Z$ be a rack; we describe all the amalgamated sums $X=Y\cup Z$
for $Y=\{0,1\}$ the trivial rack. Denote $\Aut(Y)=\{+=\id,-\}$.
Let $\sigma$, $\tau$ be as in Lemma \ref{segundolema}. First, $Z$
should decompose as a disjoint union of subracks $Z=Z_+\cup Z_-$,
where $\tau(Z_{\pm})=\pm$.
%Second, \eqref{decomp1} is equivalent to $\sigma_0$ and $\sigma_1$
%being the identity on $Z_+$ and having no fixed points on $Z_-$.
Second, \eqref{decomp2} is equivalent to $Z_{\pm}$ being stable by
$\sigma_0$ and $\sigma_1$; and condition \eqref{decomp3} reads
\begin{align*}
&\phi_z\sigma_0=\sigma_0\phi_z,\quad\phi_z\sigma_1=\sigma_1\phi_z,
	\quad\forall z\in Z_+, \\
&\phi_z\sigma_0=\sigma_1\phi_z,\quad\phi_z\sigma_1=\sigma_0\phi_z,
	\quad\forall z\in Z_-.
\end{align*}
Another way to describe the situation is: $Z=Z_+\cup Z_-$ a disjoint union;
let ${C_+=<\phi_x,\ x\in Z_+>}$ be the group generated by $\phi_{Z_+}$,
$C_-=<\phi_x\phi_y,\ x,y\in Z_->$, and let $C=<C_+,C_->$.
Then $[\sigma_0,C]=1\in\Sim_Z$, $\sigma_1=\phi_x\sigma_0\phi_x^{-1}$ for any
$x\in Z_-$.
%If $\sigma_1\neq\sigma_0$ and $Z_{\pm}$ both are not empty,
%then it can be seen that $Z_+\cup Z_-$ is a decomposition of $Z$.
If $Z$ is a quandle then $X$ is a quandle. If $Z$ is a crossed set then $X$ is a
crossed set iff $Z_+=Z^{\sigma_0}=\{\text{fixed points of }\sigma_0\}=Z^{\sigma_1}$.

\subsubsection*{\bf Polyhedral crossed sets.}

\begin{figure}[ht]
\begin{center}
\epsfig{file=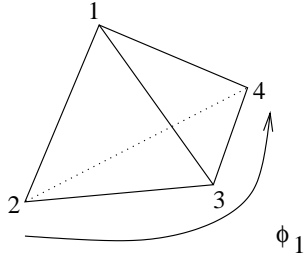,width=4cm}
\caption{Polyhedral crossed set of the tetrahedron}\label{fg:tetra}
\end{center}
\end{figure}

Let $P\subset \R^3$ be a regular polyhedron with vertices
$X=\{x_1,\ldots,x_n\}$ and center in $0$. For $1\le i\le n$,
let $T_i$ be the orthogonal linear map which fixes $x_i$ and rotates
the orthogonal plane by an angle of $2\pi/r$ with the right hand rule
(pointing the thumb to $x_i$), where $r$ is the number of edges ending in
each vertex. Then $(X,\trid)$ defined by $x_i\trid x_j=T_i(x_j)$ is a
crossed set. To see this, simply take $\g$ as the group of orthogonal
transformations of $P$ with determinant $1$ and notice that
$\{T_1,\ldots,T_n\}$ is a conjugacy class of $\g$, whose
underlying (standard) crossed set is isomorphic to $X$.
It is evident that $\gax$ is isomorphic to the group generated
by $\{T_1,\ldots,T_n\}$. It is clear that to each polyhedron
also corresponds an analogous crossed set given by the faces;
it is isomorphic to the crossed set given by the vertices
of the dual polyhedron. It follows by inspection that
the crossed sets of the vertices of the tetrahedron,
the octahedron, the dodecahedron and the icosahedron are
indecomposable, while that of the cube has two components,
each of which is isomorphic to the crossed set of the vertices
of the tetrahedron.
%Just because it is closer to the authors
%intuition, we will usually work with the faces crossed sets,
%except for the tetrahedron.
It is easy to see that in the indecomposable
cases the group $\gax$ coincides with $\g$. See Figure
\ref{fg:tetra}.

\subsubsection*{\bf Coxeter racks.}
Let $(V,\,<\,,>)$ be a vector space over a field $k$, provided with
an anisotropic symmetric bilinear form $<\,,>$.
Let $v\trid u=u-\frac{2 <u,v>}{<v,v>}v$.
Then $(V-\{0\}, \trid)$ is a rack, as well as any
subset closed for this operation. Particular cases of this are the
root systems of semisimple Lie algebras; the action
$\phi_{\alpha}$ coincides with the action of $w_{\alpha}\in W$,
the Weyl group. To turn this into a quandle one can either quotient
out by the relation $v\sim -v$, or take the conjugated quandle
$(V-\{0\}, \trid^{\iota})$ as at the end of subsection \ref{ss:rqcs}
(see \eqref{eq:frtq}).
It is easy to see that $v\trid^{\iota} u=\frac{2 <u,v>}{<v,v>}v-u$,
and then this also is a crossed set.

\subsubsection*{\bf Homogeneous crossed sets.}
Let $G$ be a finite group and let $\phi:\Sim_G\to\fun(G,\Sim_G)$
be the function given by
$$ \phi^s_x(y) = s(yx^{-1})x, \qquad s\in \Sim_G, x, y\in G.$$
Then $\phi^s_x\phi^t_x = \phi^{st}_x$, i.e., $\phi$ is a morphism
of groups, with the pointwise multiplication in\linebreak
$\fun(G,\Sim_G)$.

If, in addition, $s:G\to G$ is a group automorphism, then define
$x\trid y=\phi^s_x(y) = s(yx^{-1})x$. It is easy to see that this
makes $(G,\trid)$ into a crossed set. For instance, let us check
\eqref{ccc3}: $$ x\trid y=y \iff y=s(yx^{-1})x \iff yx^{-1} =
s(yx^{-1}) \iff xy^{-1} = s(xy^{-1}) \iff y\trid x=x. $$ We shall
say that $(G,\trid)$ is a \emph{principal homogeneous} crossed
set, and we will denote it by $(G, s)$.

Let $t: G \to G$, $t(x) = s(x^{-1})x$, so that $x\trid y=  s(y)\,
t(x)$.  It is clear that
\begin{equation}\label{fibras-afine}
\phi_x =\phi_z \iff t(x) = t(z),
\end{equation}
whence the fibers as a crossed set are the same as the fibers of
$t$. Note that $t$ is a group homomorphism if and only if
$\im (t) \subseteq Z(G)$, the center of $G$.

More generally, let $H \subset G^s$ be a subgroup, where
$G^s$ is the subgroup of elements of $G$ fixed by $s$. Then
$H\backslash G$ is a crossed set, with $Hx\trid Hy =
Hs(yx^{-1})x$; it is called a \emph{homogeneous} crossed set.

It can be shown that a crossed set $X$ is homogeneous if
and only if it is a single orbit under the action of $\gx$
\cite{jo}.

\subsubsection*{\bf Twisted homogeneous crossed sets.}
In the same vein, let $G$ be a group and $s\in\Aut(G)$. Take $x\trid y=xs(yx^{-1})$.
This is a crossed set which is different, in general, to the previous one.
Any orbit of this is called a \emph{twisted homogeneous crossed set}.

\subsubsection*{\bf Affine crossed sets.}
Let $A$ be a finite abelian group and $g:A\to A$ an isomorphism of
groups. The corresponding principal homogeneous crossed set is
called an \emph{affine} crossed set, and denoted by $(A, g)$.
(They are also called Alexander quandles, see e.g.  \cite{quandle1}).

Let $f=\id-g$, then $x\trid y =x+g(y-x) =f(x)+g(y)$.

Let us compute the orbits of $(A, g)$. Since $g=\id-f$, we have
$x\trid y=f(x-y)+y$. It is clear then by induction that
$x_1\trid(x_2\trid\cdots(x_n\trid y))\in y+f(A)$, whence $\mathcal
O_y=y+\im(f)$. In this case, by \eqref{fibras-afine},
indecomposable is equivalent to being faithful (since $A$ is finite).

We  compute now when two indecomposable affine crossed
sets are isomorphic.

\begin{lem}\label{lm:autafin}
Two indecomposable affine crossed sets $(A, g)$ and $(B, h)$ are
isomorphic if and only if there exists a linear isomorphism $T:
A\to B$ such that $Tg = hT$.

If this happens, any isomorphism of crossed sets $U :(A, g) \to
(B, h)$ can be uniquely written as $U = \tau_bT$, where $\tau_b:
B\to B$ is the translation $\tau_b(x) = x+b$ and $T: A\to B$ is a
linear isomorphism such that $Tg = hT$. \end{lem}

\pf Let $U :(A, g) \to (B, h)$ be an isomorphism of crossed sets.
Since the translations are isomorphisms of crossed sets, we
decompose $U = \tau_bT$, where $T$ is an isomorphism of crossed
sets with $T(0) = 0$. Let $f = \id -g$, $k = \id -h$. We have
$$T(f(x) + g(y)) = k(T(x)) + h(T(y)).$$ Letting $x = 0$, we see
that $Tg = hT$; letting $y = 0$, we see that $Tf = kT$. Hence
$T(f(x) + g(y)) = T(f(x)) + T(g(y))$. But $(A, g)$ is
indecomposable, thus $f$ is bijective; we conclude that $T$ is
linear.
\epf

\begin{rem}
After we finished (the first version of) the paper, Nelson gave in
\cite{ne} a classification of non indecomposable affine crossed sets.
\end{rem}

\begin{cor}
If $(A, g)$ is an indecomposable affine crossed set, then
$\Aut_{\trid}(A, g)$ is the semidirect product $A \rtimes
\Aut(A)^g$, where $\Aut(A)^g$ is the subgroup of all linear
automorphisms of $A$ such that $Tg=gT$. \hfill\qed
\end{cor}

Notice that when $(A,g)$ is not indecomposable, the corollary does
not hold. For instance, if $g=\id$, then $f=0$ and $(A,\id)$ is
trivial, whence $\Aut_{\trid}(A,g)=\Sim_A$.

The group $\Inn_{\trid}  (A, g)$ is usually smaller: it
can be easily shown that $\Inn_{\trid}  (A, g) =\im f\rtimes
\langle g\rangle$. In fact, if $a\in A$ then
$\phi_a=(f(a),g)\in A\rtimes\Aut(A)^g$. In particular,
$\Inn_\trid(A,g)$ is solvable.

\begin{rem}\label{rm:wia}
We conclude that a standard crossed set $\mathcal O$, where
$\mathcal O$ is a single non-trivial orbit of a simple non-abelian
group, can not be affine; \emph{cf.} Lemma \ref{le:gammastandard}.
\end{rem}

Similarly, let us discuss when a polyhedral crossed set $X$ is
affine. The crossed set of vertices of the tetrahedron is affine,
actually isomorphic to $(A = \Z_2 \times \Z_2, g= \mdpd 0111 )$.
It can be seen by hand that the crossed set of the vertices of the
cube is not affine. Indeed, it is easy to see that the underlying group
$A$ should be either $\Z/2\times\Z/4$ or $(\Z/2)^3$, but since
$\Aut(\Z/2\times\Z/4)$ has order $8$ and $\Inn_\trid(X)$ has $3$-torsion,
the case $A=\Z/2\times Z/4$ is impossible.
Furthermore, one can exclude the case $A=(\Z/2)^3$ by looking at automorphisms
$g\in\Aut(A)$ s.t. $g^3=\id$. For the other polyhedral crossed sets, we have that
$\Inn_{\trid}(X)$ is $\Sim_4$ for the octahedron, and $\mathbb A_5$
for the icosahedron and the dodecahedron, so that $X$ is not affine.

As a particular case, let $A=\Z_n=\Z/n\Z$, and let
$q\in\Z_n^{\times}$ such that $(1-q)\in\Z_n^{\times}$. Then we
have a structure on $\Z_n$ given by $$(\Z_n,\trid^q)\qquad
x\trid^qy=(1-q)x+qy.$$ This is an indecomposable crossed set, and
it can be seen that $(\Z_n,\trid^q)\simeq(\Z_n,\trid^{q'})\iff
q=q'$ (by Lemma \ref{lm:autafin}, since $\Aut(\Z_n)$ is abelian).

It is immediate to see that $(\Z_3,\trid^2)$ is the only
indecomposable crossed set with $3$ elements; and that any crossed
set with $3$ elements is either trivial or isomorphic to
$(\Z_3,\trid^2)$.

It is proved in \cite{egs} that any indecomposable rack of prime order $p$
is either isomorphic to $(\Z_p,\trid^q)$ for $q\in\Z_p^{\times}$ or it is
isomorphic to $(\Z_p,\trid)$ with $x\trid y=y+1$.
Indecomposable racks of order $p^2$ are classified in \cite{grp2}, in particular
it is proved there that any indecomposable quandle of order $p^2$ is affine.

\subsubsection*{\bf Affine racks.}
These are a generalization of affine quandles. Let $A$ be an abelian group,
$g \in \Aut(A)$, $f \in \End(A)$ be such that $fg=gf$ and $f(\id-g-f)=0$.
We define then on $A$ the structure of rack given by $x \trid y = f(x) + g(y)$.
It is clear that $A$ is a quandle iff $f=\id-g$. Notice that
$(\id-f)(\id+g^{-1}f)=\id$; thus $\id-f$ is an automorphism of $A$. We can
consider then the affine quandle $(A,\id-f)$. It is easy to see
that this quandle is the associated quandle for the rack just defined
(see \eqref{eq:frtq}). 
As an example, one can take $A=\Z_{p^2}$, $g=\id$, $f(x)=px$.

\subsubsection*{\bf Amalgamated unions of affine crossed sets.}
Let $(A, g)$, $(A,h)$ be two indecomposable affine crossed sets;
let $f,k$ be given by $f=\id-g$, $k=\id-h$, and let
$\sigma:(A,g)\to\gxd{A,h}$,\linebreak
$\tau:(A,h)\to\gxd{A,g}$ have the form $\sgm_x=\al(x)+\be$
(i.e., $\sgm_x(y)=\al(x)+\be(y)$) and $\tau_y=\gm+\dt(y)$
(i.e., $\tau_y(x)=\gm(x)+\dt(y)$) for certain $\be,\gm\in\Aut(A)$,
$\al,\dt\in\End(A)$. Then by Lemma \ref{lm:autafin} we must
have $\be\in(\Aut(A))^h$, $\gm\in(\Aut(A))^g$ and one can verify that
%$\sigma_x=\alpha(x)+\beta$;
%$\tau_y=\gamma+\delta(y)$. Then $\beta\in (\Aut A)^h$, $\gamma\in
%(\Aut A)^g$, and we have:
\begin{flalign*}
&\mbox{\eqref{decomp2}} && \mbox{is equivalent to }&
    & f-f\gamma-\delta\alpha=0,\qquad \delta\beta-g\delta=0;& \\
&\mbox{\eqref{decomp3}} && \mbox{is equivalent to }&
    & k-k\beta-\alpha\delta=0,\qquad \alpha\gamma-h\alpha=0;& \\
&\mbox{\eqref{decomp1}} && \mbox{is equivalent to }&
    & \alpha(x)+\beta(y)=y\iff \gamma(x)+\delta(y)=x.&
\end{flalign*}

Thus, the disjoint union $X = (A, g) \sqcup (A, h)$ is a
decomposable crossed set. For instance, let $h = -gf^{-1} = \id -
f^{-1}$. We define $\sigma$ and $\tau$ by taking $\alpha=\id$,
$\beta=g$, $\gamma=h$, $\delta=\id$.

If $(A,g)$, $(B,h)$ are non isomorphic indecomposable affine crossed
sets, then any amalgamated sum $A\cup B$ is not affine. This is a
consequence of the following easy lemma:

\begin{lem}
If $(A,g)$ is an affine crossed set, then its orbits are
isomorphic as crossed sets.
\end{lem}
\pf Let $\{x_0=0,\ldots,x_n\}$ be a full set of representatives of
coclasses in $A/\im(f)$. The orbits are $\im(f),\ \im(f)+x_1,\ldots,\
\im(f)+x_n$. Thus, $\tau_i:\im(f)+x_i\to\im(f)$, $\tau_i(y)=y-x_i$
are isomorphisms.
\epf

\subsubsection*{\bf Involutory crossed sets}
This is a discretization of symmetric spaces, which we shall
roughly present (see \cite{jo} for the full explanation). Let $S$
be a set provided with a collection of
functions $\gamma:\Z\to S$, called \emph{geodesics}, such that any
two points of $S$ belong to the image of some of them. 
Consider the affine crossed
set $(\Z,-1)$. Assume that the following condition holds:
if $x,y\in S$ belong to two geodesics $\gamma$ and
$\gamma'$, say $\gamma(n)=\gamma'(n')=x$,
$\gamma(m)=\gamma'(m')=y$, then $\gamma(n\trid
m)=\gamma'(n'\trid m')$. Then we can define on $S$ a unique binary
operation $\trid$ in such a way that the geodesics respect
$\trid$, namely, $x\trid y=\gamma(n\trid m)=\gamma(2n-m)$ for any
geodesic $\gamma$ such that $\gamma(n)=x$ and $\gamma(m)=y$. This
operation furnishes $S$ with the structure of an involutory crossed set if it
maps geodesics to geodesics; that is, if $x\trid\gamma$ is a
geodesic for any $x$ and any $\gamma$.
It can be shown that any involutory crossed set arises in this way \cite{jo}.

\subsubsection*{\bf Core crossed sets}
Let $G$ be a group. The \emph{core } of $G$ is the crossed set
$(G, \trid)$, where  $x\trid y=xy^{-1}x$. The core is an
involutory crossed set. If $G$ is abelian, its core is the affine
crossed set with $g=-\id$. More generally, one can define the core
of a Moufang loop.

\subsubsection*{\bf The free quandle of a set.}
Let $C$ be a set and let $F(C)$ be the free group generated by
$C$; let ${\mathcal O}_c$ denote the orbit of $c\in C$ in $F(C)$.
We claim that the standard crossed set $X_C := \bigcup_{c\in C}
{\mathcal O}_c$ is the free quandle on the set $C$.

For, let $\psi: C \to (X, \trid)$ be any function and let $\Psi:
F(C) \to \gax$ be the unique group homomorphism extending $c
\mapsto \phi_{\psi(c)}$. We define then $\widehat{\psi}: X_C \to
X$ by $$ \widehat{\psi}(y) = \Psi(u)(c), \qquad \text{if } y =
ucu^{-1}, c\in C. $$ The well-definiteness of $\widehat{\psi}$ is
a consequence of the following facts about orbits in free groups:

\begin{enumerate}
\item If $c, d\in C$ and ${\mathcal O}_c = {\mathcal O}_d$, then $c=d$.
\item The centralizer of $c\in C$ is $\langle c\rangle$.
\end{enumerate}

It is not difficult to see that $\widehat{\psi}$ is indeed a
morphism of quandles extending $\psi$, and the unique one.
It is easy to see that $X_C$ is also the free crossed set generated
by $C$. It is not, however, a free rack.

%%%%%%%%%%%%%%%%%%%%%%%%%%%%%%%%%%%%%%%%%%%%%%%%%%%%%%%%%%%%%%%
%%%%%%%%%%%%%%%%%% Seccion  Extensions %%%%%%%%%%%%%%%%%%%%%%%%
%%%%%%%%%%%%%%%%%%%%%%%%%%%%%%%%%%%%%%%%%%%%%%%%%%%%%%%%%%%%%%%
\section{Extensions}\label{sn:2}

%%%%%%%%%%%%%%%%%%%%%%%%%%%%%%%%%%%%%%%%%%%%%%%%%%%%%%%%%%%%%%%
%%%%%%%%%%%%%%%%%% SubSeccion  dynamical %%%%%%%%%%%%%%%%%%%%%%
%%%%%%%%%%%%%%%%%%%%%%%%%%%%%%%%%%%%%%%%%%%%%%%%%%%%%%%%%%%%%%%
\subsection{Extensions with dynamical cocycle}

We now discuss another way of constructing racks (resp. quandles,
crossed sets), generalizing Example \ref{prodcs}. The proof of the
following result is essentially straightforward.

\begin{lem}\label{lm:nac}%non-abelian cocycle
Let $X$  be a rack and let $S$ be a non-empty set. Let
$\al: X \times X \to \fun(S\times S, S)$ be a function, so that
for each $i,j \in X$ and $s, t\in S$ we have an element
$\al_{ij}(s, t) \in S$. We will write $\al_{ij}(s):S\to S$ the
function $\al_{ij}(s)(t)=\al_{ij}(s,t)$. Then $X\times S$ is a rack with
respect to $$(i, s) \trid (j, t) = (i\trid j, \al_{ij}(s, t))$$ if
and only if the following conditions hold:
\begin{flalign}
\label{cocycle0}& \al_{ij}(s) \text{ is a bijection};&\\
\label{cocycle3} & \al_{i,j\trid k}(s, \al_{jk}(t, u))
	= \al_{i \trid j, i \trid k}(\al_{ij}(s, t), \al_{ik}(s, u))
	\ \forall i,j, k \in X,\ s,t,u \in S.&
\end{flalign}
\hspace{2cm}in other words, $\al_{i,j\trid k}(s)\al_{j,k}(t)
	=\al_{i\trid j,i\trid k}(\al_{i,j}(s,t))\al_{i,k}(s)$.

\noindent If $X$ is a quandle, then $X\times S$ is a quandle iff further
\begin{flalign}
\label{cocycle1} &\al_{ii}(s, s) = s
	\quad\text{for all $i  \in X$ and $s \in S$.}&
\end{flalign}
If $X$ is a crossed set, then $X\times S$ is a crossed set iff further
\begin{flalign}
\label{cocycle2} & \al_{ji}(t, s) = s
\text{ whenever } i\trid j = j \text{ and  }  \al_{ij}(s, t) = t
	\quad\forall i,j \in X,\ s,t\in S. &
\end{flalign}
\hfill\qed
\end{lem}

\begin{defn}\label{df:dyncoc}
If these conditions hold we say that $\al$ is a \emph{dynamical cocycle}
and that $X\times S$ is an extension of $X$ by $S$; we shall denote it by
$X\times_{\al} S$.
When necessary, we shall say that $\al$ is a \emph{rack (resp. quandle,
crossed set) dynamical cocycle} to specify that we require it to
satisfy \eqref{cocycle0}$+$\eqref{cocycle3} (resp. 
\eqref{cocycle0}$+$\eqref{cocycle3}$+$\eqref{cocycle1},
\eqref{cocycle0}$+$\eqref{cocycle3}$+$\eqref{cocycle1}$+$\eqref{cocycle2}).
The presence of the parameter $s$ justifies the name of ``dynamical".
\end{defn}

Assume that $X$ is a quandle and let $\al$ be a quandle dynamical cocycle.
For $i\in X$ consider $s\trid_it := \al_{ii}(s,t)$. It is immediate
to see that $(S,\trid_i)$ becomes a quandle $\forall i\in X$.
Then \eqref{cocycle3}, when $ j =k$, says:
\begin{equation}
\label{cocycle4}  \al_{i,j}(s, t \trid_j u) = \al_{ij}(s, t)
\trid_{i\trid j} \al_{ij}(s, u), \qquad \forall s,t,u\in S.
\end{equation}
In other words, the map $\al_{ij}(s)$ is an isomorphism of quandles
$\al_{ij}(s): (S, \trid_j) \to (S, \trid_{i\trid j})$.

The projection $X\times_{\al} S \to X$ is clearly a
morphism. Conversely, it turns out that projections of indecomposable
racks (resp. quandles, crossed sets) are always extensions.
Before going over this, we state a technical lemma for further use.

\begin{lem}\label{lm:esext}
Let $(X,\trid)$ be a quandle which is a disjoint union
$X=\coprod_{i\in Y}X_i$ such that there exists $\bar\trid:Y\times
Y\to Y$ with $X_i\trid X_j=X_{i\bar\trid j}$. Suppose that
$\card(X_i)=\card(X_j)$ $\forall i,j$ (this holds for instance if
$X$ is indecomposable). Then $(Y,\bar\trid)$ is a quandle.

Furthermore, take $S$ a set such that $\card(S)=\card(X_i)$ and
for each $i\in Y$ set $g_i:X_i\to S$ a bijection. Let $\al:Y\times
Y\to\fun(S\times S,S)$ be given by $\al_{ij}(s,t):=g_{i\bar\trid
j}(g_i^{-1}(s)\trid g_j^{-1}(t))$. Then $\al$ is a dynamical
cocycle and $X\simeq Y\times_{\al}S$.
\end{lem}
\pf This follows without troubles from Lemma \ref{lm:nac}.\epf

\begin{rem}
\begin{enu}
\item Within the hypotheses of the lemma, if $X$ is indecomposable
then so is $Y$.
\item The whole lemma can be stated in terms of racks.
\item In order to state the lemma in terms of crossed sets, it is
necessary to further assume that $(Y,\bar\trid)$ is a crossed set,
i.e., that it satisfies \eqref{ccc3}.
\end{enu}
\end{rem}

\begin{cor}\label{co:projection}
Let $(X,\trid)$, $(Y,\bar\trid)$ be quandles (resp. racks, crossed
sets). Let $f:X\to Y$ be a surjective morphism such that the fibers
$f^{-1}(y)$ all have the same cardinality (this happens for
instance if $X$ is indecomposable). Then $X$ is an extension
$X=Y\times_{\al} S$. \hfill\qed
\end{cor}

Let $X$ be a rack. Let $\al:X\times X\to \fun(S\times S,S)$ be a
dynamical cocycle and let 
$\gamma:X\to \Sim_S$ be a function. 
Define $\al':X\times X\to \fun(S\times S,S)$ by
\begin{equation}\label{eq:2cdch}%2 cociclos dinamicos cohomologos
\al'_{i,j}(s,t)
=\gamma_{i\trid j}(\al_{i,j}(\gamma_i^{-1}(s),\gamma_j^{-1}(t))),
\qquad\text{i.e., }\al_{ij}'(s)
=\gamma_{i\trid j}\al_{ij}(\gamma_i^{-1}(s))\gamma_j^{-1}.
\end{equation}
Then $\al'$ is a dynamical cocycle and we have an isomorphism of
racks $T:(X\times_{\al} S)\to(X\times_{\al'} S)$ given by
$T(i,s)=(i,\gamma_i(s))$.\newline
Conversely, if there is an isomorphism of racks
$T:(X\times_{\al} S)\to(X\times_{\al'} S)$ which commutes with
the canonical projection $X\times S\to X$ then there exists $\gamma:
X\to \Sim_S$ such that $\al$ and $\al'$ are related as in \eqref{eq:2cdch}.

\begin{defn}
We say that $\al$ and $\al'$ are \emph{cohomologous} if and only if 
there exists $\gamma:
X\to \Sim_S$ such that $\al$ and $\al'$ are related as as in \eqref{eq:2cdch}.
\end{defn}

\begin{exmp}\label{ex:extcubo}
Let $Y$ be the crossed set given by the faces of the cube. Then
$Y$ is the disjoint union of the subsets made out of the pairs of
opposite faces. This union satisfies the hypotheses of Lemma
\ref{lm:esext} and, being indecomposable, the quotient
$(X,\bar\trid)$ is isomorphic to the crossed set $(\Z_3,\trid^2)$.
\end{exmp}

\begin{exmp}\label{ex:extensionafin}
Let $(A,g)$ be an affine crossed set. Suppose that there exists a
subgroup $B\subseteq A$ invariant by $g$; let $\overline g$ be the
induced automorphism of $A/B$. Consider the affine crossed set
$(A/B, \overline g)$;  the projection $(A,g)\fllad{\pi}(A/B, \overline
g)$ is a morphism of crossed sets. Corollary \ref{co:projection}
applies and we see that $A$ is an extension of $A/B$.
\end{exmp}

More examples appear in \cite{chns} by means of group extensions.
They are used to color twist-spun knots.

%%%%%%%%%%%%%%%%%%%%%%%%%%%%%%%%%%%%%%%%%%%%%%%%%%%%%%%%%%%%%%%
%%%%%%%%%%%%%%%%%% SubSeccion  constant  %%%%%%%%%%%%%%%%%%%%%%
%%%%%%%%%%%%%%%%%%%%%%%%%%%%%%%%%%%%%%%%%%%%%%%%%%%%%%%%%%%%%%%
\subsection{Extensions with constant cocycle}\label{ssecc}

Let $X$ be a rack.
Let $\be:X \times X \to \Sim_S$. We say that $\be$ is a
\emph{constant rack cocycle} if
\begin{flalign}
\label{cocycle7} &\be_{i,j\trid k}\be_{j,k}=\be_{i\trid j,i\trid k}\be_{i,k}.& \\
\intertext{If $X$ is a quandle, we say that $\be$ is a
	\emph{constant quandle cocycle} if it further satisfies}
\label{cocycle5} &\be_{ii}=\id, \quad\forall i\in X.& \\
\intertext{If $X$ is a crossed set, we say that $\be$ is a
	\emph{constant crossed set cocycle} if it further satisfies}
\label{cocycle6} & \be_{ji}=\id
\text{ whenever } i\trid j = j \text{ and  } \be_{ij}(t) = t
	\text{ for some $t\in S$.}&
\end{flalign}

We have then an extension $X\times_{\be} S:=X\times_{\al} S$,
taking $\al_{ij}(s,t)=\be_{ij}(t)$. Note that $\trid_i$ is trivial
for all $i$, and the fiber $F_{\phi_{(i,s)}}=F_{\phi_i}\times S$.

We shall say in this case that the extension is \emph{non-abelian}.
It is clear that an extension $X\times_\al S$ is non-abelian if and only if
$\al_{ij}(s)=\al_{ij}(t)\quad\forall s,t\in S$, $\forall i,j\in X$.

\begin{defn}\label{de:dynamico-cohomol}
Let $\gamma:X\to \Sim_S$ be a function and let $\be$ be a constant
cocycle. Define $\be':X\times X\to \Sim_S$ by
$$\be'_{i,j}=\gamma_{i\trid j}\,\be_{i,j}\,\gamma_j^{-1}.$$
Then we have an isomorphism $T:(X\times_{\be} S)\to(X\times_{\be'}S)$
given by $T(i,s)=(i,\gamma_i(s))$.
In this case, we shall say that $\be$ and $\be'$ are \emph{cohomologous}.
\end{defn}

The use of the word ``cocycle" is not only suggested by its analogy
with group $2$-cocycles, which describe extensions: there is a
general definition of abelian cohomology (see Section \ref{sn:4})
for which this is its natural non-abelian counterpart. The use of
the word ``cocycle" in the phrase ``dynamical cocycle" stands on the
same basis.

For $X\times_{\be} S$ a non-abelian extension, let $\psi:X\to
\Inn_{\trid}  (X\times_{\be} S)$ be given by $\psi_i=\phi_{(i,t)}$
for an arbitrary $t\in S$; that is, $\psi_i(j,s)=(i \trid j,
\be_{i,j}(s))$. Then $\psi(X)$ generates $\Inn_{\trid}(X\times_{\be} S)$.
Let $H_i=\{h\in \Inn_\trid(X\times_\be S)\ |
	\ h(i,s)\in i\times S\ \forall s\in S\}$.

\begin{defn}\label{de:transitive}
Assume that $X$ is indecomposable. A constant cocycle $\be:X\times
X\to\Sim_S$ is \emph{transitive} if for some $i\in X$, the group
$H_i$ acts transitively on $i\times S$. Note that this definition
does not depend on $i$.
\end{defn}

We have seen that all the fibers of an indecomposable
rack (resp. quandle, crossed set) have the same cardinality;
we provide now a precise description of an indecomposable rack
(resp. quandle, crossed set). Recall the map $\phi:Y\to\gay$
from Definition \ref{de:ddp}.

\begin{prop}
Let $Y$ be an indecomposable rack (resp. quandle, crossed set),
let $X=\phi(Y)$ and let $S$ be a set with the cardinality of the
fibers of $\phi$. Then we have an isomorphism
$T:Y\to X\times_{\be}S$ for some constant cocycle $\be$.

Conversely, a non-abelian extension $X \times_{\be} S$ is
indecomposable if and only if $X$ is indecomposable and $\be$ is
transitive.
\end{prop}

\pf Choose, for each $x\in X$, a bijection $g_x:F_x=\phi^{-1}(x)\to S$.
We have then a bijection $T:Y\to X\times S$,
$T(i)=(\phi_i,g_{\phi_i}(i))$. We define, for $x,y\in X$ and $s\in S$,
$$\be_{xy}(s)=g_{x\trid y}(T^{-1}(x,s)\trid T^{-1}(y,s)).$$
It is straightforward to see that $\be$ is a constant cocycle
and that $T$ is an isomorphism.

The second part is clear.
\epf

\begin{exmp}
Let $X=\{1,2,3,4\}$ be the tetrahedral crossed set defined in
Subsection \ref{ss:examples} and let $S=\{a,b\}$. Then
$\Sim_S=\{\id,\sigma\}\simeq C_2=\{\pm 1\}$. There is a
non-trivial $2$-cocycle $\be:X\times X\to S$ given by
$$\begin{cases} \be(x,y)=1 &\text{ if } x=1 \text{ or } y=1 \text{
or } x=y, \\ \be(x,y)=-1 &\text{ otherwise.}
\end{cases}$$
Let $\psi:X\to \Inn_{\trid}  (X\times_{\be} S)$ be as in the
paragraph preceding Definition \ref{de:transitive}. It is clear
that $\psi(1)\psi(2)\in H_4$, and $\psi(1)\psi(2)(4\times
a)=\psi(1)(3\times b)=4\times b$, whence $\be$ is transitive and
$X\times_{\be}S$ is indecomposable.
\end{exmp}

A way to construct cocycles, which resembles the classification of
Yetter--Drinfeld modules over group algebras, is the following one:
\begin{exmp}
Let $X$ be an indecomposable (finite) rack, $x_0\in X$ a fixed element,
$G=\gax$, and let $H=G_{x_0}$ be the subgroup of the inner isomorphisms
which fix $x_0$. Let $Z$ be a (finite) set and $\rho:H\to\Sim_Z$ a group
homomorphism. There is then a bijection $G/H\to X$ given by $g\mapsto g(x_0)$.
Fix a (set theoretical) section $s:X\to G$; i.e.,
$s(x)\cdot x_0=x\ \forall x\in X$. This determines, for each $x,y\in X$,
an element $t_{x,y}\in H$ such that $\phi_xs(y)=s(x\trid y)t_{x,y}$.
To see this, we compute
$$s(x\trid y)^{-1}\phi_xs(y)\cdot x_0
	=s(x\trid y)^{-1}\phi_x\cdot y
	=s(x\trid y)^{-1}\cdot(x\trid y)=x_0.$$
Then it is straightforward to see that $\beta:X\times X\to\Sim_Z$,
$\beta_{x,y}=\rho(t_{x,y})$ is a constant cocycle. Explicitly, this
defines an extension $X\times_{\beta}Z$ as
$$(x,z)\trid(x',z')=(x\trid x', \rho(s(x\trid x')^{-1}\phi_xs(x'))(z')).$$
\end{exmp}

Even if $X$ is a quandle, this is not a quandle in general, since it
does not necessarily satisfy \eqref{cocycle5}. Let us compute when
\eqref{cocycle5} is satisfied (suppose that $X$ is a quandle)
$$\beta_{x,x}=\rho(s(x\trid x)^{-1}\phi_xs(x))
	=\rho(s(x)^{-1}\phi_xs(x))=\rho(\phi_{s(x)^{-1}\cdot x})
	=\rho(\phi_{x_0}).$$
Then $X\times_{\beta}Z$ is a quandle iff $X$ is a quandle and
$\rho(\phi_{x_0})=1\in\Sim_Z$.

\begin{rem}
It is easy to see that for polyhedral quandles and affine quandles the
group $H$ is generated by $\phi_{x_0}$, and then we can not construct
non-trivial extensions in this way. 
\end{rem}

%%%%%%%%%%%%%%%%%%%%%%%%%%%%%%%%%%%%%%%%%%%%%%%%%%%%%%%%%%%%%%%
%%%%%%%%%%%%%%%% SubSeccion  X-modules %%%%%%%%%%%%%%%%%%%%%%%
%%%%%%%%%%%%%%%%%%%%%%%%%%%%%%%%%%%%%%%%%%%%%%%%%%%%%%%%%%%%%%%
\subsection{Modules over a rack}

Throughout this subsection, $\rs$ will denote the category of
racks. All the constructions below can be performed in 
the category $\qs$ of quandles, or in  the category $\cs$ of crossed sets.

It is clear that finite direct products exist in the category
$\rs$, \emph{cf.} Example \ref{prodcs}. Then we can
consider group objects in $\rs$; they are determined by the
following Proposition.

\begin{prop}\label{pr:grobj} A group object in $\rs$ is given by
a triple $(G,s,t)$, where $G$ is a group, $s\in\Aut(G)$,
$t:G\to Z(G)$ is a group homomorphism, and
\begin{itemize}
\item $st=ts$,
\item $s(x)x^{-1}t(x)\in\ker(t)\quad\forall x\in G$.
\end{itemize}
The rack structure is given by $x\trid y=t(x)s(y)$.

A group object $(G,s,t)$ in $\rs$ is a quandle iff $t(x)=s(x)^{-1}x$
$\forall x\in G$, iff it is an homogeneous quandle (hence crossed set).

A group object $(G,s,t)$ in $\rs$ is abelian iff it is an affine rack.
\end{prop}

\pf The second and third statements follow from the first
without difficulties (notice that the image of $t$ lies in the center).

A group object in $\rs$ is a triple $(G, \cdot, \trid)$ with $(G, \cdot)$
a group and $(G, \trid)$ a rack, such that the multiplication $\cdot$ is a
morphism of racks. Let $s, t: G \to G$, $s(x) = 1\trid x$,
$t(x) = x\trid 1$. Then both $s$ and $t$ are group homomorphisms since
$s(xy)=1\trid(xy)=(1\cdot1)\trid(x\cdot y)=(1\trid x)\cdot(1\trid y)=s(x)s(y)$,
and analogously for $t$. Furthermore,
\begin{align*}
x\trid y &=(1\cdot x)\trid (y\cdot 1)=(1\trid y)\cdot(x\trid 1)
	=s(y)t(x)\\
& = (x\trid 1) \cdot (1\trid y) =t(x)s(y).
\end{align*}
Then, $s$ must be an isomorphism, and then $t(x)$ is central $\forall x$.
Last,
\begin{align*}
x\trid(y\trid z) &= t(x)\cdot st(y)\cdot s^2(z), \\
(x\trid y)\trid (x\trid z) &= t^2(x)\cdot ts(y)\cdot st(x)\cdot s^2(z),
\end{align*}
whence $t(x)\cdot st(y)=t^2(x)\cdot ts(y)\cdot st(x)$. Taking $x=1$
we see that $st=ts$. Taking $y=1$ we see that $t(t(x)s(x)x^{-1})=1$.
The converse is not difficult.
\epf

Let us now consider the ``comma category" $\rsx$ over a fixed
rack $X$; recall that the objects of $\rsx$ are the maps
$f:Y \to X$ and the  arrows between $f:Y \to X$ and $g:Z \to X$
are the commutative triangles, \emph{i.e.} the maps $h: Y \to Z$
such that $gh = f$.

Since equalizers  exist in $\rs$ (they are just the
set-theoretical equalizers with the induced $\trid$),  $\rs$
has finite limits. It follows that $\rsx$ also has finite limits.
We are willing to determine all abelian group objects in $\rsx$. 

\begin{defn} Let $X$ be a rack and let $A$ be an abelian group.
A \emph{structure of $X$-module} on $A$ consists of the following
data:

\begin{itemize}
\item a family $(\eta_{ij})_{i,j\in X}$ of automorphisms of $A$, and
\item a family $(\tau_{ij})_{i,j\in X}$ of endomorphisms of $A$,
\end{itemize}
such that the following axioms hold:
\begin{flalign}
\label{xmod2}& \eta_{i,j\trid k}\, \eta_{j,k} = \eta_{i\trid j,i\trid k} \, \eta_{i,k},& \\
\label{xmod3}& \tau_{i,j\trid k} = \eta_{i\trid j,i\trid k} \, \tau_{i,k}
	+ \tau_{i\trid j,i\trid k} \, \tau_{i,j}, &\\
\label{xmod5} & \eta_{i,j\trid k} \, \tau_{j,k}
	= \tau_{i\trid j,i\trid k} \, \eta_{i,j}.& \\
\intertext{If $X$ is a quandle, a \emph{quandle structure of $X$-module}
	on $A$ is a structure of $X$-module which further satisfies}
\label{xmod1}& \eta_{ii} +  \tau_{ii} = \id. & \\
\intertext{If $X$ is a crossed set, a \emph{crossed set structure of $X$-module}
	on $A$ is a quandle structure of $X$-module which further satisfies}
\label{xmod4}& \text{if } i\trid j=j\text{ and }(\id-\eta_{ij})(t)=\tau_{ij} (s)
	\text{ for some  } s,t \text{ then }(\id-\eta_{ji})(s)=\tau_{ji} (t),&
\end{flalign}

An \emph{$X$-module} is an abelian group $A$ provided with a structure of
$X$-module.
\end{defn}

\begin{rem} Taking $j =k$ in \eqref{xmod2}, one gets in presence of
\eqref{ccc2} the suggestive equality:% que me significa?
\begin{equation}\label{xmod2bis}
\eta_{i,j}\, \eta_{j,j} = \eta_{i\trid j,i\trid j} \, \eta_{i,j}.
\end{equation}
\end{rem}

Let $A$ be an $X$-module. We define $\al_{ij}: A\times A\to A$ by
$$ \al_{ij} (s,t) := \eta_{ij} (t) + \tau_{ij} (s). $$

\begin{thm}
\begin{enumerate}
\item\label{cucua}
	$\al_{ij}$ is a dynamical cocycle, hence we can form the rack
	$Y = X\times _{\al} A$.
\item The canonical projection $p:Y\to X$ is an abelian group in $\rsx$.
\item If $p: Y\to X$ is an  abelian group in $\rsx$ and $X$ is
	indecomposable, then $Y \simeq X\times _{\al} A$ for some
	$X$-module $A$ and $p$ is the canonical projection.
\item If $X$ is a quandle and $A$ is a quandle $X$-module,
	then the preceding statements are true in the category of quandles.
	Same for crossed sets.
\end{enumerate}
\end{thm}

\pf
\begin{enu}
\item Condition \eqref{cocycle0} follows since $\eta_{ij}$ is
an automorphism.
The left hand side of \eqref{cocycle3} is
$$\eta_{i,j\trid k}\, \eta_{j,k} (u) + \eta_{i,j\trid k} \,
\tau_{j,k} (t) + \tau_{i,j\trid k}(s) $$ and the right hand side
of \eqref{cocycle3} is $$\eta_{i\trid j,i\trid k} \, \eta_{i,k}
(u) + \tau_{i\trid j,i\trid k} \, \eta_{i,j} (t) + \eta_{i\trid
j,i\trid k} \, \tau_{i,k} (s) + \tau_{i\trid j,i\trid k} \,
\tau_{i,j}(s).$$ Thus, \eqref{cocycle3} follows from
\eqref{xmod2}, \eqref{xmod3} and \eqref{xmod5}.
\item
Let $\sigma: X \to X\times _{\al} A$, $\sigma(i) = (i,0)$,
$i\in X$; it is clearly a morphism of racks.
Let $+:Y\times_X Y \to Y$, $(i,a)+(i,b) = (i,a+b)$; it is clearly a
morphism in $\rsx$. It is not difficult to verify that $(Y, +)$ is
an abelian group in $\rsx$ with identity element $\sigma$.
\item
Let  $p: Y\to X$ be an  abelian group in $\rsx$, $X$ an
arbitrary rack. We have morphisms in $\rsx$ $+:Y\times_X Y \to Y$
and $\sigma: X \to Y$ satisfying the axioms of abelian group.
In particular, we have:

\medskip
\begin{enu2}
\item\label{gaem1}
The existence of $\sigma$ implies that $p$ is surjective; if $i\in X$,
$A_i := p^{-1}(i)$ is an abelian group with identity $\sigma(i)$.
\item\label{gaem2}
$(a + b) \trid (c + d) = (a\trid c) + (b\trid d) \in A_{i\trid j}$,
if $a,b \in A_i$, $c,d \in A_j$.
\item\label{gaem3}
The map $h_{ij}: A_j \to A_{i\trid j}$, $h_{ij}$ the restriction
of $\phi_{\sigma(i)}$, is an isomorphism of abelian groups, for
all $i, j\in X$.
\end{enu2}

\medskip
Assume now that $X$ is indecomposable. Then all the abelian groups
$A_i$ are isomorphic, by (\ref{gaem3}). Fix an abelian group $A$ provided
with group isomorphisms $\gamma_i: A_i \to A$. Define $\al_{ij}:
A\times A \to A$, $i, j\in X$, by $$ \al_{ij}(s,t):= \gamma_{i
\trid j}(\gamma_i^{-1}(s)\trid \gamma_j^{-1}(t)). $$ We claim that
$\al_{ij}(s,t)=\eta_{ij}(t) + \tau_{ij}(s)$, where $\eta_{ij}(t) =
\al_{ij}(0,t)$ and $\tau_{ij}(s) = \al_{ij}(s,0)$; this follows
without difficulty from (\ref{gaem2}), since the $\gamma$'s are linear. Now,
$\eta_{ij}(t) = \gamma_{i\trid j}h_{ij}\gamma_j^{-1}$ is a linear
automorphism, whereas $\tau_{ij}$ is linear by (\ref{gaem2}). We need
finally to verify conditions \eqref{xmod2}--\eqref{xmod5};
this is done reversing the arguments in part (\ref{cucua}). 
\item
Condition \eqref{cocycle1} amounts in the present case to \eqref{xmod1},
whereas condition \eqref{cocycle2} amounts to \eqref{xmod4}.
\end{enu}
\epf

\begin{rem}
Assume now that $X$ is a \emph{non-indecomposable}
quandle and keep the notation of the proof. Then
$A_i$ is a subquandle of $Y$, indeed an abelian group in $\qs$.
By Proposition \ref{pr:grobj}, $A_i$ is affine, with respect to
some $g_i\in \Aut(A_i)$.
\end{rem}

\begin{exmp} If $(A, g)$ is an affine crossed set, then it is an
$X$-module over any rack $X$ with $\eta_{ij} = g$, $\tau_{ij}=f=\id-g$.
We shall say that $A$ is a \emph{trivial} $X$-module if $g =\id$,
that is when it is trivial as crossed set.
\end{exmp}

\begin{exmp} Let $X$ be a trivial quandle, let $(A_i)_{i\in X}$ be a
family of affine quandles and let $Y$ be the disjoint sum of the
$A_i$'s, with the evident projection $Y\to X$. Then $Y\to X$ is an
abelian group in $\qsx$ which is not an extension of $X$ by any
$X$-module, if the $A_i$'s are not isomorphic.
\end{exmp}

We now show that the category of $X$-modules, $X$ an
indecomposable quandle or rack, is abelian with enough
projectives. Actually, it is equivalent to the category of modules
over a suitable algebra.

\begin{defn}\label{defrackalg}
The \emph{rack algebra} of a rack $(X,\trid)$ is the $\Z$-algebra $\ral X$
presented by generators $(\eta_{ij}^{\pm 1})_{i,j\in X}$ and
$(\tau_{ij})_{i,j\in X}$, with relations
$\eta_{ij} \eta_{ij}^{-1} = \eta_{ij}^{- 1}\eta_{ij} = 1$,
\eqref{xmod2}, \eqref{xmod3} and \eqref{xmod5}.

The \emph{quandle algebra} of a quandle $(X,\trid)$ is the $\Z$-algebra
$\qal X$ presented by generators $(\eta_{ij}^{\pm 1})_{i,j\in X}$
and $(\tau_{ij})_{i,j\in X}$, with relations $\eta_{ij}
\eta_{ij}^{- 1} = \eta_{ij}^{- 1}\eta_{ij} = 1$,
\eqref{xmod2}, \eqref{xmod3}, \eqref{xmod5} and \eqref{xmod1}.
\end{defn}

It is evident that the category of $X$-modules, $X$ a rack, is
equivalent to the category of modules over $\ral X$; therefore, it
is abelian with enough injective and projective objects.
Same for quandle $X$-modules and $\qal X$.

The algebra $\qal X$ is augmented  with augmentation
$\varepsilon:\qal X \to \Z$, $\varepsilon(\eta_{ij})=1$,
$\varepsilon(\tau_{ij}) = 0$. The algebra $\ral X$ is also
augmented, composing $\varepsilon$ with the projection
$\ral X\to\qal X$.

There are various interesting quotients of  the algebra
$\qal X$.

First, the quotient of $\qal X$ by the ideal generated by the
$\tau_{ij}$'s is isomorphic to the group algebra of the group
$\Lambda(X)$ presented by generators $(\eta_{ij})_{i,j\in X}$
with relations  \eqref{xmod2}.

Next, consider the following elements of the group
algebra $\Z G_X$ (see Definition \ref{gx}):
\begin{equation}\label{defrestmod}
\eta_{ij} := i, \quad \tau_{ij} := 1 - (i\trid j).
\end{equation}

It is not difficult to see that this defines a surjective algebra
homomorphism $\qal X \to \Z G_X$; in particular any $G_X$-module
has a natural structure of $X$-module.

\begin{defn}\label{def:restmod}
If $X$ is a crossed set and $M$ is any $G_X$-module then $M$ also
satisfies the extra condition \eqref{xmod4}. We shall say that $M$
is a \emph{restricted $X$-module}.
\end{defn}

\begin{defn} Let $X$ be a rack and let $A$ be an $X$-module.
A \emph{$2$-cocycle on $X$ with values in $A$} is a collection
$(\kappa_{ij})_{i, j\in X}$ of elements in $A$ such that
\begin{equation}\label{generaldoscociclo}
\eta_{i,j\trid k}(\kappa_{jk}) +  \kappa_{i,j\trid k}
	= \eta_{i\trid j,i\trid k}(\kappa_{ik})
	+ \tau_{i\trid j,i\trid k}(\kappa_{ij})
	+ \kappa_{i\trid j,i\trid k}\quad \forall i, j, k\in X.
\end{equation}
Two $2$-cocycles $\kappa$ and $\kappa'$ are \emph{cohomologous} iff
$\exists f:X\to A$ such that
\begin{equation}\label{eq:2cch}%2 cociclos cohomologos
\kappa'_{ij}=\kappa_{ij}-\eta_{ij}(f(j))
	+f(i\trid j)-\tau_{ij}(f(i)).
\end{equation}
\end{defn}
As remarked earlier, the reader can find in section \ref{sn:4} a complex
which justifies these names.

\begin{prop}\label{affinemodule-pr}
\begin{enu}
\item Let $X$ be a rack and consider functions $\eta:X\times X\to\Aut(A)$,\linebreak
$\tau:X\times X\to\End(A)$, $\kappa:X\times X\to A$.
Let us define a map $\al$ by
\begin{equation}\label{affinemodule-eq}
\al_{i,j} (a,b) = \eta_{i,j} (b) + \tau_{i,j} (a) + \kappa_{ij},
\qquad a,b \in A.
\end{equation}
Then the following conditions are equivalent:
\begin{enu2}
\item $\al$ is a dynamical cocycle.
\item $\eta,\tau$ define a structure of an $X$-module on $A$ and
	$(\kappa_{ij})$ is a $2$-cocycle, i.e. it satisfies
	\eqref{generaldoscociclo}.
\end{enu2}

\item Let $\kappa$ and $\kappa'$ be $2$-cocycles and let $\alpha$,
$\alpha'$ be their respective dynamical cocycles. If $\kappa$ and
$\kappa'$ are cohomologous then $\alpha$ and $\alpha'$ are cohomologous.
\end{enu}
\end{prop}
\pf Straightforward.
For the second part, if $f$ is as in \eqref{eq:2cch}, then take
$\gamma:X\to\Sim_A$ by $\gamma_i(s)=f(i)+s$.
It is easy to verify \eqref{eq:2cdch}.
\epf

\begin{defn}\label{affinemodule-def}
If $X$ is a rack, $A$ is an $X$-module and $(\kappa_{ij})_{i,j\in X}$
is a $2$-cocycle on $X$ with values in $A$, then the extension
$X\times_{\al} A$, where $\al$ is given by \eqref{affinemodule-eq},
is called an \emph{affine module over $X$}.
By abuse of notation, this extension will be denoted 
$X\times_{\kappa} A$.
\end{defn}

All the constructions and results in this section
are valid more generally over a fixed commutative ring $R$.

%%%%%%%%%%%%%%%%%%%%%%%%%%%%%%%%%%%%%%%%%%%%%%%%%%%%%%%%%%%%%%%
%%%%%%%%%%%%%%%% Seccion  Simple racks  %%%%%%%%%%%%%%%%%%%%%%%
%%%%%%%%%%%%%%%%%%%%%%%%%%%%%%%%%%%%%%%%%%%%%%%%%%%%%%%%%%%%%%%
\section{Simple racks}\label{sn:3}

%%%%%%%%%%%%%%%%%%%%%%%%%%%%%%%%%%%%%%%%%%%%%%%%%%%%%%%%%%%%%%%
%%%%%%%%%%%%%%%% SubSeccion  Faithrul C.S  %%%%%%%%%%%%%%%%%%%%
%%%%%%%%%%%%%%%%%%%%%%%%%%%%%%%%%%%%%%%%%%%%%%%%%%%%%%%%%%%%%%%
\subsection{Faithful indecomposable crossed sets}

To classify indecomposable racks, we may first consider faithful
indecomposable crossed sets (recall that a faithful rack is
necessarily a crossed set), and next compute all possible extensions.

\begin{prop} Let $X$ be an indecomposable finite  rack. Then $X$ is
isomorphic to an extension $Y\times_{\al} S$, where $Y$ is a
faithful indecomposable crossed set. Furthermore $Y$ can be chosen uniquely
with the property that any surjection $X\to Z$ of racks,
with $Z$ faithful, factorizes through $Y$. \end{prop}

\pf Consider the sequence $X \to X_1 := \phi(X) \to X_2 := \phi
(X_1) \dots $. Since $X$ is finite the sequence stabilizes, say at
$Y = X_n$, which is clearly faithful and indecomposable. By
Corollary \ref{co:projection}, $X \simeq Y\times_{\al} S$. Now let
$\psi: X\to Z$ be a surjection, with $Z$ faithful. By Lemma
\ref{le:gammafuntor}, it gives a surjection $\psi_1: X_1\to
Z$, and so on.
\epf

Faithful indecomposable crossed sets can be characterized as follows.

\begin{prop}
\begin{enu}
\item\label{pipu1}
If $X$ is a faithful indecomposable crossed set, then there
exists a group $G$ and an injective morphism of crossed sets
$\varphi: X \to G$, such that $Z(G)$ is trivial and  $\varphi(X)$
is a single orbit generating $G$ as a group. Furthermore, $G$ is
unique up to isomorphisms with these conditions.
\item\label{pipu2}
If $X$ is a single orbit in a group $G$ with $Z(G)$ trivial
and  $X$ generates $G$, then $X$ is a faithful indecomposable
crossed set. 
\item
There is an equivalence of categories between
\begin{enu2}
\item The category of faithful indecomposable crossed sets, with
surjective morphisms.
\item The category of pairs $(G, {\mathcal O})$, where
$G$ is a group with trivial center and $\mathcal O$ is an orbit
generating $G$; a morphism $f: (G, {\mathcal O}) \to (K, {\mathcal U})$
is a group homomorphism $f: G \to K$ such that $f({\mathcal O}) = \mathcal U$.
\end{enu2}
\end{enu}
\end{prop}

\pf
\begin{enu}
\item Existence: take $G=\gax$; uniqueness: by Lemma \ref{le:gammastandard}.
\item Immediate.
\item Follows from (\ref{pipu1}), (\ref{pipu2}) and Lemma \ref{le:gammafuntor}.
\end{enu}
\epf

We shall  say that a projection (= surjective homomorphism)
$\pi: X \to Y$ of racks is \emph{trivial} if $\card Y$ is either 1 or $\card X$.

\begin{defn} A rack $X$ is \emph{simple} if it is not trivial and
%$\card X > 1$, it is different from the trivial rack with $2$ elements, and
any projection of racks $\pi: X \to Y$  is  trivial.
\end{defn}

A decomposable rack has a projection onto the trivial rack with two elements; 
it follows that a simple rack is indecomposable.

Let $X$ be a simple rack with $n$ elements; then $\phi(X)$ has
only one point, or $\phi$ is a bijection. In the first case, $X$
is a permutation rack defined by a cycle of length $n$; in the
second case $X$ is a crossed set (and necessarily $\card X > 2$).

It is not difficult to see that a permutation rack 
with $n$ elements defined by a cycle of length $n$ is simple if and only if
$n$ is prime.

\begin{prop}\label{simplesprimerresultado}
Let $X$ be an indecomposable faithful crossed set, which
corresponds to a pair $(G,\mathcal{O})$. Then $X$ is simple if and
only if any quotient of $G$ different from $G$ is cyclic.
\end{prop}

\pf
Assume $X$ is simple. If $\pi:G\to K$ is an epimorphism of groups,
then either $\pi(X)$ is a single point or $\pi|_X$ is bijective.
If $\pi(X)$ is a point, this point generates $K$, and then $K$
must be cyclic. If $\pi|_X$ is bijective, take $h\in\ker(\pi)$;
then $hxh^{-1}=x$ $\forall x\in X$, whence $h\in Z(G)$. This means
that $\pi$ is bijective. 

Assume now that any non-trivial quotient
of $G$ is cyclic. Let $\pi:X\to Y$ be a surjective morphism of
crossed sets, then $G\simeq\gax\to\gay$ is an epimorphism, so that
either it is a bijection (then $X\simeq Y$) or $\gay$ is cyclic.
In the last case, $\phi(Y)$, being indecomposable, is a point;
hence $Y$, being also indecomposable, is a point.
\epf

\begin{exmp}
Let $X$ be the crossed set of the faces of the cube. We can
realize $X$ as the orbit given by $4$-cycles inside $\Sim_4$,
since $X$ is faithful and $\gax\simeq\Sim_4$. Considering the
quotient $\Sim_4\to\Sim_4/K\simeq\Sim_3$, where $K$ is the Klein
subgroup, we see that $X$ is not simple. Indeed, this gives $X$ as
the same extension $\Z_3\times_{\al}\Z_2$ as in Example
\ref{ex:extcubo}. Also, if $X'$ is the orbit given by the $6$
transpositions in $\Sim_4$, we see taking the same quotient
$\Sim_4/K$ that $X'$ is an extension $\Z_3\times_{\al'}\Z_2$.
\end{exmp}

\begin{exmp}
Since the only proper quotient of $\Sim_n$ ($n\ge 5$) is cyclic,
we see that if $X$ is a non-trivial orbit of $\Sim_n$ then either
\begin{itemize}
\item $X$ generates $\Sim_n$ and then it is simple, or
\item $X\subset\mathbb{A}_n$.
\end{itemize}
If $X\subset\mathbb{A}_n$, then it might fail to be an orbit in $\mathbb{A}_n$.
This would happen if and only if the centralizers of the elements
of $X$ lie inside $\mathbb{A}_n$ (because the order of the orbits is the
ratio between the order of the group and the order of the
centralizers). In this case, thus, $X$ decomposes as a union of
two orbits, which are isomorphic via conjugation by any element in
$\Sim_n\setminus\mathbb{A}_n$ (an example of this case arises for $n=5$
and $X$ the $5$-cycles). On the other hand, if the centralizers of
the elements of $X$ are not included in $\mathbb{A}_n$ then $X$ is
indecomposable, and hence simple by the proposition.
\end{exmp}

%%%%%%%%%%%%%%%%%%%%%%%%%%%%%%%%%%%%%%%%%%%%%%%%%%%%%%%%%%%%%%%
%%%%%%%%%%%%%%%%%%%%%% clasificacion %%%%%%%%%%%%%%%%%%%%%%%%%%
%%%%%%%%%%%%%%%%%%%%%%%%%%%%%%%%%%%%%%%%%%%%%%%%%%%%%%%%%%%%%%%

\subsection{Classification of simple racks}\label{ssn:sr}

We  characterize now simple racks in group-theoretical
terms. We  first  classify finite groups $G$ such that $Z(G)$
is trivial and $G/N$ is cyclic for any normal non-trivial
subgroup of $G$. We are grateful to R. Guralnick for help in this
question.

To fix notation, if $G$ acts on $H$, we put $H\rtimes G$ the semidirect
product with structure
$$(h,g)(h',g')=(h(g\cdot h'),gg').$$

\begin{thm}[Guralnick]\label{clasigur2}
Let $G$ be a non-trivial finite group such that $Z(G)$ is trivial and $G/H$ is
cyclic for any normal non-trivial subgroup $H$. Then there are a
simple group $L$, a positive integer $t$ and a finite
cyclic group $C=<x>$ in $\Aut (N)$, where $N := L^t=L\times\cdots\times L$ ($t$
times), such that one of the following possibilities hold.
\begin{enumerate}
\item\label{clasiguri}
    $L$ is abelian, so that $N$ is elementary abelian of order $p^t$,
    $x$ is not trivial and it acts irreducibly on $N$. Furthermore,
    $G\simeq N\rtimes C$.
%   and $x$ acts irreducibly on $N$; if $t = 1$, $x$ is non-trivial.
\item\label{clasigurii}
    $L$ is simple non-abelian,
%   $x$ is non-trivial if $t>1$,
%   $C$ is non-trivial if $t>1$,
    $G=NC\simeq (N\rtimes C)/Z(N\rtimes C)$ and $x$ acts by
    \begin{equation}\label{acciondex2}
    x\cdot(l_1,\ldots,l_t)=(\theta(l_t),l_1,\ldots,l_{t-1})
    \end{equation}
    for some $\theta\in\Aut(L)$.
\end{enumerate}
Conversely, all the groups in (\ref{clasiguri})
or (\ref{clasigurii}) have the desired properties.
Furthermore, two groups in either of the lists are isomorphic
if and only if the corresponding groups $L$ are isomorphic,
the corresponding integers $t$ are equal and the corresponding
automorphisms $x$ define, up to conjugation, the same element
$\Out (N) = \Aut (N) / \Int (N)$.
\end{thm}

Before proving the Theorem, we observe that:
\begin{ite}
\setlength{\itemsep}{8pt}
\item If $G$ is a group which is not abelian and such that $G/H$ is abelian
	for any normal non-trivial subgroup $H$, then $G$ has a unique	
	minimal non-trivial normal subgroup, namely $[G,G]$.
\item If $G$ is a finite group such that $G/N$ is cyclic for any normal
	non-trivial subgroup $N$, then ``$Z(G)$ is trivial" is equivalent
	to ``$G$ is non-abelian". 
\item Case (\ref{clasigurii}) covers the case where $G$
	is non-abelian simple ($t=1$, $C$ is trivial).
\item In case (\ref{clasiguri}), identify $L$ with $\fp$ and the automorphism
	$x$ with $T\in \GL(t, \fp)$. Then $x$ acts irreducibly if and only if
	the characteristic polynomial of $T$ is irreducible, hence equals
	the minimal polynomial of $T$. In this case, if $n = \ord x$ and if
	$d$ divides $n$, $1 \neq d \neq n$, then $\ker(T^d-\id)=0$; this
	implies that $N - 0$ is a union of copies of $C$. Hence $|C|$ divides
	$p^t-1$. Clearly, we may assume that $x$ acts by the companion matrix
	of an irreducible polynomial in $\fp[X]$ of degree $t$.
\item Let $N$, $C$ be two groups, $C$ acting on $N$ by group automorphisms. 
	The center $Z(N \rtimes C)$ is given by
	$Z(N\rtimes C) = \{(n,c)\ |\ 
		c\in Z(C),\ n\in N^C,\ c\cdot m=n^{-1}mn\ \forall m\in N\}$.
	In particular, if  $p:(N\rtimes C)\to (N\rtimes C)/Z(N\rtimes C)$
	is the projection then   $p(N) \simeq N  / Z(N)^C$.
\end{ite}

\pf 
\emph{Step I.} We first show that the groups described in the Theorem have
the desired properties. 

Let $L$, $N$, $C$, $G$ be as in case (\ref{clasiguri}).
By the irreducibility of the action of $C$, being $x$
non-trivial, we see that $Z(G)$ is trivial.
We claim that any non-trivial normal subgroup $M$
of $G$ contains $N$.  For, $M\cap N$ is either trivial or $M\cap N=N$. If
$a\in G$ and $m\in M$, $m\neq e$; then $[a,m]\in M\cap [G,G]\subseteq M\cap N$.
Thus $M \cap N$ is non-trivial, since otherwise $m \in Z(G)$, proving the claim.
Hence any non-trivial quotient of $G$, being a quotient of $C$, is cyclic,
and $G$ satisfies the requirements of the theorem.

Let now $L$, $N$, $C$, $G$ be as in case (\ref{clasigurii}).  We identify
$N$ with its image in $G$. 
We claim next that $Z(G)$ is trivial. Let $(n,c) \in N\rtimes C$
be such that $p(n,c)=nc\in Z(G)$. It is easy to see that 
$n\in N^x$ and $c$ acts on $N$ by conjugation by $n^{-1}$.
Thus $(n,c)\in Z(N\rtimes C)$ and $nc=1$ in $G$.

Any normal subgroup $P$ of $N$ is of the form $\prod_{j\in J} L_j$,
for some subset $J$ of $\{1, \dots, t\}$; if $P$ is also $x$-stable
then either $P$ is trivial or equals $N$, because $x$ permutes the
copies of $L$. As in case (\ref{clasiguri}), we conclude that any non-trivial
normal subgroup $M$ of $G$ contains $N$; hence any non-trivial
quotient of $G$ is cyclic.

\medskip
\emph{Step II.} Let $G$ be a finite group with a minimal normal
non-trivial subgroup $N$, and assume that $G/N$ is cyclic.
Then there exists a simple subgroup
$L$ of $N$, and a subgroup $C = <x>$ of $G$ such that
$N = L \times \dots \times L$ ($t$ copies) and $G = NC$.

Indeed, let $x\in G$ be such that its class generates $G/N$ and
let $C$ be the subgroup generated by $x$; then $G = NC$. Let $L$
be a minimal normal non-trivial subgroup of $N$.
Then $L_i : = x^{i-1}Lx^{-i+1}$ is also a minimal normal subgroup
of $N$ and $(L_1\cdots L_j) \cap L_{j + 1}$ is either trivial or $L_{j+1}$.
Let $t$ be the smallest positive integer such that
$(L_1\cdots L_t) \cap L_{t+1}=L_{t+1}$.
Then $L_1 \cdots L_{t} \simeq L_1 \times \dots \times  L_{t}$ is
a normal subgroup of $N$ and is stable by conjugation by $x$;
so it is normal in $G$ and therefore equal to $N$. If $S$ is a normal
subgroup of $L$, then $S \simeq S \times 1 \times \dots \times 1$
is normal in $N$, and by minimality of $L$, $L$ is simple.

\medskip
\emph{Step III.}
We now show that any group $G$ satisfying the requirements of the
theorem is either as in case (\ref{clasiguri}) or (\ref{clasigurii}).
We may assume that $G$ is not simple. We keep the notation
of Step II and assume then that $Z(G)$ is trivial and that any proper quotient of
$G$ is cyclic.

If $N$ is abelian then $N \cap C \subseteq Z(G)$ is trivial,
whence $G=N\rtimes C$. Furthermore, $x$ should act irreducibly
since any subgroup $P\subset N$ which is $x$-stable is normal.
Hence, $G$ is as in case (\ref{clasiguri}).

If $L$ is not abelian and $t>1$ we have, for $i>1$,
$[xL_tx^{-1},L_i]=x[L_t,L_{i-1}]x^{-1}=1$, from where $x$ sends
$L_t$ isomorphically to $L_1$, and $x$ acts as in case
(\ref{clasigurii}) of the statement.
Consider the projection $p:N\rtimes C \to G$. Since $Z(G)$ is trivial,
$Z(N\rtimes C)\subset\ker(p)$. Let now $(n,c)\in\ker(p)$; then
$c=n^{-1}$ in $G$, whence $c$ acts by conjugation on $N$ by $n^{-1}$.
Thus $(n,c)\in Z(N\rtimes C)$, and we are done.

\medskip
\emph{Step IV.}
We prove the uniqueness statement.
$N$ is unique since, as remarked, $N=[G,G]$. Now, $L$ is unique by Jordan--H\"older,
then $t$ is unique. Since $x$ was chosen modulo $N$, $x$ and $x'$ give rise
to the same group if they coincide in $\Out(N)$. Furthermore, since any
automorphism $G\to G$ must leave $N$ invariant, $x$ is unique up to conjugation
in $\Out(N)$.
\epf

\begin{rem}\label{remarkorbitas}
Let $N$, $C$ be finite groups with $C$ acting on $N$ by group automorphisms, 
and let $G = N \rtimes C$. If $(m,z), (n,y) \in G$, then 
$$(m, z)(n,y)(m,z)^{-1} = (m(z\cdot n)(z y z^{-1}\cdot m^{-1}), z y z^{-1}).$$ 
When $C$ is abelian, it follows that
\begin{equation}\label{clasedeconj}
{\mathcal C}(n,y) = \bigcup_{z\in C} {\mathcal O_y}(z\cdot n)\times\{y\},
\end{equation}
where ${\mathcal C}$ stands for conjugacy class, and ${\mathcal O_y}$
for the orbit under the action of $N$ on itself given by 
\begin{equation}\label{acciontorcida}
m\rightharpoonup_y n := m \, n \, (y\cdot m^{-1}). 
\end{equation} 
Note that $\bigcup_{z\in <y>} {\mathcal O_y}(z\cdot n) = {\mathcal O_y}(n)$.
For, $n^{-1}\rightharpoonup_y n =  y\cdot n$, and the claim follows.
Note also that $m\rightharpoonup_y n$ is \emph{not} the same as 
$m\trid n$.
\end{rem}

We can now state the classification of simple racks. We begin by the
following important theorem. The proof uses
\cite[Lemma 8]{egs}, which is in turn based on the classification
of simple finite groups.

\begin{thm}\label{clasipotprimo} 
Let $(X, \trid)$ be a simple crossed set and let $p$ be a prime number.
Then the following are equivalent.

\begin{enumerate}
\item\label{tuta1} $X$ has $p^t$ elements, for some $t\in \N$.
\item\label{tuta2} $\gax$ is solvable.
\item\label{tuta3} $\gax$ is as in case (\ref{clasiguri}) of Theorem \ref{clasigur2}.
\item\label{tuta4} $X$ is an affine crossed set $({\fp}^t, T)$ where $T\in
	\GL(t, \fp)$ acts irreducibly.
\end{enumerate}
\end{thm}

\pf
\begin{list}{}{}
\item[(\ref{tuta1} $\implies$ \ref{tuta2})]
This is \cite[Lemma 8]{egs}.
\item[(\ref{tuta2} $\implies$ \ref{tuta3})]
This is Proposition \ref{simplesprimerresultado}
plus Theorem \ref{clasigur2}.
\item[(\ref{tuta3} $\implies$ \ref{tuta4})]
This follows from the preceding discussion.
Since $\phi(X)\subset\Inn_\trid(X)$ is a conjugacy class, by \eqref{clasedeconj}
we have $\phi(X)=N\times\{x^r\}$ for some $r$. Since $\phi(X)$ generates
$\Inn_\trid(X)$, $r$ must be coprime to the order of $x$. Take $y=x^r$
and call $T\in\GL(t,\fp)$ the action of $y$ (which is also irreducible). We have
$$(m,y)\trid(n,y)=(m,y)(n,y)(m,y)^{-1}=(Tn+(\id-T)m,y).$$
%we have $(m, T^d)(n,T)(m, T^d)^{-1} = (T^d\cdot n + (\id - T)\cdot m, T)$,
%thus $m\trid n = T\cdot n + (\id - T)\cdot m$, that is, $X$ is affine.
\item[(\ref{tuta4} $\implies$ \ref{tuta1})]
Clear.
\end{list}
\epf

\begin{cor}\label{clasirackpot}
The classification of simple racks with  
$p^t$ elements, for some  prime number $p$ and $t\in \N$, is the following:
\begin{enumerate}\itemsep 6pt
\item\label{clrkp1}
Affine crossed sets $({\fp}^t, T)$, where $T$ is the companion matrix of an 
irreducible monic  polynomial in $\fp[X]$ different from $X-1$ and $X$.
\item The permutation rack corresponding to the cycle $(1, 2, \dots, p)$ if $t=1$.
\hfill\qed
\end{enumerate}
\end{cor}

\begin{rem}
Keep the notation of Step III in Theorem \ref{clasigur2}. If $L$ is abelian,
then $x$ does not necessarily send $L_t$ to $L_1$, as wrongly stated in
\cite[Lemma 4 (ii)]{jo2}. This explains why in \cite[Lemma 6]{jo2},
only irreducible polynomials of the form $X^t - a$ appear; while,
as we have seen, this restriction is not necessary.
\end{rem}

\medskip
\begin{thm}\label{clasipqr}
Let $(X, \trid)$ be a  crossed set whose cardinality
is divisible by at least two different primes.
Then the following are equivalent.

\begin{enumerate}\itemsep 6pt
\item\label{ququ1}
$X$ is simple.
\item\label{ququ2}
There exist a non-abelian simple group $L$, a positive integer $t$ 
and  $x \in \Aut (L^t)$, where $x$ acts by \eqref{acciondex2} 
for some $\theta\in\Aut(L)$, such that
$X = {\mathcal O_x}(n)$ is an orbit of the action $\rightharpoonup_x$
of $N=L^t$ on itself as in \eqref{acciontorcida} ($n \neq m^{-1}$ if $t=1$
and $x$ is inner, $x(p) = mpm^{-1}$). Furthermore, $L$ and $t$ are unique,
and $x$ only depends on its conjugacy class in $\Out (L^t)$.
If $m, p\in X$ then
\begin{equation}\label{tridsimplenoab}
m \trid p = m  x(pm^{-1}).
\end{equation}
\end{enumerate}
\end{thm}

\pf
\begin{list}{}{}
\item[(\ref{ququ1} $\implies$ \ref{ququ2})]
By Theorem \ref{clasipotprimo}, $\gax$ is as in case (\ref{clasigurii})
of Theorem \ref{clasigur2}. Therefore, we have $L$, $t$ and $x\in\Aut(L^t)$. 
Let $G=  (N\rtimes C)/Z(N\rtimes C) \simeq \gax$, and
$\widetilde G= N\rtimes C$ and let $p: \widetilde G\to G$ be the projection.
If ${\mathcal C}(n,y)$ is a conjugacy class in $\widetilde G$ then
$p({\mathcal C}(n,y))$ is a conjugacy class in $G$, and it is not difficult
to see that $p:{\mathcal C}(n,y)\to p({\mathcal C}(n,y))$ is an isomorphism
of crossed sets. Then $X$ has the structure of $\mathcal{C}(n,y)$ given by
\eqref{acciontorcida}. Now, $ny\in G$ must be such that $p(\mathcal{C}(n,y))$
generates $G$. Since the subgroup generated by $p(\mathcal{C}(n,y))$ is
invariant, we know by the proof of Theorem \ref{clasigur2} that it is
either trivial or it contains $N$. It is trivial if $(n,y)\in Z(N\rtimes C)$,
i.e., if $t=1$ and $y$ acts on $N$ by conjugation by $n^{-1}$; thus we must
exclude this case. This case excluded, $y$ must generate $C$, and then we can
take $x=y$ in the proof of Theorem \ref{clasigur2}, whence $X$ is as in
\eqref{tridsimplenoab}.
%On the other hand, since the subgroup
%of $G$ generated by $p({\mathcal C}(n,y))$ is normal, it generates
%$G$ if and only if $y$ generates $C$ (when $t> 1$, or $t=1$ and $x$ is not inner)
%or $n\neq 1$ (when $t=1$ and $x = \id$). We can assume that
%$X \simeq  {\mathcal C}(n,x)$, without loss of generality (just
%take $x=y$ in the proof of Theorem \ref{clasigur2}).
%We conclude by Remark \ref{remarkorbitas}; \eqref{tridsimplenoab}
%follows since ${\mathcal C}(n,x)$ is standard. 
\item[(\ref{ququ2} $\implies$ \ref{ququ1})]
Similar. 
\end{list}
\epf

We restate the previous results as: if $X$ is a simple rack then either
\begin{enumerate}
\item $|X|=p$ a prime, $X\simeq\fp$ a permutation rack, $x\trid y=y+1$.
\item $|X|=p^t$, $X\simeq({\fp}^t, T)$ is affine, as in 
	Corollary \ref{clasirackpot} (\ref{clrkp1}).\label{rprp2}
\item $|X|$ is divisible by at least two different primes, and $X$ is
	twisted homogeneous, as in Theorem \ref{clasipqr}.
\end{enumerate}
Compare this with \cite[Thm. 7]{jo2}.

The simple crossed sets in (\ref{rprp2}) can be alternatively 
described as $(X, \trid^a)$ where $X\simeq \fq$, $q= p^t$;
$a\in \fq$ generates $\fq$ over $\fp$ and 
$x\trid^a y=(1-a)x + ay$. 
It follows easily that $\Aut (X, \trid^a)$ is the
semidirect product $\fq \rtimes \fq^{\times}$.
 
It is natural to ask how many different simple crossed sets with
$q= p^t$ elements there are. This is a well-known elementary
result. For, if $I(n)$ denotes the 
number of monic  irreducible polynomials  in $\fp[X]$ with
degree $n$, then $\sum_{d \vert n} d I(d) = p^n$. Thus
$$I(n)=\frac{1}{n} \sum_{d|n} \mu\left(\frac nd\right)p^d,$$
where $\mu$ is the M\"obius function.

%%%%%%%%%%%%%%%%%%%%%%%%%%%%%%%%%%%%%%%%%%%%%%%%%%%%%%%%%%%%%%%
%%%%%%%%%%%%%%%%%%% SubSeccion Cohomology %%%%%%%%%%%%%%%%%%%%%
%%%%%%%%%%%%%%%%%%%%%%%%%%%%%%%%%%%%%%%%%%%%%%%%%%%%%%%%%%%%%%%
\section{Cohomology}\label{sn:4}

\subsection{Abelian cohomology}

Let $(X,\trid)$ be a rack. We define now a cohomology theory which
contains all cohomology theories of racks known so far.
We think that this cohomology can be computed by some cohomology theory
in the category of modules over $X$.

For a sequence of elements $(x_1,x_2,\ldots x_n)\in X^n$ we will denote
$$[x_1\cdots x_n]=x_1\trid(x_2\trid(\cdots (x_{n-1}\trid x_n)\cdots)).$$
Notice that if $i<n$ then
$$[x_1\cdots x_i] \trid[x_1\cdots\widehat{x_i}\cdots x_n]
	=[x_1\cdots x_n].$$
\begin{defn}\label{df:ddc}%definicion del complejo
Let $*\in X$ be a fixed element (which is important only in degrees $0$ and $1$).
Let $\ral X$ be the rack algebra of $X$ (see Definition \ref{defrackalg}) and let,
for $n\ge 0$, $C_n(X)=\ral XX^n$,
i.e., the free left $\ral X$-module with basis $X^n$ ($X^0=\{*\}$ is a singleton).
For $n\ge 1$, let $\partial=\partial_n:C_{n+1}(X)\to C_n(X)$ be the $\ral X$-linear map
defined on the basis by
\begin{equation}\label{dfhomgral}
\begin{split}
\text{For $n\ge 1$:}\qquad
&\partial (x_1,\ldots,x_{n+1})
        = \sum_{i=1}^n(-1)^i\eta_{[x_1\cdots x_i],[x_1\cdots\widehat{x_i}\cdots x_{n+1}]}
                (x_1,\ldots,\widehat{x_i},\ldots,x_{n+1}) \\
&\hspace{3cm}-\sum_{i=1}^n(-1)^i(x_1,\ldots,x_{i-1},x_i\trid x_{i+1},\ldots,x_i\trid x_{n+1}) \\
&\hspace{3cm}-(-1)^{n+1}\tau_{[x_1\cdots x_n],[x_1\cdots x_{n-1}x_{n+1}]}(x_1,\ldots,x_n). \\
\text{For $n=0$:}\qquad
&\partial(x)=-\tau_{*,*^{-1}\trid x}*.
\end{split}
\end{equation}
\end{defn}

\begin{lem}
$(C_{\bullet}(X), \partial)$ is a complex.
\end{lem}
\pf
We decompose $\partial_n=\sum_{i=1}^{n+1}(-1)^i\partial_n^i$, where
\begin{align*}
\partial_n^i(x_1,\ldots,x_{n+1})
        &= \eta_{[x_1\cdots x_i],[x_1\cdots\widehat{x_i}\cdots x_{n+1}]}
                (x_1,\ldots,\widehat{x_i},\ldots,x_{n+1}) \\
        &\hspace{1cm}-(x_1,\ldots,x_{i-1},x_i\trid x_{i+1},\ldots,x_i\trid x_{n+1})
                \hspace{1cm}\mbox{for $i\le n$} \\
\partial_n^{n+1}(x_1,\ldots,x_{n+1})
        &=-\tau_{[x_1\cdots x_n],[x_1\cdots x_{n-1}x_{n+1}]}(x_1,\ldots,x_n),
                \hspace{1cm}\mbox{for $i=n+1>1$,} \\
\partial_0^1(x) &=-\tau_{*,*^{-1}\trid x}*.
                \hspace{5cm}\mbox{for $n=0$.}
\end{align*}
Then, it is straightforward to verify that
$$\partial_{n-1}^i\partial_n^j=\partial_{n-1}^{j-1}\partial_n^i
        \qquad\mbox{for $1\le i<j\le n+1$},$$
and thus
\begin{align*}
\partial_{n-1}\partial_n
	&=\hspace{-2mm}\sum_{1\le i<j\le n+1}\hspace{-2mm}(-1)^{i+j}\partial^i\partial^j
	+\hspace{-2mm}\sum_{1\le j\le i\le n}\hspace{-2mm}(-1)^{i+j}\partial^i\partial^j \\
	&=\hspace{-2mm}\sum_{1\le i<j\le n+1}\hspace{-2mm}(-1)^{i+j}\partial^{j-1}\partial^i
	+\hspace{-2mm}\sum_{1\le j\le i\le n}\hspace{-2mm}(-1)^{i+j}\partial^i\partial^j
	=0
\end{align*}
\epf

We are now in position to define rack (co)homology.

\begin{defn} Let $X$ be a rack.
Let $A=(A,\eta,\tau)$ be a left $X$-module and take $C^n(X,A)=\Hom_{\ral X}(C_n(X),A)$,
and the differential $d=\partial^*$. By the lemma, this is a cochain complex.
We define then
$$H^n(X,A)=H^n(C^{\bullet}(X,A)).$$
If $A$ is a right $X$-module (i.e., a right $\ral X$-module), we define
$$C_n(X,A)=A\otimes_{\ral X} C_n(X)\text{\ \ and\ \ }H_n(X,A)=H_n(C_\bullet(X,A)).$$
\end{defn}

\begin{rem} Low degree cohomology can be interpreted in terms of extensions:
a $2$-cocycle is the same as a function $\kappa:X\times X\to A$ which
satisfies \eqref{generaldoscociclo}; and two $2$-cocycles are cohomologous
if and only if they satisfy \eqref{eq:2cch}.
\end{rem}

\begin{rem} Let $X$ be a quandle; replacing the rack algebra $\ral X$
by the quandle algebra $\qal X$ in Definition \ref{df:ddc} leads to
a quandle cohomology theory that has as a particular case the
\emph{quandle cohomology} $H^\bullet_Q(X,A)$ in \cite{quandle1}.
\end{rem}

We consider now particular cases of this definition.

\subsection{Cohomology with coefficients in an abelian group}
Recall that $\ral X$ has an augmentation $\ral X\to\Z$ given by
$\eta_{i,j}\mapsto 1$, $\tau_{i,j}\mapsto 0$.
Then, any abelian group $A$ becomes an $X$-module. The complexes
$C_{\bullet}(X,A)$, $C^{\bullet}(X,A)$ coincide thus with previous
complexes found in the literature (see for instance
\cite{fr,quandle1,gr}). We recall them for later use:
\begin{equation}\label{eqcca}%complejo de cohomologia abeliana
\begin{split}
&C_n(X,A)=A \otimes_\Z\Z X^n,\qquad
	C^n(X,A)=\Hom_\Z(\Z X^n,A)\simeq\fun(X^n,A) \\
&\partial(a\otimes (x_1,\ldots,x_{n+1}))
        =\sum_{i=1}^n(-1)^i\big(a\otimes (x_1,\ldots,\widehat{x_i},\ldots,x_{n+1}) \\
        &\hspace{6cm}-(a\otimes (x_1,\ldots,x_{i-1},x_i\trid x_{i+1},
		\ldots,x_i\trid x_{n+1})\big) \\
&df(x_1,\ldots,x_{n+1})
        =\sum_{i=1}^n(-1)^i\big(f(x_1,\ldots,\widehat{x_i},\ldots,x_{n+1}) \\
        &\hspace{6cm}-f(x_1,\ldots,x_{i-1},x_i\trid x_{i+1},\ldots,x_i\trid x_{n+1})\big)
\end{split}
\end{equation}
Here, $\partial_0:C_1(X,A)\to C_0(X,A)$ vanishes.
Notice that $H^1(X, A)= A^{\pi_0(X)}$, where $\pi_0(X)$ is the set of
$\gax$-orbits in $X$.
\begin{exmp}\label{exp-cocycle}
Let $(X,\trid)$ be a crossed set,
let $A$ be an abelian group (denoted additively), and let $f$ be a $2$-cocycle
with values in $A$. Let $B: X\times X\to \Sim_A$ be given by
$$
B_{ij}(a) = a + f_{ij}, \qquad i,j\in X, \quad a\in A.
$$
Then $B$ is a constant rack cocycle (i.e., it satisfies \eqref{cocycle7}).
It is a constant quandle cocycle (i.e., it satisfies \eqref{cocycle5}) iff
$f_{ii}=0\ \forall i\in X$. The definition of ``quandle cocycle" is thus
seen to be the same for these kind of extensions as that in \cite{quandle1}.
The cocycle $B_{ij}$ satisfies \eqref{cocycle6} iff $f_{ji} = 0$
whenever $i\trid j =j$ and $f_{ij} = 0$.
\end{exmp}

\begin{lem} $H^2(X,G)\simeq\Hom(H_2(X,\Z),G)$ for any abelian group $G$.
\end{lem}
\pf This follows from the ``Universal Coefficient Theorem"
since $H_1(X,\Z)$ is free (cf. \cite[Prop. 3.4]{cjks}).
\epf

\begin{rem}%\label{rmmop}
Some second and third cohomology groups are computed in \cite{ln}; some others
in \cite{mo}. See \cite{sur} for tables of the computations done so far.
In particular, we excerpt from \cite{mo} that $H^2(X,\Z)=\Z$ if $X=(\Z/p,\trid^q)$
where $p$ is a prime and $1\neq q\in(\Z/p)^\times$.
\end{rem}

\begin{lem} Let $X$ be the disjoint sum of the indecomposable crossed sets $Y$ and $Z$.
Then $H^2(X,\Z) \simeq H^2(Y,\Z) \oplus H^2(Z,\Z)\oplus \Z^2$. \end{lem}
\pf Let $f \in Z^2(X,\Z)$, i.e.,
$f_{j,k}+f_{i,j\trid k}=f_{i\trid j,i\trid k}+f_{i,k}$ $\forall i,j,k\in X$.
Then one has $f_{i, j\trid k} = f_{i, k} = f_{t\trid i, k}$ for all $i, t\in Y$,
$j,k\in Z$, since the actions of $Z$ on $Y$ and viceversa are trivial.
Therefore $f_{i, k} = f_{t, j}$ for all  $i, t\in Y$, $j,k\in Z$,
being both $Y$ and $Z$ indecomposable. We conclude that
$Z^2(X,\Z) \simeq Z^2(Y,\Z) \oplus Z^2(Z,\Z)\oplus \Z^2$.
Now, if $g: X\to \Z$ then $\delta g$ is really only a function on
$(Y\times Y) \cup (Z\times Z)$, and the claim follows. \epf

\subsection{Restricted modules}
Consider a restricted $X$-module $A$, see Definition \ref{def:restmod}. This is the
same as a $G_X$-module, $G_X$ the enveloping group of $X$. We get then the
complex $(C^{\bullet}(X,A),d)\simeq(\fun(X^{\bullet},A),d)$, with
\begin{align*}
df(x_1,\ldots,x_{n+1})
        =&\sum_{i=1}^n(-1)^i[x_1\cdots x_i]f(x_1,\ldots,\widehat{x_i},\ldots,x_{n+1}) \\
        &-\sum_{i=1}^n(-1)^if(x_1,\ldots,x_{i-1},x_i\trid x_{i+1},\ldots,x_i\trid x_{n+1}) \\
        &-(-1)^{n+1}(1-[x_1\cdots x_{n+1}])f(x_1,\ldots,x_n)
\end{align*}
As a particular case, suppose that $X$ is any quandle and $(A,g)$ is an affine crossed
set. Take $\Lambda=\Z[T,T^{-1}]$ the ring of Laurent polynomials. Then $A$ becomes a
$\Lambda$-module, and a fortiori a $G_X$ module by $x\cdot a=Ta=g(a)$ $\forall x\in X,\ a\in A$
(since for any quandle there is a unique algebra map $\Z G_X\to\Lambda$, $x\mapsto T$).
Then $\eta_{i,j}$ acts by $g$ and $\tau_{i,j}$ acts by $f=1-g$ on $A$. We get in this way
the complex considered in \cite{ces}.

\begin{lem}\label{etingof} If $X$ is a rack and 
$A$ is an abelian group with trivial action then
$$H^2(X,A) = H^1(G_X, \fun(X,A)).$$
Here the space $\fun(X,A)$ of all functions
from $X$ to $A$ is a trivial left $G_X$-module;
the right action is given by $(f\cdot x)(z)=f(x \trid z)$.
\end{lem}
\pf (Sketch, see \cite{eg}.)
Consider the map $\Phi:X\to G_X$, $x\mapsto x$.
Let $f:G_X\to\fun(X,A)$ and take $r(f):X\times X\to A$,
where $r(f)(x,y)=f(\Phi(x))(y)$. It is easy to see that this map
gives a morphism $H^1(G_X,\fun(X,A))\to H^2(X,A)$. On the other
hand, let $g\in Z^2(X,A)$. This gives a map $g':X\to\fun(X,A)$,
$g'(x)(y)=g(x,y)$. Recall that for a right
$G_X$-module $M$, a map $\pi:G_X\to M$ is a $1$-cocycle
iff the map $\hat \pi:G_X\to G_X\ltimes M$ given by 
$z\mapsto (z,\pi(z))$ is a homomorphism of groups.
Denote $M=\fun(X,A)$. We have a map $\xi_g: X\to G_X\ltimes M$ given by 
$\xi_g(x)=(x,g'(x))$. So we need to show that $\xi_g$ extends to a 
homomorphism $G_X\to G_X\ltimes M$. But the group $G_X$ is generated by 
$X$ with relations $xy=(x\trid y)x$. Thus, we only need to check that the
$\xi_g(x)$'s satisfy the same relations, and it is straightforward to see
that this is equivalent to $dg=0$.
Now, it is easy to see that this map is the inverse of $r$.
\epf

%%%%%%%%%%%%%%%%%%%%%%%%%%%%%%%%%%%%%%%%%%%%%%%%%%%%%%%%%%%%%%%
%%%%%%%%%%%%%%%% SubSeccion non principal Cohomology %%%%%%%%%%
%%%%%%%%%%%%%%%%%%%%%%%%%%%%%%%%%%%%%%%%%%%%%%%%%%%%%%%%%%%%%%%

\subsection{Non-principal cohomology}\label{ssnpc}

Let $X=\sqcup_{i\in I}X_i$ be a decomposition of the rack $X$.
It is possible then to decompose the complex $C_{\bullet}(X,\Z)$
of \eqref{eqcca} into a direct sum
$$C_n(X)=\bigoplus_{i\in I}C_n^i(X),\qquad
C_n^i(X)=\Z(X^{n-1}\times X_i)\simeq\Z X^{n-1}\otimes\Z X_i.$$
For each $i\in I$, let $A_i$ be an abelian group. We denote by $A_I$ this collection.
Then we take $C^{\bullet}(X,A_I)$ the complex
$$C^{\bullet}(X,A_I)=\bigoplus_{i\in I}\Hom_\Z(C_n^i(X),A_i)$$
Notice that if $A=A_i$ $\forall i\in I$, then this complex is the same as
$C^{\bullet}(X,A)$ in \eqref{eqcca}.

%%%%%%%%%%%%%%%%%%%%%%%%%%%%%%%%%%%%%%%%%%%%%%%%%%%%%%%%%%%%%%%
%%%%%%%%%%%%%%%% SubSeccion  non ab. Cohomology %%%%%%%%%%%%%%%
%%%%%%%%%%%%%%%%%%%%%%%%%%%%%%%%%%%%%%%%%%%%%%%%%%%%%%%%%%%%%%%

\subsection{Non-abelian cohomology}

Let $(X,\trid)$ be a crossed set and let $\Gamma$ be a group.
We define:
\begin{align}
H^{1}(X, \Gamma) &= Z^{1}(X, \Gamma)
= \{\gamma: X\to \Gamma\ |\ \gamma_{i\trid j} = \gamma_j,
\quad \forall i, j\in X\}, \\
Z^{2}(X, \Gamma) &= \{\beta: X\times X \to \Gamma\ |\ 
\be_{i,j\trid k}\, \be_{j,k} = \be_{i\trid j,i\trid k} \, \be_{i,k}
\quad \forall i, j, k\in X\}.
\end{align}
The elements of $Z^{2}(X, \Gamma)$ shall be called non-abelian $2$-cocycles with
coefficients in $\Gamma$.\newline
If $\beta, \widetilde{\beta}: X\times X \to \Gamma$ we set
$\beta \sim \widetilde{\beta}$ if and only if there exists
$\gamma: X\to \Gamma$ such that
\begin{equation}
\widetilde{\beta}_{ij} = \left(\gamma_{i\trid
j}\right)^{-1}\beta_{ij}\gamma_{j}.
\end{equation}
It is easy to see that $\sim$ is an equivalence relation, and that
$Z^{2}(X, \Gamma)$ is stable under $\sim$.
We define
\begin{equation}\label{eq:ddcnc}%def 2coc no conmutativo
H^{2}(X, \Gamma) = Z^{2}(X, \Gamma)/\sim.
\end{equation}

\begin{exmp} Let $S$ be a non-empty set; then $H^{2}(X, \Sim_S)$ parameterizes
isomorphism classes of constant extensions of $X$ by $S$, as in Subsection \ref{ssecc}.
\end{exmp}
\begin{rem}
Though obvious, we point out that $H^{2}(X, \Sim_S) = H^{2}(X, \Z_2)$
when $S$ has only two elements.
\end{rem}

%This remark talks about quandle cohomology; not defined here.
%\begin{rem}
%Similar to the abelian case (see \ref{rmmop}), it is proved in \cite{grp2}
%that if $(X,\trid)=(\Z/p, \trid^q)$, $p$ a prime and $1\neq q\in\Z/p^*$ then
%$H^2_Q(X,G)$ is trivial for any group $G$.
%\end{rem}

%%%%%%%%%%%%%%%%%%%%%%%%%%%%%%%%%%%%%%%%%%%%%%%%%%%%%%%%%%%%%%%
%%%%%%%%%%%%%%%% SubSeccion | non principal - non ab. %%%%%%%%%
%%%%%%%%%%%%%%%%%%%%%%%%%%%%%%%%%%%%%%%%%%%%%%%%%%%%%%%%%%%%%%%

\subsection{Nonabelian non-principal cohomology}

We combine the theory of non-principal cohomology in Subsection \ref{ssnpc}
with that of non-abelian cohomology: let $X=\sqcup_{i\in I}X_i$ be
a decomposition of the rack $X$ and for each $i\in I$ let $\Gamma_i$
be a group. Let us denote $\Gamma_I$ this collection. We consider
\begin{multline*}
Z^2(X,\Gamma_I)=\{f=\sqcup_if_i:X\times X_i\to \Gamma_i\ |\  \\
	f_i(x,y\trid z)f_i(y,z)=f_i(x\trid y,x\trid z)f_i(x,z)\ 
	\forall x,y\in X, z\in X_i\}.
\end{multline*}
As usual, if $f,\widetilde f\in Z^2(X,\Gamma_I)$, we say that $f\sim \widetilde f$
iff $\exists g=\sqcup_ig_i:X_i\to\Gamma_i$ such that
$$\widetilde f_i(x,y)=g_i(x\trid y)^{-1}f_i(x,y)g_i(y)\qquad\forall x\in X,y\in X_i.$$

The importance of such a theory becomes apparent in theorem \ref{thydnp}.
\begin{defn}
For $i\in I$, let be given a positive integer $n_i$ and a subgroup
$\Gamma_i\subset\GL(\kk,n_i)$. Let $f=\sqcup_if_i:X\times X_i\to\Gamma_i$.
Take $V=\bigoplus_{i\in I}X_i\times\kk^{n_i}$ a vector space and consider
the linear isomorphism
$$c^f:V\otimes V\to V\otimes V,\qquad
	c^f((x,a)\otimes (y,b))=(x\trid y,f_i(x,y)(b))\otimes(x,a),$$
	
$x\in X_j$, $y\in X_i$, $a\in\kk^{n_j}$, $b\in\kk^{n_i}$.
\end{defn}

\begin{thm}\label{thydnp}%yetter-drinfeld y no-principales.
\begin{enu}
\item Let $X,\Gamma_I,V,c^f$ be as in the definition.
Then $c^f$ satisfies the Braid Equation if and only if $f\in Z^2(X,\Gamma_I)$.
\item Furthermore, if $f\in Z^2(X,\Gamma_I)$, then there exists a group
$G$ such that $V$ is a Yetter--Drinfeld module over $G$. In particular, the braiding
of $V$ as an object in $\ydskg$ coincides with $c^f$. If $\Gamma_i$ is
finite $\forall i\in I$ and $X$ is finite, then $G$ can be chosen to be finite.
\item Conversely, if $G$ is a finite group and $V\in\ydskg$, then there
exist $X=\sqcup X_i$, finite groups $\Gamma_I$, $f\in Z^2(X,\Gamma_I)$ such that
$V$ is given as in the definition and the braiding $c\in\Aut(V\otimes V)$ in the
category $\ydskg$ coincides with $c^f$. Here, $X$ can be chosen to be a crossed set.
\item If $f\sim\widetilde f$ and $(V,c^f)$, $(\widetilde V,c^{\widetilde f})$
are the spaces associated to $f,\widetilde f$, then they are isomorphic as braided
vector spaces, i.e., there exists a linear isomorphism
$\gm:\tilde V\to V$ such that
$(\gm\otimes\gm)c^{\widetilde f}=c^f(\gm\otimes\gm):
	\widetilde V\otimes\widetilde V\to V\otimes V$.
\end{enu}
\end{thm}

\pf
\begin{enu}
\item Straightforward.
\item It follows from \cite[2.14]{grfrns}. The finiteness of $G$ follows from the
fact that the group $G\subset\GL(V)$ can be chosen to be the group generated by the
maps $(y,b)\mapsto(x\trid y,f(x,y)b)$, which is contained in the product
$\prod_{i\in I}G_i$, where $G_i=\Sim_{X_i}\times\Gamma_i$.
\item It follows from the structure of the modules in $\ydskg$. Indeed, if
$V=\oplus_iM(g_i,\rho_i)$ (see \cite{gr} for the notation), $g_i\in G$,
$G_i=\{xg_i=g_ix\}$ the centralizer of $g_i$, $\{h^i_1,\ldots,h^i_{s_i}\}$
a set of representatives of left coclasses $G/G_i$, and $t^{iu}_{jv}\in G_u$
defined by $h^i_jg_i(h^i_j)^{-1}h^u_v=h^u_{v'}t^{iu}_{jv}$; then take
$X=\sqcup_iX_i$, $X_i=\{h^i_1,\ldots,h^i_{s_i}\}$, and
$f(h^i_j,h^u_v)=\rho_u(t^{iu}_{jv})$.
\item It is straightforward to verify that if
$\widetilde f_i(x,y)=g_i(x\trid y)^{-1}f_i(x,y)g_i(y)$ $\forall x\in X,\ y\in X_i$,
then the map $(x,a)\mapsto (x,g_i(x)(a))\ (x\in X_i)$ is an isomorphism
of braided vector spaces.
\end{enu}
\epf

Notice that any Yetter--Drinfeld module over a group algebra can be constructed
by means of a crossed set, and one does not need the more general setting of
quandles, nor racks for it. However, racks may give easier presentations than
crossed sets for some braided vector spaces.

%%%%%%%%%%%%%%%%%%%%%%%%%%%%%%%%%%%%%%%%%%%%%%%%%%%%%%%%%%%%%%%
%%%%%%%%%% SubSeccion Braided vs of set-th. type %%%%%%%%%%%%%%
%%%%%%%%%%%%%%%%%%%%%%%%%%%%%%%%%%%%%%%%%%%%%%%%%%%%%%%%%%%%%%%
\section{Braided vector spaces}\label{sn:5}

We have seen that it is possible to build a braided vector space
$(\kk X, c^\qg)$ from a rack $(X, \trid)$  and a $2$-cocycle $\qg$,
\emph{cf.} Theorem \ref{thydnp} and \cite{gr}. It
turns out that the braided vector space does \emph{not} determine
the rack. We now present a systematic way of
constructing examples of different racks, with suitable cocycles, 
giving rise to  equivalent  braided vector spaces.
We consider affine modules over a rack $X$, that is
extensions of the form $X\times_{\kappa} A$, where $A$
is an abelian $X$-module; see Definition \ref{affinemodule-def}. 
If the cocycle $\qg$ is chosen in a convenient way, we can change
the basis ``\`a la Fourier" and obtain a braided vector space arising
from a set-theoretical solution of the QYBE. This solution is in turn
related to the braided vector space arising from the derived rack of
the set-theoretical solution of the QYBE. 

Subsection \ref{sub-setth} is an exposition
of the relevant facts about set-theoretical solutions needed in this paper.
In subsection \ref{sub-bvs-setth} we discuss 
braided vector spaces
arising from set-theoretical solutions.
In subsection \ref{subs-fourier} we present the general method
and discuss several examples.

%%%%%%%%%%%%%%%%%%%%%%%%%%%%%%%%%%%%%%%%%%%%%%%%%%%%%%%%%%%%%%%
%%%%%%%%%%%%%%%% SubSeccion set-theoretical %%%%%%%%%%%%%%%%%%%
%%%%%%%%%%%%%%%%%%%%%%%%%%%%%%%%%%%%%%%%%%%%%%%%%%%%%%%%%%%%%%%

\subsection{Set-theoretical solutions of the QYBE}\label{sub-setth}

There is a close relation between racks and set-theoretical
solutions of the Yang-Baxter equation, or, equivalently, of the
braid equation. It was already observed by Brieskorn \cite{Bk}
that racks provide solutions of the braid equation. On the other
hand, certain set-theoretical solutions of the braid equation produce
racks. This is proved in \cite{s,lyz1}, which belong to a series
of papers (see \cite{egs,ess,lyz2,lyz3}) devoted to set-theoretical
solutions of the braid equation and originated in a question by
Drinfeld \cite{Dr}. We give here the definitions necessary to us.

Let $X$ be a non-empty set and let $S: X\times X \to X\times X$ be
a bijection. We say that $S$ is a \emph{set-theoretical solution of the
braid equation} if
$(S\times \id)(\id \times S)(S\times \id)
	=(\id \times S)(S\times \id)(\id \times S)$.
We shall briefly say that $S$ is ``a solution" or that $(X, S)$ is a
\emph{braided set}. A trivial example of a solution
is the transposition $\tau: X\times X \to X\times X$, $(x,y)\mapsto(y,x)$.
It is well-known that $S$ is a solution if and only if
$R = \tau S: X\times X \to X\times X$ is a solution of the set-theoretical
quantum Yang-Baxter equation. If $(X, S)$ is a braided set, there is
an action of the braid group $\mathbb B_n$ on $X^n$, the standard
generators $\sigma_i$ acting by $S_{i, i+1}$, which means, as usual,
that $S$ acts on the $i,i+1$ entries.

In particular, a finite braided set gives rise to a finite quotient of
$\mathbb B_n$ for any $n$, namely
the image of the group homomorphism $\rho^n: \mathbb B_n \to 
\mathbb S_{X^n}$ induced by the action.

\begin{lem} \cite{Bk}. Let $X$ be a set and let $\trid: X \times X \to X$
be a function. Let
\begin{equation}\label{eq:calt}%como actua la trenza
c:X\times X\to X\times X,\quad c(i,j)=(i\trid j,i).
\end{equation}
Then $c$ is a solution if and only if $(X, \trid)$ is a rack.
\end{lem}

\pf It is easy to check that $c$ is a bijection if and only if
\eqref{ccc1} holds, and that it satisfies the braid equation
if and only if \eqref{ccc4} holds.
\epf

\begin{defn} Let $X$, $\widetilde X$ be two non-empty sets and let
$$S: X\times X \to X\times X,\quad
\widetilde S:\widetilde X\times\widetilde X\to\widetilde X\times\widetilde X$$
be two bijections. We say that $(X, S)$ and $(\widetilde X, \widetilde
S)$ are \emph{equivalent} if there exists a family of bijections
$T^n: X^n \to \widetilde X^n$ such that $T^n S_{i, i+1} =
\widetilde S_{i, i+1}T^n$, for all $n \ge 2$, $1\le i \le n-1$.
\end{defn}

If $(X, S)$ and $(\widetilde X, \widetilde S)$ are equivalent and
$(X, S)$ is a solution, then $(\widetilde X, \widetilde S)$ is
also a solution and the $T^n$'s intertwine the corresponding
actions of the braid group $\mathbb B_n$.

\begin{defn} \cite{s,lyz1}
Let $(X, S)$ be a solution and let $f, g: X \to \fun(X, X)$ be given by
\begin{equation}\label{trenzas-notac}
S(i,j) = (g_i(j), f_j (i)).
\end{equation}
The solution (or the braided set) is \emph{non-degenerate} if the images
of $f$ and $g$ lie inside $\Sim_X$.
\end{defn}

\begin{prop}\cite{s,lyz1}\label{slyz1}
Let $S$ be a non-degenerate solution with the notation in
\eqref{trenzas-notac} and define $\trid$ by
\begin{equation}\label{eq:rackasoc}
i\trid j = f_{i} \left(g_{f_j^{-1}(i)} (j)\right).
\end{equation}
\begin{enu}
\item One has
\begin{flalign}\label{rack-cs1}
&f \text{ preserves $\trid$, \emph{i.e.} } f_i(j\trid
k) = f_i(j) \trid f_i(k);& \\ \label{rack-cs2} & f_i f_j  =
f_{f_i(j)}f_{g_j(i)}, \qquad \forall i,j\in X.&
\end{flalign}
\item If $c$ is given by \eqref{eq:calt}, then $c$ is a solution;
we call it the \emph{derived solution of $S$}. The  solutions $S$ and
$c$ are equivalent, and $(X,\trid)$ is a rack.
\item Let $(X,\trid)$ be a rack and let $f:X \to \Sim_X$. We
define $g:X \to \Sim_X$ by 
\begin{equation}\label{tresveinticinco} 
g_i(j) = f_{f_j(i)}^{-1} (f_j(i) \trid j). 
\end{equation}
 Let $S: X\times X \to X\times X$ be given by
\eqref{trenzas-notac}. Then $S$ is a solution if and only if
\eqref{rack-cs1}, \eqref{rack-cs2} hold. If this happens, the
solutions $S$ and $c$ are equivalent, and $S$ is non-degenerate.
\end{enu}
\end{prop}

\pf
\begin{enu}
\item It is not difficult.
\item It is enough to show that $S$
and $c$ are equivalent; automatically, $c$ is a solution and {\it
a fortiori} $(X,\trid)$ is a rack. Let $T^n: X^n \to X^n$ be
defined inductively by $$ T^2(i, j) = (f_j(i), j), \quad T^{n+1}
= Q_n (T^n \times \id), $$ where $Q_n(i_1, \dots, i_{n+1}) =
(f_{i_{n+1}} (i_{1}), \dots, f_{i_{n+1}} (i_{n}), i_{n+1})$. One
verifies using \eqref{rack-cs1} and \eqref{rack-cs2} that $T^n S_{i, i+1} = c_{i,
i+1}T^n$, as needed.
\item Straightforward.
\end{enu}
\epf

Note that \eqref{tresveinticinco} is equivalent to 
\begin{equation}\label{tresveintiseis}
 g_{f_j^{-1}(h)}(j) = f^{-1}_h (h \trid j). 
\end{equation}

\begin{rem} Let $(X, S)$ be a non-degenerate solution and let
$\trid$ be defined by \eqref{eq:rackasoc}.
Then $(X,\trid)$ is a quandle if and only if
\begin{equation}
f_i^{-1}(i) = \left(g_{f_i^{-1}(i)} (i)\right), \qquad \forall
i\in X;
\end{equation}
if this holds, it is a crossed set if and only if
\begin{equation}
f_i^{-1}(j) = \left(g_{f_j^{-1}(i)} (j)\right) \quad \implies
\quad f_j^{-1}(i) = \left(g_{f_i^{-1}(j)} (i)\right)
	\quad\forall i,j\in X.
\end{equation}
\end{rem}

Let $(X, S)$, $(\widetilde X, \widetilde S)$ be two non-degenerate 
braided sets, with corresponding maps $f$, $g$,
resp. $\widetilde f$, $\widetilde g$. A function 
$\varphi: X\to \widetilde X$ is a morphism of braided sets if and only
if
\begin{align}\label{mor-cs1}
\widetilde g_{\varphi(i)}\varphi(j) &= \varphi g_{i}(j),
\\ \label{mor-cs2}
\widetilde f_{\varphi(i)}\varphi(j) &= \varphi f_{i}(j), \qquad \forall i,j\in X.
\end{align}

It can be shown that $\varphi$ is a morphism of braided sets if and only
$\varphi$ is a morphism of the associated racks and \eqref{mor-cs2} holds.
One may say that a non-degenerate braided set is \emph{simple}
if it admits no non-trivial projections. It follows that 
any solution associated to a simple crossed set is simple, but the converse
is not true as the following example shows:
take a set $X$ with $p$ elements, $p$ a prime, and a cycle $\mu$
of length $p$. Then $S(i, j) = (\mu(j), \mu^{-1}(i))$ is simple
but the associated rack is trivial.

\begin{defn} \cite{ess}.
Let $(X,S)$ be a solution and let $S^2(i,j) = (G_i(j), F_j (i))$.
The group $G_X$, resp. $A_X$, is the quotient of the free group
generated by $X$ by the relations $ij = g_i(j)f_j (i)$, resp. $f_j
(i) j = F_j (i) f_j (i)$, for all $i,j\in X$.

If $(X,\trid)$ is a rack and $c$ is the corresponding solution,
then $G_X$ is the group already defined in \ref{gx}, and coincides
with $A_X$.
\end{defn}

\subsection{Braided vector spaces of set-theoretical type}
\label{sub-bvs-setth}

We now describe how set-theoretical solutions of the QYBE plus a
$2$-cocycle give rise to braided vector spaces. We begin by the case
of solutions arising from a rack.

Let $(X, \trid)$ be a rack and let $q\in Z^2(X,\kk^{\times})$; so
that
\begin{equation}\label{cociclo-explicito}
q_{i,j\trid k} q_{j,k} = q_{i\trid j,i\trid k} q_{i,k} \quad
\forall i, j, k\in X.
\end{equation}
Then, by \ref{thydnp}, the space $V=\kk X$ has a structure of a
Yetter--Drinfeld module over a group whose braiding is given by
$c^q:\kk X\otimes\kk X\to\kk X\otimes\kk X$,
$$c^q(i\otimes j)=q_{i,j} i\trid j\otimes i, \qquad i,j\in X.$$

Let $X=\{x_1,\ldots,x_n\}$ be a set, let $S: X \times X \to X
\times X$ be a bijection and let $F: X\times X\to\kk^\times$ be a
function. Let $\kk X$ denote the vector space with basis $X$ and
define
\begin{equation}\label{eq:bvsf}
\begin{split}
&S^F: \kk X \otimes \kk X\to \kk X \otimes \kk X \\
&S^F(i \otimes j) = F_{i,j} S(i, j) = F_{i,j}\, g_i j \otimes f_j i, \qquad i,j \in X,
\end{split}
\end{equation}
where we use the notation in \eqref{trenzas-notac}.

\begin{lem}
\begin{enumerate}
\item $S^F$ is a solution of the braid equation if and only if
 $(X, S)$ a solution and
\begin{equation}\label{cociclo-braideq}
F_{i,j} F_{f_ji, k} F_{g_ij, g_{f_ji} k} = F_{j,k} F_{i, g_j k}
F_{f_{g_jk} i, f_kj}, \qquad i,j, k \in X.
\end{equation}
\item Assume that \eqref{cociclo-braideq} holds. Then $S^F$ is
\emph{rigid} if and only if $S$ is non-degenerate.
\end{enumerate}
\end{lem}

\pf
\begin{enu}
\item Straightforward.
\item Rigidity is equivalent to $c^{\flat}:V^*\otimes V\to V\otimes
V^*$ being an isomorphism, where
$$c^{\flat}=(\mbox{ev}_V\otimes\id_{V\otimes V^*})
(\id_{V^*}\otimes c\otimes\id_{V^*}) (\id_{V^*\otimes V}\otimes
\mbox{ev}^*_V).$$
Assume for simplicity that $F = 1$.
Let $(\delta_i)_{i\in X}$ be the basis of $V^*$ dual to $X$. 
Then $c^{\flat}(\delta_i \otimes j) = \sum_{h: g_j(h) = i} f_h(j) \otimes h$.
Hence, if $c^{\flat}$ is an isomorphism then $g_j$ is bijective for
all $j$, for $c^{\flat}(\delta_i \otimes j) = 0$ if $i$ is not in the image
of $g_j$. Now, if $f_h (j) = f_h (k)$  then 
$c^{\flat}(\delta_{g_j(h)} \otimes j) = c^{\flat}(\delta_{g_j(k)} \otimes j)$,
which implies $j = k$. 

Conversely, if $S$ is non-degenerate then $c^{\flat}$ is
an isomorphism with inverse
$(c^{\flat})^{-1}(r\otimes \delta_s)
	=\delta_{g_{f_s^{-1}(r)}(s)} \otimes f_s^{-1}(r)$. 
\end{enu}
\epf

\begin{defn}
Let $(X, S)$ be a non-degenerate solution and let $F: X\times X\to
\kk^\times$ be a function such that \eqref{cociclo-braideq} holds.
We say that the braided vector space $(\kk X, S^F)$ is \emph{of
set-theoretical type}.
\end{defn}

By results of Lyubashenko and others, 
a braided vector space $(\kk X, c^F)$ of
set-theoretical type can be realized as a Yetter-Drinfeld module
over some Hopf algebra $H$.
See for example \cite{Tk}.

\begin{exmp}\label{exa-quandlefromyd} Let $\Gamma$ be a finite group.
Let $x\in \Gamma$, let $\mathcal O$ be the conjugacy class
containing $x$ and let $\rho: \Gamma_x\to \Aut W$ be a finite
dimensional representation of $\Gamma_x$ the centralizer of $x$. We
choose a numeration $\{p_1 = x, p_2, \dots, p_r\}$ of $\mathcal O$
and fix elements $g_1, g_2, \dots, g_r$ in $\Gamma$ such that $g_i
x g_i^{-1} =p_i$. Then $$M(x, \rho) := \Ind_{\Gamma_x}^{\Gamma} W
\simeq\kk \mathcal O \otimes W \simeq \oplus_{1\le i \le r} g_i
\otimes W$$ is a Yetter-Drinfeld module over $\kk \Gamma$ with the
coaction $\delta(g_i \otimes w) = p_i \otimes g_i \otimes w$, and the
induced action; that is $h\cdot(g_i\otimes w)=g_j\otimes t\cdot w$,
where $j$ and $t\in \Gamma_x$ are uniquely determined by
$hp_ih^{-1} = p_j$, $hg_i = g_j t$. In particular, given $i, k\in
\{1, \dots, r\}$, let us denote by $j_{ik}\in \{1, \dots, r\}$,
$t_{ik}\in \Gamma_x$ the elements uniquely determined by
\begin{equation}
p_kp_ip_k^{-1} = p_{j_{ik}}, \qquad p_kg_i = g_{j_{ik}} t_{ik}.
\end{equation}
We can then express the braiding in a compact way; write for
simplicity $g_i w = g_i \otimes w$. If $u, w\in W$ and $i, k\in
\{1, \dots, r\}$ then
\begin{equation}
c(g_k w\otimes g_i u)=g_{j_{ik}} t_{ik}\cdot u\otimes g_k w.
\end{equation}
We know from Theorem \ref{thydnp} that this braided vector space
can be presented with the crossed set $\{p_1,\ldots,p_r\}$ and a
non-abelian $2$-cocycle with values in $\GL(W)$. We now show that
under a suitable assumption we can present it with a (larger) rack
and an abelian $2$-cocycle with values in $\kk^\times$:
assume that there exists a basis $w_1, w_2, \dots, w_R$ of
$W$ such that
\begin{equation}
h\cdot w_s = \chi_s(h) w_{\sigma_h(s)}.
\end{equation}
for some group homomorphism $\sigma: \Gamma_x \to \Sim_r$ and some
map $\chi: \{1, \dots, r\} \times \Gamma_x \to \kk^{\times}$
satisfying $$ \chi_s(th) = \chi_s(h)\chi_{\sigma_h(s)}(t), \qquad
1\le s \le r, \quad t,h \in \Gamma_x. $$ Then
\begin{equation}
c(p_k w_q \otimes p_i w_s) = \chi_s(t_{ik}) \; p_{j_{ik}}
w_{\sigma_{p_k}(s)} \otimes p_k w_q.
\end{equation}
That is, the braided vector space $(M(x, \rho), c)$ is of rack type.
\end{exmp}

We now introduce a relation between braided vector
spaces weaker than isomorphism but useful enough to deal with Nichols
algebras; for example braided vector spaces related by a twisting
are t-equivalent as below.

\begin{defn}\label{t-equivalence-def}
We say that two braided vector spaces $(V, c)$ and $(W, d)$
are \emph{t-equivalent} if there is
a collection of linear isomorphisms $U^n: V^{\otimes n} \to W^{\otimes n}$
intertwining the corresponding representations of the braid group $\mathbb B_n$,
for all $n\ge 2$. The collection $(U^n)_{n\ge 2}$ is called a \emph{t-equivalence}.
\end{defn}

\begin{exmp}
Let  $(\kk X, S^F)$ be a braided vector space of set-theoretical type
(see \eqref{eq:bvsf}). Let $(X, c)$ be the derived solution; define
$q_{ij} = F_{f_j^{-1}(i), j}$. If
\begin{equation}\label{t-eq-trivial}
q_{f_ki, f_kj} =q_{ij} \qquad \forall i,j,k\in X.
\end{equation}
then the collection of maps $T^n: X^n \to X^n$ defined in the proof of
Proposition \ref{slyz1} induce a t-equivalence between $(\kk X, S^F)$
and $(\kk X, c^q)$. Indeed, 
computing only the coefficients, we have
\begin{align*}
\textup{coeff}(T^n S^F_{h, h+1} ( i_1 \otimes \dots \otimes i_n))
& = F_{i_h, i_{h+1}}
%f_{i_n} f_{i_{n-1}} \dots f_{i_2}i_1 \otimes f_{i_n} f_{i_{n-1}}
%\dots f_{i_3}i_2 \dots \otimes f_{i_n} i_{n-1} \otimes i_n,
\\
\textup{coeff}(c^q_{h, h+1} T^n  ( i_1 \otimes \dots \otimes i_n))
& = q_{f_{i_n} f_{i_{n-1}} \dots f_{i_{h+1}}i_h,
	f_{i_n} f_{i_{n-1}} \dots f_{i_{h+2}}i_{h+1}}
%\\ & \qquad f_{i_n} f_{i_{n-1}} \dots f_{i_2}i_1 \otimes 
%f_{i_n} f_{i_{n-1}} \dots f_{i_3}i_2 \dots \otimes f_{i_n} i_{n-1}
%\otimes i_n;
\end{align*}
and the equality holds by \eqref{t-eq-trivial}.
\end{exmp}

%%%%%%%%%%%%%%%%%%%%%%%%%%%%%%%%%%%%%%%%%%%%%%%%%%%%%%
%%%%%%%%%%%%%%%% SubSeccion Fourier %%%%%%%%%%%%%%%%%%
%%%%%%%%%%%%%%%%%%%%%%%%%%%%%%%%%%%%%%%%%%%%%%%%%%%%%%

\subsection{Fourier transform}\label{subs-fourier}
We consider a rack $(X, \trid)$, a finite abelian $X$-module $A$,
a dynamical cocycle $\al: X\times X \to \fun(A\times A, A)$,
and a $2$-cocycle ${\mathfrak q}$ on the rack $X\times_{\al} A$.
Let $\widehat A$ be the group of characters of $A$.
We define a family of elements
\begin{equation}
(i, \psi) := \sum_{a\in A} \psi(a) (i, a) \in \kk(X\times_{\al}
A), \qquad i\in X, \quad \psi\in \widehat A.
\end{equation}

We want to know under which conditions there exists a family of
scalars $F_{ij}^{\psi, \phi}$ such that

\begin{equation}
c^{\mathfrak q}\left((i, \psi) \otimes (j, \phi) \right)= F_{ij}^{\psi, \phi}
(i\trid j, \vartheta) \otimes (i, \nu), \qquad i, j\in X, \quad
\psi, \phi\in \widehat A,
\end{equation}
for some $\vartheta, \nu \in \widehat A$. 
Our main result in this
direction, and one of the main results in this paper, is Theorem
\ref{mainfourier} below. 

In what follows, we shall assume that the extension $X\times_\al A$
is an affine module over $X$, \emph{cf.} Definition \ref{affinemodule-def}.
That is, $\al$ is given by
\begin{equation}\label{fourier-cociclodinamico}
\al_{i,j} (a,b) = \eta_{i,j} (b) + \tau_{i,j} (a) + \kappa_{ij},
\end{equation}
where $\eta_{i,j} \in \Aut(A)$, $\tau_{i,j} \in \End(A)$
define the $X$-module structure on $A$, and 
$\kappa_{ij} \in A$. We denote $Y_{\kappa} := X\times_{\kappa}A$.
We shall also write ${\mathfrak q}_{i,j}^{a, b} = 
{\mathfrak q}_{(i, a),(j, b)}$.
We begin by the following result.

\begin{lem}\label{auxfourier}
Let  ${\mathfrak q}$ be given by
\begin{equation}
\label{fourier-cociclo-condicion} {\mathfrak q}_{i,j}^{a, b} 
= \chi_{i,j} (b) \; \mu_{i,j} (a) \; q_{ij},
\end{equation}
where 
$\chi_{i,j}, \mu_{i,j} \in \widehat A$, $q_{ij} \in \kk^\times$.
Then ${\mathfrak q}$ is a  $2$-cocycle if and only if
\begin{align}
\label{fourier-cociclo-condicion1} 
\chi_{i,j \trid k} (\kappa_{j,k}) \; q_{i,j \trid k}\; q_{j, k} 
& = \chi_{i\trid j,i\trid k} (\kappa_{i,k}) \; 
\mu_{i\trid j,i\trid k} (\kappa_{i,j}) \; q_{i\trid j, i \trid k}\; q_{i, k}, \\
\label{fourier-cociclo-condicion2}
\chi_{i,j \trid k} (\eta_{j,k} (a)) \;  \chi_{j,k} (a) &=
\chi_{i\trid j,i\trid k} (\eta_{i,k} (a)) \;  \chi_{i,k} (a), \\
\label{fourier-cociclo-condicion3}
\mu_{i,j \trid k} (a)&= \chi_{i\trid j,i\trid k} 
(\tau_{i,k} (a)) \; \mu_{i\trid j,i\trid k} (\tau_{i,j} (a))
\;\mu_{i, k} (a), \\
\label{fourier-cociclo-condicion4}
\mu_{i\trid j,i\trid k} (\eta_{i,j} (a)) &= 
\chi_{i,j\trid k} (\tau_{j,k} (a)) \; \mu_{j,k} (a).
\end{align}
for all $i, j, k\in X$, $a \in A$. \end{lem}

\pf Writing explicitly down the cocycle condition on ${\mathfrak q}$,
one gets an equality of functions from $A\times A \times A$ to $\kk$.
Specialization at $(0,0,0)$ implies \eqref{fourier-cociclo-condicion1};
then, specialization at $(a,0,0)$,  $(0,a,0)$, $(0,0,a)$
implies the other conditions. The converse is similar.
\epf

\begin{thm}\label{mainfourier}
If the $2$-cocycle ${\mathfrak q}$ is given
by
\eqref{fourier-cociclo-condicion} with
$\chi_{i,j}, \mu_{i,j} \in \widehat A$, $q_{ij} \in \kk^\times$,
then
\begin{equation}
\label{fourier-formula} c^{\mathfrak q}\left((i, \psi) \otimes (j, \phi)
\right)= F_{ij}^{\psi, \phi} (i\trid j, \vartheta) \otimes (i,
\nu),
\end{equation}
for all $i, j\in X$, $\psi, \phi\in \widehat A$, where
\begin{align}
\label{fourier-cociclo1} F_{ij}^{\psi, \phi} &=  q_{ij}\;
\phi(\widetilde\kappa_{ij})^{-1} \; \chi_{i,j}(\widetilde
\kappa_{i,j})^{-1}, \\ \label{fourier-cociclo2}
 \vartheta &= (\phi \chi_{i,j}) \circ \eta_{i,j}^{-1}, \\
\label{fourier-cociclo3} \nu &=  \psi\; \mu_{i,j}  \;
((\phi\chi_{i,j})\circ \widetilde \tau_{i,j}) ^{-1}.
\end{align}
Here $\widetilde \tau_{i,j} (a) =
\eta_{i,j}^{-1}(\tau_{i,j} (a))$, $\widetilde \kappa_{ij} =
\eta_{i,j}^{-1}(\kappa_{ij})$.
\end{thm}

\pf We compute
\begin{align*}
c^{\mathfrak q}&\left((i, \psi) \otimes (j, \phi) \right) \\
&= \sum_{a, b \in A} \psi(a) \phi(b) {\mathfrak q}_{i,j}^{a, b}\ 
	(i \trid j, \eta_{i,j} (b)+\tau_{i,j} (a)+\kappa_{ij}) \otimes (i, a) \\
&= \sum_{a, c \in A} \psi(a) \phi(\eta_{i,j}^{-1}(c) - \widetilde \tau_{i,j} (a)
	- \widetilde \kappa_{ij}) {\mathfrak q}_{i,j}^{a, \eta_{i,j}^{-1}(c)
	- \widetilde \tau_{i,j} (a) - \widetilde \kappa_{ij}}\ 
	(i \trid j, c) \otimes (i, a) \\
&= \sum_{a, c \in A} q_{ij}\;\phi(\widetilde\kappa_{ij})^{-1} \; \chi_{i,j}
	(\widetilde \kappa_{i,j})^{-1} \quad \psi(a)\; \phi
	(\widetilde \tau_{i,j} (a))^{-1} \; \mu_{i,j} (a) \; \chi_{i,j}
	(\widetilde \tau_{i,j} (a))^{-1}  \times \\
&\qquad\qquad\qquad \times \phi(\eta_{i,j}^{-1}(c)) \chi_{i,j}
	(\eta_{i,j}^{-1}(c)) \; (i \trid j, c) \otimes (i, a);
\end{align*}
where in the first equality we use \eqref{fourier-cociclodinamico};
in the second, we perform the change of variables
$c = \eta_{i,j} (b) + \tau_{i,j} (a) + \kappa_{ij}$ which gives
$$b = \eta_{i,j}^{-1} \left(c - \tau_{i,j} (a) - \kappa_{ij}\right) =
\eta_{i,j}^{-1}(c) - \widetilde \tau_{i,j} (a) - \widetilde
\kappa_{ij};$$
in the third equality we use \eqref{fourier-cociclo-condicion}
and that $\phi$ and $\chi_{ij}$ are multiplicative. The claim follows.
\epf

\begin{rem} The derived rack of the braided set underlying
\eqref{fourier-formula} is given by
$$
(i, \psi) \trid (j,\phi) = (i\trid j, [(\phi\chi_{ij}) \circ \eta_{ij}^{-1}]
\; [(\psi\chi_{i\trid j, i})\circ \widetilde \tau_{i\trid j, i}]^{-1}
\; \mu_{i\trid j, i}).$$
Note that it \emph{does not} depend on $(\kappa_{ij})$
but only on the cocycle $\qg$.
\end{rem}

\begin{exmp}\label{exafourier1}
Assume that $X$ is trivial. 
\begin{enu}
\item  ${\mathfrak q}$ given by
\eqref{fourier-cociclo-condicion} is a cocycle if and only if:
\begin{align}
\label{cociclo-fourier-trivial1}
\chi_{i,k} (\kappa_{j,k}) &= \chi_{j,k} (\kappa_{i,k}) \; 
\mu_{j,k} (\kappa_{i,j}), \\
\label{cociclo-fourier-trivial2}
\chi_{i,k} (\eta_{j,k} (s)) \;  \chi_{j,k} (s) &=
\chi_{j,k} (\eta_{i,k} (s)) \;  \chi_{i,k} (s), \\
\label{cociclo-fourier-trivial3}
1&= \chi_{j,k} (\tau_{i,k} (s)) \; \mu_{j,k} (\tau_{i,j} (s)), \\
\label{cociclo-fourier-trivial4}
\mu_{j,k} (\eta_{i,j} (s)) &= \chi_{i,k} (\tau_{j,k} (s)) \; \mu_{j,k} (s),
\end{align}
for any $i,j,k \in X$, $s\in A$.
\item Let $\sigma: X\to A$ be any function; define
$\kappa_{i,j} := \sigma_{i} - \sigma_{j}$, and
\begin{equation}
\al_{i,j} (a,b) =  b  +   \kappa_{ij},
\end{equation}
that is $\eta_{ij} = \id$, $\tau_{ij} = 0$ in
\eqref{fourier-cociclodinamico}. Then $\al$ is a non-trivial cocycle,
provided that $\sigma$ is not constant; we shall assume this in the
rest of the example and in Lemma \ref{lema-ejuno} below.
\item Furthermore, let $q_{ij} \in \kk^\times$ and let
$\omega: X\to \widehat{A}$ be any function; define
$\chi_{i,j} := \omega_{j} =: \mu_{i,j}^{-1}$, and
define ${\mathfrak q}$ by \eqref{fourier-cociclo-condicion}.  
We claim that ${\mathfrak q}$ is a cocycle. Conditions
\eqref{cociclo-fourier-trivial2}, \eqref{cociclo-fourier-trivial3}
and \eqref{cociclo-fourier-trivial4} follow
because $\eta_{ij} = \id$ and $\tau_{ij} = 0$;
condition \eqref{cociclo-fourier-trivial1} follows
from the special definition of $\kappa_{ij}$.

Then we can apply Theorem 
\ref{mainfourier}. We have
\begin{equation}
 c^{\mathfrak q}\left((i, \psi) \otimes (j, \phi)
\right)= q_{ij}(\phi \omega_{j})(\sigma_{j}) \; (\phi \omega_{j})(\sigma_{i})^{-1}
(j, \phi \omega_{j} ) \otimes (i,
\psi \omega_{j}^{-1}),
\end{equation}
for all $i, j\in X$, $\psi, \phi\in \widehat A$.
In other words, we consider the solution $(X\times \widehat A, S)$
where
\begin{equation}
S \left((i, \psi) , (j, \phi)
\right)= 
\left((j, \phi \omega_{j} ) , (i,
\psi \omega_{j}^{-1}) \right),
\end{equation}
the cocycle $F_{ij}^{\psi, \phi} 
= q_{ij}(\phi \omega_{j})(\sigma_{j}) \; (\phi \omega_{j})(\sigma_{i})^{-1}$,
and the corresponding braided vector space\linebreak
$(\kk(X\times \widehat A), S^F)$. The associated rack is 
$$
(i, \psi) \trid (j, \phi) = (j, \phi \omega_{j}\omega_{i}^{-1}).
$$
\item Assume now that $\omega_{i} = \omega \in \widehat{A}$, for all $i$.
Hence the associated rack is trivial.
Let $$Q_{ij}^{\psi, \phi} =  q_{ij} \phi(\sigma_j) \phi(\sigma_i)^{-1}$$
and let $(\kk (X\times \widehat A), c^Q)$ be the associated braided vector
space. Let $T^n: (X\times \widehat A)^n \to (X\times \widehat A)^n$
be as in the proof of Proposition \ref{slyz1}.
In our case, we have
$$
T^n \left((i_1, \psi_1), \dots, (i_j, \psi_j), \dots, (i_n, \psi_n)\right) =
\left((i_1, \psi_1 \omega^{1 - n}), \dots, (i_j, \psi_j \omega^{j - n}), 
\dots, (i_n, \psi_n)\right).
$$
\end{enu}
\end{exmp}

\begin{lem}\label{lema-ejuno} 
The braided vector spaces $(\kk (X\times \widehat A), S^F)$ and 
$(\kk (X\times \widehat A), c^Q)$ are t-equivalent.
\end{lem}

\pf Let $p_h$ be defined inductively by $p_1 = 1$, $p_{h+1} = p_h + n - h$,
and let 
$$
\lambda_{i_1, \dots, i_n} = \prod_{1\le h \le n} \omega^{p_h} (\sigma_{i_h}).
$$
We shall show that the  map $U^n: \kk(X\times \widehat A)^{\otimes n} \to 
\kk(X\times \widehat A)^{\otimes n}$, 
$U^n \left((i_1, \psi_1) \otimes \dots \otimes
(i_n, \psi_n)\right) =  \lambda_{i_1, \dots, i_n}  T^n \left((i_1, \psi_1)
\otimes \dots \otimes (i_n, \psi_n)\right)$, 
satisfies $U^n S^F_{j, j+1} = c^Q_{j,
j+1}U^n$;
that is, $U^n$ intertwines the corresponding representations of the braid group.

On one hand,
\begin{align*}
U^n S^F_{j, j+1}& \left((i_1, \psi_1) \otimes \dots \otimes(i_n, \psi_n)\right) \\
&= q_{i_{j}i_{j + 1}}(\psi_{j + 1} \omega)(\sigma_{i_{j + 1}}) 
	\; (\psi_{j + 1} \omega)(\sigma_{i_{j}})^{-1} \times \\
&\quad\times U^n \left((i_1, \psi_1 \omega^{1 - n})\otimes \dots\otimes
	(i_{j + 1}, \psi_{j + 1} \omega), (i_{j}, \psi_{j} \omega^{- 1}) \dots,
	(i_n, \psi_n)\right) \\
&= q_{i_{j}i_{j + 1}}(\psi_{j + 1} \omega)(\sigma_{i_{j + 1}}) \; (\psi_{j + 1} \omega)
	(\sigma_{i_{j}})^{-1} \lambda_{i_1, \dots, i_{j + 1}, i_{j}, \dots,  i_n} \times \\
&\quad\times \left((i_1, \psi_1 \omega^{1 - n})\otimes \dots\otimes
	(i_{j + 1}, \psi_{j + 1} \omega^{j + 1 - n}), (i_{j}, \psi_{j} \omega^{j  - n})
	\otimes \dots\otimes (i_n, \psi_n)\right);
\end{align*}
whereas, on the other hand,
\begin{align*}
c^Q_{j, j+1} U^n &\left((i_1, \psi_1) \otimes \dots \otimes (i_n, \psi_n)\right) \\
&= \lambda_{i_1, \dots, i_n} c^Q_{j,
j+1}
\left((i_1, \psi_1 \omega^{1 - n})\otimes \dots\otimes (i_j, \psi_j \omega^{j - n}),
(i_{j + 1}\otimes \psi_{j + 1} \omega^{j + 1 - n})\otimes 
\dots\otimes (i_n, \psi_n)\right)
\\
&= \lambda_{i_1, \dots, i_n} q_{i_{j}i_{j + 1}} 
\psi_{j + 1} \omega^{j + 1 - n}(\sigma_{i_{j + 1}})
\psi_{j + 1} \omega^{j + 1 - n}(\sigma_{i_{j}})^{-1} \times
\\ &\quad\times 
\left((i_1, \psi_1 \omega^{1 - n})\otimes \dots\otimes 
(i_{j + 1}, \psi_{j + 1} \omega^{j + 1 - n}), (i_j, \psi_j \omega^{j - n})\otimes 
\dots\otimes (i_n, \psi_n)\right).
\end{align*}
The equality holds since
$$\omega (\sigma_{i_{j + 1}}) 
\;  \omega(\sigma_{i_{j}})^{-1} \omega^{p_j} (\sigma_{i_{j+ 1}})
\omega^{p_{j + 1}} (\sigma_{i_{j }})
= \omega^{p_j} (\sigma_{i_{j}}) \omega^{p_{j + 1}} (\sigma_{i_{j + 1}})  
 \omega^{j + 1 - n}(\sigma_{i_{j + 1}})
 \omega^{-j - 1 + n}(\sigma_{i_{j}}),$$
by definition of the $p_h$'s.
\epf

\begin{exmp}\label{exafourier2}
Assume that $\eta_{ij} = \id$ for all $i,j\in X$. 
\begin{enu}
\item\label{exf21}
A family $(\tau_{ij})$ defines a quandle structure of $X$-module on $A$
if and only if $\tau_{ii} = 0$, $\tau_{j,k} = \tau_{i\trid j, i\trid k}$ and 
\begin{equation}\label{eq0-cubo-perm}
\tau_{i, j\trid k} =  \tau_{i,  k} + \tau_{i, j} \tau_{j, k},
\end{equation}
for all $i,j, k\in X$. Given such a family, 
$(\kappa_{ij})$ is a $2$-cocycle if and only if  
\begin{equation}\label{eq1-cubo-perm}
\kappa_{j, k} + \kappa_{i, j\trid k} = \kappa_{i, k} 
+ \tau_{j,k}(\kappa_{i,  j}) + \kappa_{i\trid j,i\trid  k},
\end{equation}
\item We shall consider the following family of examples:
$X = (\Z/3, \trid^2)$ is the unique simple crossed set with $3$ elements;
$A$ is a finite abelian group of exponent $2$;
$\tau_{ij} = \id - \delta_{ij}$, $i,j \in X$.
It is not difficult to verify that  
$(\tau_{ij})$ satisfies the conditions in (\ref{exf21}).
We fix $a\in A$, and set
$$
\kappa_{ij} =  (1 - \delta_{ij})\;a, \qquad \wkappa_{ij} = \delta_{i, j+ 1}\; a,
$$
$i,j \in X$. Both families $(\kappa_{ij})$ and $(\wkappa_{ij})$
satisfy \eqref{eq1-cubo-perm}; we denote the corresponding extensions
by $Y = X \times_{\kappa} A$, $\widetilde Y = X \times_{\wkappa} A$.
We shall assume that $a \neq 0$. 
Indeed, $Y = X \times_{\kappa} A$ (for $a \neq 0$) is isomorphic to
$Y = X \times_{0} A$ (case $a = 0$); just consider the function
$f:X\to A$, $f(i) = a$ for all $i$, 
and check that $\kappa$ is cohomologous to $0$, \emph{cf.}
\eqref{eq:2cch}.

Analogously, we denote $\widehat Y = X \times_{0} \widehat A$.
\item If $A = \Z/2$, $a=1$, $Y$ is isomorphic
to the crossed set of transpositions in $\mathbb S_4$,
via the identification
\begin{align*}
&(12) = (0, 0), \quad (34) = (0, 1), \quad (13) = (1, 0),\\
&(24) = (1, 1), \quad (14) = (2, 0), \quad (23) = (2, 1).
\end{align*}
Furthermore, 
$\widetilde Y$ is isomorphic
to the crossed set of $4$-cycles in $\mathbb S_4$
(that is, the faces of the cube),
via the identification
\begin{align*}
(1234) &= (0, 0), \quad (1324) = (1, 0), \quad
(1243) = (2, 0), \\
(1432) &= (0, 1), \quad (1423) = (1, 1),  \quad (1342) = (2, 1).
\end{align*}
The crossed sets $Y$ and $\widetilde Y$ are not isomorphic.
\item We consider now the cocycle
${\mathfrak q}_{ij}^{ab} = q_{ij} \in \kk^\times$; that is
$\chi_{i,j} = \mu_{i,j} = \varepsilon$, for all $i,j$,
in  \eqref{fourier-cociclo-condicion}.  
By Theorem \ref{mainfourier}, we have in the braided vector spaces
$(\kk Y, c^{\qg})$ and $(\kk \widetilde Y, \widetilde c^{\qg})$
the equalities
\begin{align*}
c^{\mathfrak q}\left((i, \psi) \otimes (j, \phi) \right)
& = q_{ij}\; \phi ((1 - \delta_{ij})\;a) \,
(i\trid j, \phi) \otimes (i, \psi \; \phi^{1 - \delta_{ij}}), \\
\widetilde c^{\mathfrak q}\left((i, \psi) \otimes (j, \phi) \right)
& = q_{ij}\; \phi ((\delta_{i, j + 1})\;a) \,
(i\trid j, \phi) \otimes (i, \psi \; \phi^{1 - \delta_{ij}}),
\end{align*}
for all $i, j\in X$, $\psi, \phi\in \widehat A$.
That is, we have isomorphisms with braided vector spaces\linebreak
$(\kk (X \times \widehat A), S^{F})$,
respectively $(\kk (X \times \widehat A), \widetilde  S^{F})$.
In both cases, the associated rack is given by
$(i, \psi) \trid (j, \phi) = (i\trid j, \phi \psi^{1- \delta_{ij}})$;
this is the crossed set $\widehat Y$. Let 
$Q_{ij}^{\psi\phi} = q_{ij}$. 
\end{enu}
\end{exmp}

\begin{lem}\label{lema-ejdos} 
\begin{enumerate}
\item\label{ejf21}
The braided vector spaces $(\kk Y, c^{\qg})$ and 
$(\kk \widehat Y, c^Q)$ are t-equivalent.
\item\label{ejf22}
The braided vector spaces $(\kk \widetilde Y, \widetilde c^{\qg})$
and $(\kk \widehat Y, c^Q)$ are t-equivalent.
\item The braided vector spaces $(\kk Y, c^{\qg})$ and 
$(\kk \widetilde Y, \widetilde c^{\qg})$ are t-equivalent.
\end{enumerate}
\end{lem}

\pf
\begin{enu}
\item By Theorem \ref{mainfourier}, it is enough to show that the  map
\begin{align*}
&U^n: \kk(X\times \widehat A)^{\otimes n} \to \kk(X\times \widehat A)^{\otimes n}, \\
&U^n \left((i_1, \psi_1) \otimes\dots\otimes (i_n,\psi_n)\right)
	=(\psi_1 \dots \psi_n)(a)\,  T^n \left((i_1, \psi_1)
		\otimes\dots\otimes (i_n, \psi_n)\right),
\end{align*}
satisfies $U^n S^F_{j, j+1} = c^Q_{j,
j+1}U^n$. This is a straightforward computation.  
\item Let $\Gamma_n$ be the image of the group homomorphism
$\rho^n: \mathbb B_n \to\mathbb S_{X^n}$ induced by the rack
structure on $X$. Let $\Lambda_n := \sum_{g\in \Gamma_n} g$
be a non-normalized integral of the Hopf algebra $\kk \Gamma_n$.
The group $\Gamma_n$ acts on the vector space
$\fun (X^n, \kk X)$ in the usual way; let
$\eta^n = \Lambda_n\cdot \delta_1$ where $\delta_1$
is the function  $\delta_1 (i_1, \dots, i_n) = i_1$.
We write
$$\eta^n = \sum_{k \in K^n}  \eta^n_k,$$
where $\eta^n_k$ is actually a function from $X^n$ to $X$
and $K^n$ is an index set. Let 
\begin{align*}
R^n_1(i_1, \dots, i_{n}) 
&= \left(
\sum_{k \in K^{n-1}} \delta_{i_1 \trid \eta_k^n(i_2, \dots, i_n), \, i_1 + 1}
\right) a, \\
R^n_t(i_1, \dots, i_{n}) 
&= R^n_1(i_t, i_t \trid i_1, i_t \trid i_2, \dots, i_t \trid i_{t-1},
i_{t+1}, \dots,  i_{n}) , \quad t \ge 2.
\end{align*}
We consider the map $U^n: \kk(X\times \widehat A)^{\otimes n} \to 
\kk(X\times \widehat A)^{\otimes n}$, 
$$U^n \left((i_1, \psi_1) \otimes \dots \otimes
(i_n, \psi_n)\right) =  \psi_1(R^n_1) \dots \psi_n(R^n_n) 
\, T^n \left((i_1, \psi_1)
\otimes \dots \otimes (i_n, \psi_n)\right),$$
By Theorem \ref{mainfourier}, it is enough to show 
that $U_n$ satisfies $U^n S^F_{j, j+1} = c^Q_{j, j+1}U^n$. 
A straightforward computation shows that this is 
equivalent to the following set of identities:

\begin{align}
\label{dificil1}
R^n_t(i_1, \dots, i_{n}) &= 
R^n_t(i_1, \dots, i_h\trid i_{h+1}, i_h, \dots, i_{n}), \qquad h\neq t, t+1; \\
\label{dificil2}
R^n_t(i_1, \dots, i_{n}) &= 
R^n_{t + 1}(i_1, \dots, i_t\trid i_{t+1}, i_t, \dots, i_{n}); 
\end{align}
\begin{multline}
\label{dificil3}
R^n_t(i_1, \dots, i_{n}) = \delta_{i_t, i_{t+1} + 1}a
+ R^n_{t - 1}(i_1, \dots, i_{t - 1}\trid i_{t}, i_{t - 1}, \dots, i_{n})
\\
\qquad\qquad + \left(1 -   \delta_{i_t, i_{t - 1}}\right)
R^n_t(i_1, \dots, i_{t - 1}\trid i_{t}, i_{t - 1}, \dots, i_{n}).
\end{multline}

Now, equation \eqref{dificil1} for $t=1$ follows from the invariance
of the integral, whereas for $t>1$ follows from the definition
and the case $t=1$. Clearly, \eqref{dificil2} follows from the definition
also. Finally, \eqref{dificil3} can be shown by induction
on $n$ and $t$.
\item follows from (\ref{ejf21}) and (\ref{ejf22}).
\end{enu}
\epf

%%%%%%%%%%%%%%%%%%%%%%%%%%%%%%%%%%%%%%%%%%%%%%%%%%%%%%%%%%%%%%%
%%%%%%%%%%%%%%%% Seccion  Nichols alg. %%%%%%%%%%%%%%%%%%%%%%%%
%%%%%%%%%%%%%%%%%%%%%%%%%%%%%%%%%%%%%%%%%%%%%%%%%%%%%%%%%%%%%%%
\section{Nichols algebras and pointed Hopf algebras}\label{sn:6}

\subsection{Definitions and tools}

The Nichols algebra of a rigid braided vector space $(V,c)$
can be defined in various different ways, see for example \cite{AG,AS2}. 
We retain the following one.
If we consider the symmetric group $\Sim_n$ and the braid group $\Tre_n$
with standard generators $\{\tau_1,\ldots,\tau_{n-1}\}$
and $\{\sgm_1,\ldots,\sgm_{n-1}\}$ respectively,
then the so-called Matsumoto section for the canonical
projection $\Tre_n\to\Sim_n$, $\sgm_i\mapsto\tau_i$, is the
set-theoretical function defined on $x\in\Sim_n$ by the recipe
(i) write $x=\tau_{i_1}\cdots\tau_{i_l}$ in a shortest possible way, and
(ii) replace the $\tau_i$'s by $\sgm_i$'s, i.e., $M(x)=\sgm_{i_1}\cdots\sgm_{i_l}$.
Then,
$$\toba(V) = \oplus_{n \ge 0} \toba^n(V) = \kk \oplus V \oplus 
\left(\oplus_{n \ge 2} T^{n}(V)/ \ker Q_n\right),$$ where 
$Q_n = \sum_{x \in \mathbb S_n} M(x)$
is the so-called ``quantum symmetrizer". 
This presentation of the Nichols algebra immediately implies:

\begin{lem}\label{basic} If $(V,c)$ and $(\widetilde V,\widetilde c)$ are t-equivalent
braided vector spaces (\emph{cf.} Definition \ref{t-equivalence-def})
then the corresponding Nichols algebras $\toba(V)$
and $\toba(\widetilde V)$ are isomorphic as graded vector spaces.
In particular, one has finite dimension, resp. finite GK-dimension, if
and only if the other one has.\hfill\qed
\end{lem}

For a subspace $J\subseteq T(V)$ we say that it is
\emph{compatible with the braiding} if $c(V\otimes J)=J\otimes V$ and
$c(J\otimes V)=V\otimes J$.

\begin{lem}\label{lm:icct}%ideal compatible con la trenza
Let $(X,\trid)$ be a rack, let $q\in\kk^\times$ and let $\qg\equiv q$ be
the cocycle $\qg_{ij}=q\ \forall i,j\in X$.
Let $(V=\kk X,c^{\qg})$ be the associated braided space and let $J\subseteq T(V)$
be a subspace. Notice that $T(V)$ is an $\gax$-comodule
algebra  with the structure induced by
$V\to\kk\gax\otimes V$, $x\mapsto\phi_x\otimes x$ for $x\in X$. Furthermore,
$T(V)$ is an $\gax$-module algebra via $\trid$. If $J$ is $\gax$-homogeneous
and $\gax$-stable, then it is compatible with the braiding.
\end{lem}
\pf
It is sufficient to prove that for $x\in X$ we have
$c(J\otimes\kk x)\subseteq V\otimes J$ and $c(\kk x\otimes J)\subseteq J\otimes \kk x$.
The first inclusion is a consequence of the homogeneity of $J$, the second
one is a consequence of the stability of $J$.
\epf

Finite dimensional Nichols algebras (as well as any finite dimensional
graded rigid braided Hopf algebra) satisfy a Poincar\'e duality: let $n$ be
the degree of the space of integrals (it is easy to see that the space
of integrals is homogeneous with respect to the $\Z$-grading),
then $\dim\toba^r(V)=\dim\toba^{n-r}(V)$ for all $r\in\Z$. Furthermore,
since $\toba(V)$ is concentrated in positive degrees, we have that
$\toba^m(V)=0$ for $m>n$; since $\dim\toba^0(V)=1$ we have that
$\dim\toba^n(V)=1$; since $\toba(V)$ is generated by $\toba^1(V)$, we have
that $\dim\toba^r(V)\neq 0$ for $0\le r\le n$.
We call $n$ the \emph{top degree of $\toba(V)$}.
Choose then a nonzero integral $\int$. There is a non-degenerate bilinear
pairing (which is the same as that in the proof of the Poincar\'e duality)
given by $(x|y)=\lambda$ if $xy=\lambda\int+\text{ terms of degree $<n$}$.
These facts, first encountered by Nichols \cite{n},  give a powerful
strategy for computing Nichols algebras. We state this strategy
after the following definition.

\begin{defn}\label{de:jyq}
For $r\ge 2$, let $J_r$ be the ideal generated by $\oplus_{i=2}^r\ker(Q_i)$.
Let $\wtoba_r(V)=T(V)/J_r$, which has a projection $\wtoba_r\to\toba(V)$.
It is not difficult to see that $\oplus_{i=2}^r\ker(Q_i)$ is a coideal which
is compatible with the braiding, whence $\wtoba_r(V)$ is a braided Hopf algebra.
Moreover, it is graded, it is generated by its elements in degree $1$,
and in degree $0$ it is $1$-dimensional. Then it fulfills the
same properties above about Poincar\'e duality as $\toba(V)$.
\end{defn}

\begin{thm}\label{strate}
\begin{enumerate}
\item\label{str1}
Suppose that $\wtoba_r(V)$ vanishes in degree $2r+1$.
Then $\wtoba_r(V) = \toba(V)$.
\item\label{str2}
Let $J\subseteq \ker(T(V)\to\toba(V))$ be an ideal which is also a coideal
and is compatible with the braiding. Suppose that $T(V)/J$ is finite
dimensional, it has top degree $n$ and $\dim\toba^n(V)\neq 0$.
Then $T(V)/J=\toba(V)$.
\item
Suppose that $\wtoba_r(V)$ is finite dimensional, it has top degree $n$ and
$\dim\toba^n(V)\neq 0$. Then $\wtoba_r(V)=\toba(V)$.
\end{enumerate}
\end{thm}

\pf
\begin{enu}
\item Follows from Poincar\'e duality: $\dim\wtoba^i_r(V)=\dim\toba^i(V)$ for
$0\le i\le r$ and the top degree of $\toba(V)$ is $\le 2r$.
\item
This is so thanks to the non-degenerate bilinear form of $T(V)/J$: let $\int$
be an integral in $T(V)/J$.
%then $\int$ projects to an integral in $\toba(V)$.
If $0\neq x\in \ker(T(V)/J\to\toba(V))$, then there exists $y\in T(V)/J$ such
that $xy=\int$. But then $\int\in\ker(T(V)/J\to\toba(V))$, which implies that
$\im(T^n(V)/J\to\toba(V))=0$, a contradiction.
\item follows from (\ref{str2}).
\end{enu}
\epf

In the examples we present in subsections \ref{ssn:ce1} and \ref{ssn:ce2}, 
we have computed the quotients $T(V)/J$
finding Gr\"obner bases, with the help of \cite{opal}. We also used a program
in Maple with subroutines in C to find the dimensions of the ideals $J_r$ in degree
$r$ for small $r$'s. Generators of $J_r$ have been found by hand, using differential
operators (see below). Thus, we have used part (\ref{str1}) of the Theorem.
However, we shall use part (\ref{str2}) in the proofs.

The best way to prove that certain elements vanish (or not) in
a Nichols algebra is given by the differential operators
$\partial_{x^*}:\toba^n(V)\to\toba^{n-1}(V)$ ($x^*\in V^*$).
These are skew derivations. When $V=\kk X$ is given by a rack, we consider
the basis $\{x\in X\}$ of $V$ and $\{x^*\}$ its dual basis; we extend
$x^*$ to $\toba(V)$ by $x^*(\alpha)=0$ if $\alpha\in\toba^n(V)$, $n\neq 1$.
We put then $\partial_x=\partial_{x^*}=(\id\otimes x^*)\circ\Delta$
(here $\Delta$ is the comultiplication in $\toba(V)$).
%then, defining
%$\beta_x(\alpha)=(\id\otimes x^*)c(x\otimes\alpha)$, we put
%$\partial_x=\partial_{x^*}=x^*$ in degree $1$ and extend it by
%$\partial_x(\alpha\alpha')=\partial_x(\alpha)\beta_x(\alpha')
%	+\alpha\partial_x(\alpha')$.
It can be proved that for $\alpha\in\toba^n(V)$ ($n\ge 2$) we have $\alpha=0$ if
and only if $\partial_{x^*}(\alpha)=0\ \forall x^*\in V^*$ (cf. \cite{n,grfrns}).
We consider analogously defined derivations $\partial_x$ in the algebras $T(V)$,
$\wtoba_r(V)$, $T(V)/J$ for $J$ an ideal as above.

%%%%%%%%%%%%%%%%%%%%%%%%%%%%%%%%%%%%%%%%%%%%%%%%%%%%%%%%%%%%%%%
%%%%%%%%%%%%%%%% SubSeccion  twisting of diagonal %%%%%%%%%%%%%
%%%%%%%%%%%%%%%%%%%%%%%%%%%%%%%%%%%%%%%%%%%%%%%%%%%%%%%%%%%%%%%
\subsection{Some  calculations of Nichols algebras related to braided 
vector spaces of diagonal group type}

In this subsection we compute examples of Nichols algebras 
of braided vector spaces arising from Example \ref{exafourier1}.
We first recall some results on Nichols algebras of braided spaces
of diagonal group type, i.e., $(V,c)=(\kk X,c^\qg)$, where $X$ is
a trivial rack.

\begin{prop}\label{cartan} Let $(V,c)$ be a braided
vector space with $V = \kk x_1 \oplus \cdots \oplus \kk x_{\theta}$ and 
$$
c(x_i \otimes x_j) = q_{ij} x_j \otimes x_i, \qquad 1 \le i,j \le \theta.
$$

\begin{enumerate}
\item Assume that $q_{ij} = q$ for all $i, j$. Thus

\begin{itemize}
\item (Nichols) If $q= - 1$, then $\toba (V) \simeq \Lambda (V)$, 
hence  $\dim \toba (V) = 2^{\theta}$.

\item (Lusztig; see \cite{AS2}) 
If $q$ is a primitive third root of 1, then  $\dim \toba (V) = 27$
when $\theta = 2$; and $\dim \toba (V) = \infty$ when $\theta > 2$.

\item (Lusztig; see \cite{AS2}) 
If $\ord q > 3$ and $\theta \ge 2$, then   $\dim \toba (V) = \infty$.
\end{itemize}
\item\cite{AD}\label{cart2}.
Assume that $q_{ii}=-1$ $\forall i$, $q_{ij}\in\{\pm 1\}\ \forall i,j$.
For $i\neq j$, set $A_{ij}\in\{0,-1\}$ such that $q_{ij}q_{ji}=(-1)^{A_{ij}}$.
Set also $A_{ii}=2$. Then $(A_{ij})_{1 \le i,j \le \theta}$
is a simply laced generalized Cartan matrix. Thus

\begin{itemize}
\item  If the components of the Dynkin diagram
corresponding to $(A_{ij})$ are of type $A_m$ (not necessarily the 
same $m$), then  $\dim \toba (V) = 2^{\vert \Phi^+ \vert}$, where
$\Phi^+$ is the set of positive roots corresponding to $(A_{ij})$.

\item 
If the the Dynkin diagram corresponding to $(A_{ij})$
contains a cycle, then $\dim\toba(V)=\infty$. \hfill\qed
\end{itemize}
\end{enumerate}
\end{prop}

\begin{conj}\cite{AD}. Same notation as in (\ref{cart2}) above. Then 
$\dim \toba (V) = 2^{\vert \Phi^+ \vert}$ if $(A_{ij})$ is of
finite type, and $\dim \toba (V) = \infty$ otherwise.
\end{conj}
Let us consider the  Dynkin diagrams:
$$
\begin{array}{cccccccccccc}
&& \bullet &&& &&& \bullet && \\
&& \vert   &&&&&& \vert   &&\\
\bullet & \textbf{------} & \bullet & \textbf{------} & \bullet &,&
	\bullet & \textbf{------} & \bullet & \textbf{------} & \bullet\quad .\\
&&&&&&&& \vert   && \\
&&&&&&&& \bullet && \\ \\
&& D_4 &&&&&& D_4^{(1)} &&
\end{array}
$$

We shall need to assume only the following weaker form 
of the previous conjecture: 
\begin{conj}\label{debil} Same notation as in (\ref{cart2}) above. Then 
$\dim \toba (V) = 2^{12}$ if $(A_{ij})$ is of
 type $D_4$, and $\dim \toba (V) = \infty$ if $(A_{ij})$ is of
 type $D_4^{(1)}$.
\end{conj}

We next consider a trivial rack $X$, a finite, non-trivial, abelian group
$A$, denoted multiplicatively; and a non-constant function $\sigma: X \to A$. 
We set $\kappa_{ij} := \sigma_i\sigma_j^{-1}$, $\eta_{ij}=\id$,
$\tau_{ij}=0\ \forall i,j\in X$ (\emph{cf.} Example \ref{exafourier1}).
Let $Y = X\times A$, let
$(q_{ij})_{i,j\in X}$ be a collection of scalars,
let $\omega \in \widehat{A}$, let
$\chi_{i,j} := \omega =: \mu_{i,j}^{-1}$, for all $i,j$;
define ${\mathfrak q}$ by \eqref{fourier-cociclo-condicion}.

\begin{prop}\label{cartan-aplicaciones} 
Let $(V,c)$ be the braided vector space $(\kk Y, c^\qg)$.
\begin{enumerate}\itemsep 6pt
\item\label{carapa} If $\ord q_{ii} > 3$ for some $i$, 
then   $\dim \toba (V) = \infty$.  
\item\label{carapb} If $\ord q_{ii} = 3$ for some $i$, 
then either $A \simeq C_2$, the group of order $2$ (hence $27$ divides 
$\dim \toba (V)$ if this is finite), or else $\dim \toba (V) = \infty$.
\item\label{carapc}
Assume that $q_{ii} = -1$ for all $i \in X$.
Furthermore, assume that $q_{ij}q_{ji} = (-1)^{A_{ij}}$ for all 
$i \neq j\in X$, where $A_{ij}\in \{0,-1\}$; and that $\ord \kappa_{ij} \le 2$. 
Then, if $A \not\simeq C_2$, we have $\dim\toba(V) = \infty$.
\item\label{carapd} Same hypotheses as in (\ref{carapc}); assume that 
the group $A \simeq \{\pm 1 \}\simeq C_2$.
Let $X_{\pm} = \{i\in X: \sigma_i = \pm 1 \}$. 
\emph{Assume that Conjecture \ref{debil} is true}. 
Then
\begin{itemize}
\item  If $\card X_+ = 1$, $\card X_- \le 3$ and 
$q_{ij}q_{ji} = 1$ for all $i \neq j \in X_-$,
then $\dim \toba (V)$ is finite.
\item  If $\card X_- = 1$, $\card X_+ \le 3$ and 
$q_{ij}q_{ji} = 1$ for all $i \neq j \in X_+$,
then $\dim \toba (V)$ is finite.
\item In all other cases,  $\dim \toba (V) = \infty$.
\end{itemize}
\end{enumerate}
\end{prop}

\pf
Let $\widehat Y = X\times \widehat A$ and let 
$Q^{\psi, \phi}_{ij} = q_{ij}\phi(\kappa_{ij}^{-1})$.
By Theorem \ref{mainfourier} and Lemma \ref{lema-ejuno}, 
it is enough to consider the braided vector space 
$(W, c) = (\kk \widehat Y, c^Q)$.

Then (\ref{carapa}) and (\ref{carapb}) follow from Proposition \ref{cartan}, 
applied to the the subspace $\kk \widehat Y_i$,
where $Y_i = \{(i,\psi): \psi \in \widehat A \}$.
The divisibility claim in (\ref{carapb}) follows from Proposition
\ref{pr:afr} below.

We prove (\ref{carapc}). 
There exists $i\neq j$ such that $\kappa_{ij}$ has order $2$
(recall that $\sigma$ is not constant). Let
$E = \{ \phi \in \widehat A: \phi(\kappa_{ij}) = 1 \}$
and 
$F = \{ \phi \in \widehat A: \phi(\kappa_{ij}) = -1 \}$;
clearly, $\card E = \card F$. Assume that $\card E > 1$,
and let $\phi_1 \neq \phi_2 \in E$, $\psi_1 \neq \psi_2 \in F$.
There are two possibilities:
\begin{itemize}\itemsep 6pt
\item  If $q_{ij}q_{ji} = 1$, then 
$(i, \phi_1)$, $(i, \phi_2)$, $(j, \psi_1)$, $(j, \psi_2)$
span a subspace $U$ of $W$ with $c(U\otimes U) = U\otimes U$;
the associated Cartan matrix is of type $A_3^{(1)}$, hence
$\dim \toba (V) = \infty$.
\item  If $q_{ij}q_{ji} = -1$, then 
$(i, \phi_1)$, $(i, \phi_2)$, $(j, \phi_1)$, $(j, \phi_2)$
span a subspace $U$ of $W$ with $c(U\otimes U) = U\otimes U$;
the associated Cartan matrix is of type $A_3^{(1)}$, hence
$\dim \toba (V) = \infty$.
\end{itemize}

We prove (\ref{carapd}). Let $i\neq j \in X_+$; then $\kappa_{ij} = 1$. 
Let us denote $\widehat A =\{ \varepsilon, \sgn\}$. 
If $q_{ij}q_{ji} = -1$, then 
$(i, \varepsilon)$, $(i, \sgn)$, $(j, \varepsilon)$, $(j, \sgn)$
span a braided vector subspace of Cartan type with matrix
$A_3^{(1)}$; by Proposition \ref{cartan} (\ref{carapb}), $\dim \toba (V) = \infty$.

We assume then that $q_{ij}q_{ji} = 1$ for all $i\neq j \in X_+$,
and also for all $i\neq j \in X_-$. Since $\kappa$ is non-trivial,
both $X_+$ and $X_-$ are non-empty. Let $i \in X_+$ and consider the 
vector subspace $U$ spanned by $(\{i\} \cup X_-)\times \widehat A$. If 
$\card X_- > 3$, then the Cartan matrix of the braiding of $U$ contains 
a principal submatrix of type $D_4^{(1)}$. If Conjecture \ref{debil}
is true, then $\dim \toba (V) = \infty$. Hence, we can assume
that $\card X_{\pm} \le 3$. Also, if $\card X_{\pm}=2$ then
the Cartan matrix of  $W$ contains a cycle; 
by Proposition \ref{cartan} (\ref{carapb}), $\dim \toba (V) = \infty$.
The only cases left are 
$\card X_+ = 1$, $\card X_- \le 3$, or
$\card X_- = 1$, $\card X_+ \le 3$. In these cases,
the Cartan matrices of $W$ are of finite type; either
$A_2 \times A_2$, or $A_3 \times A_3$, or $D_4 \times D_4$.
This concludes the proof of (\ref{carapd}). \epf

\begin{rem}
\begin{enu}
\item In part (\ref{carapb}) of the proposition, if  $\ord q_{ii} = 3$, 
for some $i\in X$, and  $\card A =  2$, then our 
present knowledge of Nichols algebras of diagonal type does
not allow to obtain any general conclusion on $\dim \toba(V)$.
\item In part (\ref{carapc}) of the proposition, if  $\ord q_{ij}q_{ji} > 2$, 
or $\ord \kappa_{ij} > 2$ for some $i, j\in X$, then our 
present knowledge of Nichols algebras of diagonal type does
not allow to obtain any conclusion on $\dim \toba(V)$.
\end{enu}
\end{rem}

%%%%%%%%%%%%%%%%%%%%%%%%%%%%%%%%%%%%%%%%%%%%%%%%%%%%%%%%%%%%%%%
%%%%%%%%%%%%%%%% SubSeccion Examples twisting   %%%%%%%%%%%%%%%
%%%%%%%%%%%%%%%%%%%%%%%%%%%%%%%%%%%%%%%%%%%%%%%%%%%%%%%%%%%%%%%
\subsection{Concrete realizations of pointed Hopf algebras computed
with Fourier transform}\label{ssn:ce1}
Here we give examples of groups with Yetter--Drinfeld modules as in
Proposition \ref{cartan-aplicaciones}. This in turn produces new examples of pointed
Hopf algebras with non-abelian group of grouplikes. We also give a new
pointed Hopf algebra using Example \ref{exafourier2}.

For $n,m\in\N$, let $F=C_{4n}$, $G=C_{4m}$ be cyclic groups. Denote by $x$ a generator of
$F$ and by $y$ a generator of $G$. Let $F$ act on $G$ by $y\prec x=y^{2m+1}$,
and $G$ act on $F$ by $y\succ x=x^{2n+1}$. We can consider then the group
$F\bowtie G$, which coincides with $F\times G$ as a set and whose multiplication
is defined by
$$(x^iy^j)(x^ky^l)=x^i(y^j\succ x^k)(y^j\prec x^k)y^l
	=x^{i+k+jk(2n)}y^{j+l+jk(2m)}.$$
Notice that the center $Z(F\bowtie G)$ is generated by $x^2,y^2$. The conjugacy
classes of $F\bowtie G$ have then cardinality $1$ and $2$, and they are
$$\{x^{2i}y^{2j}\},\ 
\{x^{2i+1}y^{2j},x^{2i+2n+1}y^{2j+2m}\},\ 
\{x^{2i}y^{2j+1},x^{2i+2n}y^{2j+2m+1}\},\ 
\{x^{2i+1}y^{2j+1},x^{2i+2n+1}y^{2j+2m+1}\}.$$
Take $\chi$ defined by $\chi(x)=\chi(y)=-1$, and the Yetter--Drinfeld
modules $V_i=M(x^{2i+1},\chi)$, $W_j=M(y^{2j+1},\chi)$.
We have then the following examples:
\begin{enumerate}\itemsep 12pt
\item $V=V_0\oplus W_0$. By Proposition \ref{cartan-aplicaciones},
	$V$ is t-equivalent to a space of type $A_2\times A_2$,
	and then $\dim\toba(V)=2^6$.
	We have a family of link-indecomposable pointed Hopf algebras
	$\toba(V)\#\kk(F\bowtie G)$ of dimension $2^6\times 16mn=2^{10}mn$.
	The smallest example of this family is $n=m=1$ with dimension $2^{10}$.
	Another way to realize this example is over the dihedral group $\Die_4$
	of order $8$, as described in \cite[6.5]{MiS} and \cite[5.2]{gr}.
\item $V=V_0\oplus V_1\oplus W_0$. By Proposition \ref{cartan-aplicaciones},
	$V$ is t-equivalent to a space of type $A_3\times A_3$, and then
	$\dim\toba(V)=2^{12}$.
	We have a family of link-indecomposable pointed Hopf algebras
	$\toba(V)\#\kk(F\bowtie G)$ of dimension $2^{12}\times 16mn=2^{16}mn$,
	for $n\ge 2$. The smallest example of this family is $n=2,\ m=1$, and
	then $\dim\toba(V)\#\kk(F\bowtie G)=2^{12}\cdot 32=2^{17}$.
\item $V=V_0\oplus V_1\oplus V_2\oplus W_0$. By Proposition
	\ref{cartan-aplicaciones}, $V$ is t-equivalent to a space
	of type $D_4\times D_4$. Assuming Conjecture \ref{debil},
	we have $\dim\toba(V)=2^{24}$.
	We have a family of link-indecomposable pointed Hopf algebras
	$\toba(V)\#\kk(F\bowtie G)$ of dimension $2^{24}\times 16mn=2^{28}mn$,
	for $n\ge 3$. The smallest example of this family is $n=3,\ m=1$, and
	then $\dim\toba(V)\#\kk(F\bowtie G)=2^{24}\cdot 48=2^{28}3$.
\end{enumerate}

\begin{rem}
Actually, in the examples above we have $\dim\mathcal{P}_{g,h}\le 2$ for any
$g,h$ grouplikes. If we allow bigger dimensions, then we can take always
$\Die_4$ as group.
\end{rem}

We give now the algebras obtained from Example \ref{exafourier2}. One of them
appears in \cite{MiS}, and then by Lemma \ref{lema-ejdos} the other one has
the same dimension. We give a full presentation by generators and relations
of both of them.

\subsubsection*{\bf Nichols algebra related to the transpositions in $\Sim_4$ \cite{MiS}}

Let $X= \{a, b, c, d, e, f\}$ be the standard crossed set of the
transpositions in $\Sim_4$ and consider the braided vector space
$(V, c) = (\kk X, c^\qg)$ associated to the cocycle $\qg \equiv -1$.
Here $$ a=(12), \quad b=(13), \quad c=(14), \quad d=(23),
\quad e=(24), \quad f=(34). $$

\begin{thm}\cite[6.4]{MiS}
The Nichols algebra $\toba(V)$ can be presented by generators\linebreak
$\{a, b, c, d, e, f\}$ with defining relations
\begin{equation}\label{eq:rts4}%relaciones toba S_4
\begin{split}
& a^2, \quad b^2,\quad c^2,\quad d^2, \quad e^2, \quad f^2, \\
& dc+cd,\quad eb+be,\quad fa+af,\quad \quad \\
& da+bd+ab,\quad db+ad+ba,\quad ea+ce+ac,\quad ec+ae+ca, \\
& fb+cf+bc,\quad fc+bf+cb,\quad fd+ef+de,\quad fe+df+ed.
\end{split}
\end{equation}
To obtain a basis, choose one element per row below, juxtaposing
them from top to bottom:
\begin{align*}
&(1,a) \\
&(1,b,ba) \\
&(1,c,cb,cba,ca,cab,caba,cabac) \\
&(1,d) \\
&(1,e,ed) \\
&(1,f)
\end{align*}

Its Hilbert polynomial is then
\begin{align*}
P(t) &= (1+t)(1+t+t^2)(1+t+2t^2+2t^3+t^4+t^5)(1+t)(1+t+t^2)(1+t) \\
&= (1+t)^4  (1+t+t^2)^2  (1+t^2)^2  \\ &= t^{12}  +6t^{11}
+19t^{10}  +42t^9 +71t^8 +96t^7 +106t^6 +96t^5 +71t^4 +42t^3
+19t^2 +6t +1.
\end{align*}
Its dimension is $2^63^2=576$. Its top degree is $12$. An integral is given by
$abacabacdedf$.
\end{thm}

We give an alternative proof to that in \cite[6.4]{MiS}
using Theorem \ref{strate}.

\pf
It is straightforward to see that the elements in \eqref{eq:rts4} vanish
in $\toba(V)$. One can either use differential operators or either compute
$Q_2=1+c$ on them. Using Gr\"obner bases it can be seen that if $J$ is the
ideal generated by these relations, then $T(V)/J$ has the stated basis.
Since $J$ is generated by primitive elements, it is a coideal.
Furthermore, $J$ is generated by homogeneous elements with respect to the
$\gax$-grading, and it is invariant under the $\gax$-action. By Lemma
\ref{lm:icct}, it is compatible with the braiding and then $T(V)/J$ is
a braided Hopf algebra. Last, we show that $abacabacdedf$
does not vanish in $\toba(V)$: it is straightforward to see that
$$\partial_a\partial_b\partial_a\partial_c\partial_a\partial_b
	\partial_a\partial_c\partial_d\partial_e\partial_d\partial_f
	(abacabacdedf)=1$$
Now, we conclude by Theorem \ref{strate} part (\ref{str2}).
\epf

To realize this example, one can take $G=\Sim_4$, $g=(1\,2)$,
$\chi=\sgn$. Then $V=M(g,\chi)\in\ydskg$ is isomorphic to $(\kk X,c^\qg)$. We get
a pointed Hopf algebra $\toba(V)\#\kk G$ whose dimension is
$576\times 24=2^93^3$. One can construct also a family of
link-indecomposable pointed Hopf algebras taking as group $\Sim_4\times C_m$,
where $C_m$ is the cyclic group of $m$ elements and $m$ is odd.
Let $g=(1\,2)\times t$ ($t$ is a generator of $C_m$) and let the character $\chi$
be the product $\sgn\times\eps$ ($\eps$ is the trivial character).
Then $V=M(g,\chi)$ is again isomorphic to $(\kk X,c^\qg)$ and
we get a pointed Hopf algebra of dimension $2^93^3m$.

\subsubsection*{\bf Nichols algebra related to the faces of the cube}

Let $X= \{a, b, c, d, e, f\}$ be the polyhedral crossed set of the
faces of the cube (that is, the 4-cycles in $\Sim_4$), where
$\{a,f\}$, $\{b,e\}$, $\{c,d\}$ are the pairs of opposite faces and
$a\trid b=c$. Consider the braided vector space $(V,c)=(\kk X,c^\qg)$
associated to the cocycle $\qg \equiv -1$.

By Lemma \ref{basic} and Lemma \ref{lema-ejdos},
the Nichols algebra of $V$ has the same Hilbert series as that
of the preceding example. We can indeed give the precise description
of $\toba(V)$:

\begin{thm}
The Nichols algebra $\toba(V)$ can be presented by generators
$\{a, b, c, d, e, f\}$ with defining relations
\begin{equation}\label{eq:rtcu}%relaciones toba cubo
\begin{split}
& a^2, \quad b^2,\quad c^2,\quad d^2, \quad e^2, \quad f^2, \\
& ec+ce, \quad db+bd, \quad fa+af, \\
& ca+bc+ab, \quad da+cd+ac, \quad eb+ba+ae, \quad fb+ef+be, \\
& fc+cb+bf, \quad fd+dc+cf, \quad fe+ed+df, \quad ea+de+ad.
\end{split}
\end{equation}
%It is then a quadratic algebra; the relations say that opposite
%faces anticommute and there is one more relation for each vertex.
To obtain a basis, choose one element per row below, juxtaposing
them from top to bottom:
\begin{align*}
&(1,a)\\
&(1,b,ba)\\
&(1,c,cb,cba)\\
&(1,d,dc,dcb)\\
&(1,e,ed)\\
&(1,f)
\end{align*}

Its Hilbert polynomial is then
\begin{align*}
P(t) &= (1+t)(1+t+t^2)(1+t+t^2+t^3)(1+t+t^2+t^3)(1+t+t^2)(1+t) \\
&= (1+t)^4  (1+t+t^2)^2  (1+t^2)^2  \\
&= t^{12}  +6t^{11} +19t^{10} +42t^9 +71t^8 +96t^7 +106t^6 +96t^5
	+71t^4 +42t^3 +19t^2 +6t +1.
\end{align*}
Its dimension is $2^63^2=576$. Its top degree is $12$. An integral is given by
$abacbadcbedf$.
\end{thm}

\pf
Again, the elements in \eqref{eq:rtcu} are easily seen to be relations in $\toba(V)$.
Using Gr\"obner bases, it can be seen that if $J$ is the ideal generated by
these elements, then $T(V)/J$ is as stated. We conclude now using Lemma \ref{basic}
and Lemma \ref{lema-ejdos}.
\epf

To realize this example, one can take the $G=\Sim_4$, $g=(1\,2\,3\,4)$,
$\chi=\sgn$. Then $V=M(g,\chi)\in\ydskg$ is isomorphic to $(\kk X,c^\qg)$. We get
a pointed Hopf algebra $\toba(V)\#\kk G$ whose dimension is
$576\times 24=2^93^3$. Also here we get a family of link-indecomposable pointed
Hopf algebras taking the group $\Sim_4\times C_m$ ($m$ odd), $g=(1\,2\,3\,4)\times t$
($t$ a generator of $C_m$) and $\chi=\sgn\times\eps$.

%%%%%%%%%%%%%%%%%%%%%%%%%%%%%%%%%%%%%%%%%%%%%%%%%%%%%%%%%%%%%%%
%%%%%%%%%%%%%%%% SubSeccion Relations affine  %%%%%%%%%%%%%%%%%
%%%%%%%%%%%%%%%%%%%%%%%%%%%%%%%%%%%%%%%%%%%%%%%%%%%%%%%%%%%%%%%
\subsection{Some relations of Nichols algebras of affine racks}

We first present relations in Nichols algebras related to affine racks.
The relation in part (\ref{rea1}) of the following lemma is related to
\cite[5.7]{MiS}; here the rack is more general than there, there the
cocycle is more general than here. Although the relation in part
(\ref{rea3}) below has the same appearence than \cite[(5.24)]{MiS},
the racks and the elements $x,y$ are different.
%Compare with \cite[5.7]{MiS}.

\begin{lem}\label{lm:rea}%relaciones en afines
Let $(A,g)$ be an affine crossed set. Let $\qg\equiv -1$ and
$(V=\kk A,c^\qg)$ the corresponding braided vector space.

\begin{enu}
\item\label{rea1}
Let $x_1,x_2\in A$. Define inductively the elements $x_i\in A$ ($i\ge 3$) by
$x_i=x_{i-1}\trid x_{i-2}$. Let $n$ be the minimum positive integer such that
$x_2-x_1\in\ker\sum_{i=0}^{n-1}(-g)^i$. Then in $\toba(V)$ we have the relation
$$x_2x_1+x_3x_2+\cdots+x_nx_{n-1}+x_1x_n=0.$$
Furthermore, taking different pairs $(x_2,x_1)$, this is a basis of the
relations in degree $2$. In other words, consider in $A\times A$ the relation
$(x_2,x_1)\sim(a,b)$ if there exists $m\in\N$ such that
$(a,b)=(x_m,x_{m-1})$, where the $x_i$'s are defined as above.
Then $\sim$ is an equivalence relation, and the dimension of the
space of relations of $\toba(V)$ in degree $2$ (i.e., the kernel
of the multiplication $V\otimes V\to\toba(V)$) coincides with the
number of equivalence classes $|A\times A/\sim|$.
\item\label{rea2}
If $\qg\equiv q$, where $-q$ is a primitive $\ell$-th root of unity,
then the relations in degree $2$ are the chains
$$x_2x_1+(-q)x_3x_2+(-q)^2x_4x_3+\cdots+(-q)^{n-1}x_1x_n$$
such that $\ell|n$. The elements $x_i\in A$ and $n\in\N$ are defined
as in part (\ref{rea1}).
\item\label{rea3}
If $(1-g+g^2-g^3)(x-y)=0$, then in the Nichols algebra $\toba(V)$
we have the relation
$$x\,y\,x\,y+y\,x\,y\,x=0.$$
The element $xyxy+yxyx$ is homogeneous with respect to the
$\Inn_\trid(A)$-grading. Furthermore, in the algebra $\wtoba_2(V)$
(see Theorem \ref{strate}) the element $xyxy+yxyx$ is primitive.
In other words, if $Q\supseteq Q_2$ is an $\Inn_\trid(A)$-homogeneous
coideal, then the ideal generated by $Q+\kk(xyxy+yxyx)$ is also an
$\Inn_\trid(A)$-homogeneous coideal.
\item\label{rea4}
If $(1-g+g^2)(x-z)=(1-g+g^2)(y-z)=0$ then in the
Nichols algebra $\toba(V)$ we have the relation
$$x\,y\,z\,x\,y\,z+y\,z\,x\,y\,z\,x+z\,x\,y\,z\,x\,y=0.$$
The element $xyzxyz+yzxyzx+zxyzxy$ is homogeneous with respect to the
$\Inn_\trid(A)$-grading. Furthermore, in the algebra $\wtoba_2(V)$
(see Theorem \ref{strate}) the element $xyzxyz+yzxyzx+zxyzxy$ is primitive.
In other words, if $Q\supseteq Q_2$ is an $\Inn_\trid(A)$-homogeneous coideal,
then the ideal generated by $Q+\kk(xyzxyz+yzxyzx+zxyzxy)$ is also an
$\Inn_\trid(A)$-homogeneous coideal.
\end{enu}
\end{lem}
\pf
\begin{enu}
\item
%(\ref{rea1})
It is easy to see by induction that
$$x_t=\sum_{i=0}^{t-2}(-g)^i(x_2-x_1)+x_1.$$
Then $x_n\trid x_{n-1}=x_{n+1}=x_1$ and $x_1\trid x_n=x_{n+2}=x_2$.
It is easy to see that, $g$ being invertible, the chain corresponding
to $x'_1=x_t$, $x'_2=x_{t+1}$ is exactly the same. This is because
$x_{t+1}-x_t=(-g)^{t-1}(x_2-x_1)$.
This proves that the relation $\sim$ is an equivalence relation.
%Let $T=\mdpdc {1-g}g10$, $S=\mdpdc 1g1{-1}$, $\bar S=\frac 1{1+g}\mdpdc 1g1{-1}$.
%We have $(x_{i+1},x_i)^t=T(x_i,x_{i-1})^t$, $S\bar S=\bar SS=1$ 
%and $\bar STS=\mdpdc 100{-g}$.
%Then $(x_{n+2},x_{n+1})^t=S{\mdpdc 100{-g}}^n\bar S(x_2,x_1)^t=(x_2,x_1)^t$.
%This implies that $x_n\trid x_{n-1}=x_1$ and $x_1\trid x_n=x_2$.
On the other hand, the relations in degree $2$ are exactly the kernel of $1+c$.
Thus, we compute
$$(1+c)(x_2x_1+x_3x_2+\cdots+x_1x_n)
	=x_2x_1-x_3x_2+x_3x_2-x_4x_3+\cdots+x_1x_n-x_2x_1=0.$$
Observe that for $x_1=x_2=x$ the minimum $n$ is $1$ and we get the relation $x^2=0$.
To see that these relations generate all the relations in degree $2$, let us take,
for a (not necessarily primitive) $n$-th root of unit $\zeta$, the vector
$$x_2x_1+\zeta x_3x_2+\zeta^2x_4x_3+\cdots+\zeta^{n-1}x_1x_n.$$
It is clear that this vector is an eigenvector of $c$ with eigenvalue $\frac{-1}{\zeta}$.
Thus, each of the strings $(x_1,x_2,\ldots,x_n)$ occurs $n$ times, one for each $n$-th
root of unity, and each of these times with a different eigenvalue. By a dimension argument,
we have diagonalized $c$ and we have picked up the eigenspace associated to $-1$.
\item
%(\ref{rea2})
The same eigenvectors found in the previous part are eigenvectors here, though
their eigenvalues are $q/\zeta$. Thus, for a chain of length $n$ to be a relation,
one must have $q/\zeta=-1$, i.e., $-q$ must be an $n$-th root of unity.
\item
%(\ref{rea3})
Let $z=y\trid x$, $w=z\trid y$. By the previous part, we have in $\toba(V)$
$$x^2=y^2=z^2=w^2=0,\qquad yx+zy+wz+xw=0.$$
Let us apply now $\partial_y$ to the alleged relation. We get
\begin{align*}
\partial_y(xyxy+yxyx) &= xzy+xyx-zyz-yxz = x(zy+yx)-(zy+yx)z \\
	&= -x(wz+xw)+(wz+xw)z=-xwz+xwz=0.
\end{align*}
Analogously, $\partial_x(xyxy+yxyx)=0$. If $a\neq x,\ a\neq y$,
then $\partial_a(xyxy+yxyx)=0$ as well.
This shows that $xyxy+yxyx=0$ in $\toba(V)$, but we have claimed a stronger fact.
To see that $xyxy+yxyx$ is primitive modulo $J_2$, we must prove that
in $T(V)$ we have
$$\Delta(xyxy+yxyx)\in (xyxy+yxyx)\otimes 1
	+1\otimes (xyxy+yxyx)+T(V)\otimes J_2+J_2\otimes T(V).$$
Now, $\Delta$ is a graded map: $\Delta=\oplus_{n,m}\Delta_{n,m}$, where
$\Delta_{n,m}:T^{n+m}(V)\to T^n(V)\otimes T^m(V)$. We must prove then that
the images of $xyxy+yxyx$ by $\Delta_{1,3}$, $\Delta_{2,2}$ and $\Delta_{3,1}$
lie in $T(V)\otimes J_2+J_2\otimes T(V)$. The previous argument, with derivations,
shows that $\Delta_{3,1}(xyxy+yxyx)\in J_2\otimes V$. For the others, let us
introduce the following notation: if $m\in\N$, we take the basis
$\{x_1\cdots x_m\ |\ x_i\in A\ \forall i\}$ of $T^m(V)$. Let
$\{(x_1\cdots x_m)^*\ |\ x_i\in A\ \forall i\}$ be the dual basis, and
let $\partial_{x_1\cdots x_m}=(\id\otimes(x_1\cdots x_m)^*)\circ\Delta$.
These maps are skew differential operators of degree $m$ and for $m=1$
they coincide with the derivations. We have then for $W\in T^n(V)$,
$$\Delta_{n-m,m}(W)=\sum_{x_1\in A,\ldots,x_m\in A}
	\partial_{x_1\cdots x_m}(W)\otimes(x_1\cdots x_m).$$
We prove now that
$\Delta_{2,2}(xyxy+yxyx)\in T(V)\otimes J_2+J_2\otimes T(V)$.
Clearly, $\Delta_{ab}(xyxy+yxyx)=0$ unless $\{a,b\}\subseteq\{x,y\}$.
Since $xx$ and $yy$ are in $J_2$, it is sufficient to see that the image of
$xyxy+yxyx$ by $\partial_{xy}$ and $\partial_{yx}$ lies in $J_2$.
Let $t=x\trid y$ and $s=t\trid x$. Then $st+tx+xy+ys\in Q_2$. Furthermore,
$x\trid z=x\trid(y\trid x)=(x\trid y)\trid (x\trid x)=t\trid x=s$. Now, it is
straightforward to check that $\partial_{xy}(xyxy+yxyx)=st+tx+xy+ys\in J_2$.
The computation for $\partial_{yx}$ is analogous.
Finally, it is easy to see that $\partial_{xyx}(xyxy+yxyx)=\partial_{yxy}(xyxy+yxyx)=0$,
which proves that $\Delta_{1,3}(xyxy+yxyx)\in T(V)\otimes J_2+J_2\otimes T(V)$.

It remains to be proved that $xyxy+yxyx$ is $\Inn_\trid(A)$-homogeneous.
This is equivalent to prove that $\phi_x\phi_y\phi_x\phi_y=\phi_y\phi_x\phi_y\phi_x$.
Take $a\in A$. We have $\phi_x\phi_y(a)=(1-g)x+g(1-g)y+g^2a$, and then
\begin{align*}
\phi_x\phi_y\phi_x\phi_y(a)&=(1+g^2)(1-g)x+g(1-g)(1+g^2)y+g^4a \\
	&=x-g^4y+g(1-g+g^2)(y-x)+g^4a=x-g^4y+(y-x)+g^4a \\
	&=y-g^4y+g^4a
\end{align*}
Analogously, $\phi_y\phi_x\phi_y\phi_x(a)=x-g^4x+g^4a$. Then
$$\phi_x\phi_y\phi_x\phi_y(a)-\phi_y\phi_x\phi_y\phi_x(a)=(1-g^4)(y-x)=0.$$
\item
%(\ref{rea4})
Let us define the following elements in $A$:
\begin{align*}
&h=y\trid z=		(1-g)y	+gz  		\\
&s=x\trid y=(1-g)x	+gy  			\\
&t=x\trid h=(1-g)x	+y	-(1-g)z  	\\
&r=x\trid z=(1-g)x		+gz  		\\
&b=x\trid(y\trid r)=x	+y	-z  		\\
\end{align*}
One can check that
$t=s\trid r$,
$b=y\trid t$.
It is straightforward to check that any two of these elements satisfy that
their difference lies in the kernel of $1-g+g^2$. This is so because each of
these is an affine combination of $x,y,z$ whose parameters are polynomials
in $g$ (and any such polynomial leaves $\ker(1-g+g^2)$ invariant).
By the first part, we have the following relations in $\wtoba_2(V)$:
\begin{equation}\label{rr3}
\begin{split}
&x^2=y^2=z^2=h^2=s^2=t^2=r^2=0 \\
&hy+yz+zh= sx+xy+ys= tx+xh+ht= rx+xz+zr= 0 \\
&sr+rt+ts= by+yt+tb= zsz-szs= 0
\end{split}
\end{equation}
Notice that for any two of these elements, say $x_1,x_2$, if we put
$x_3=x_2\trid x_1$ we then get\linebreak
$x_2x_1+x_3x_2+x_1x_3=0$, and then, since
$x_1^2=x_2^2=0$, we have $x_1x_2x_1=x_2x_1x_2$. This explains the
relation $zsz-szs=0$ above.
As in the previous case, we must prove that the image of
$xyzxyz+yzxyzx+zxyzxy$ by $\Delta_{1,5}$, $\Delta_{2,4}$, $\Delta_{3,3}$,
$\Delta_{4,2}$ and $\Delta_{5,1}$ lies in $J_2\otimes T(V)+T(V)\otimes J_2$.
This is a very long computation, but it is straightforward and we give only
two examples: for $\Delta_{5,1}$ we apply $\partial_x$ and for $\Delta_{3,3}$
we apply $\partial_{xyx}$ and $\partial_{yxy}$.
Let us call $W=xyzxyz+yzxyzx+zxyzxy$. We have
\begin{align*}
&\partial_x(W)=-srxsr+xyzsr-yzsrx+yzxyz+zsrxs-zxyzs, \\
&\partial_{xyx}(W)=txb+yzt-ztx, \\
&\partial_{yxy}(W)=xhb-hby-zxh.
\end{align*}
It can be seen that relations \eqref{rr3} imply that the first of these elements
lies in $J_2$. The second and third elements do not lie in $J_2$; however, since modulo
$J_2$ we have $xyx=yxy$, in the image by $\Delta_{3,3}$ we have
\begin{align*}
&\partial_{xyx}(W)\otimes xyx+\partial_{yxy}(W)\otimes yxy \\
	&\hspace*{2cm}=(txb+yzt-ztx+xhb-hby-zxh)\otimes xyx\text{ modulo }T(V)\otimes J_2.
\end{align*}
Now, it can be seen that relations \eqref{rr3} imply that
$txb+yzt-ztx+xhb-hby-zxh$ lie in $J_2$.
The proof that $W$ is $\Inn_\trid(A)$-homogeneous is analogous to that of $xyxy+yxyx$
being homogeneous in part (\ref{rea3}).
\end{enu}
\epf

\begin{rem}
If $A=\fpt$, $g$ is the multiplication by $w\neq -1$, then
the minimum $n$ in part (\ref{rea1}) of the Lemma is always the order of $-w$ as a root
of unit in $\fpt$, except for $x_1=x_2=x$. We have then exactly
$p^t+\frac{p^{2t}-p^t}n$ independent relations in degree $2$, and therefore
the dimension of $\toba^2(V)$ is $\frac{n-1}n(p^{2t}-p^t)$.
If $w=-1$, we have the same result with $n=p=\operatorname{char}(\fpt)$.
Furthermore, if $n=4$ then we can apply part (\ref{rea3}) for any two
elements $x,y\in A$. If $n=3$ then we can apply part (\ref{rea4}) for
any three elements $x,y,z\in A$.
\end{rem}

%%%%%%%%%%%%%%%%%%%%%%%%%%%%%%%%%%%%%%%%%%%%%%%%%%%%%%%%%%%%%%%
%%%%%%%%%%%%%%%% SubSeccion Examples affine  %%%%%%%%%%%%%%%%%%
%%%%%%%%%%%%%%%%%%%%%%%%%%%%%%%%%%%%%%%%%%%%%%%%%%%%%%%%%%%%%%%
\subsection{Examples of Nichols algebras and pointed Hopf algebras on affine
racks}\label{ssn:ce2}
We present here two examples. In both we have relations given
by Lemma \ref{lm:rea}.

\subsubsection*{\bf Nichols algebra related to the vertices of the
	tetrahedron \cite{gr}}
Let $X=\{1,2,3,4\}$ be the polyhedral crossed set of the vertices of
the tetrahedron and consider the braided vector space $(V,c)=(\kk X,c^\qg)$
associated to the cocycle $\qg \equiv -1$.

\begin{thm}
The Nichols algebra $\toba(V)$ can be presented by generators
$\{1,2,3,4\}$ with defining relations
\begin{equation}\label{eq:rttet}%relaciones toba tetrahedro
\begin{split}
& 1^2,\quad 2^2,\quad 3^2,\quad 4^2, \\
& 31 + 23 + 12,\quad 41 + 34 + 13,\quad 42 + 21 + 14,\quad 43 + 32 + 24, \\
% Para Grobner: 212 - 121, 323 - 232
& 321321 + 213213 + 132132.
\end{split}
\end{equation}
To obtain a basis, choose one element per row below, juxtaposing
them from top to bottom ($e$ is the unit element):
\begin{align*}
&(e,1) \\
&(e,2,21) \\
&(e,321) \\
&(e,3,32) \\
&(e,4)
\end{align*}
\end{thm}

Its Hilbert polynomial is then
$P(t) = t^9 + 4 t^8 + 8 t^7 +11 t^6 +12 t^5 + 12t^4 + 11t^3 +8t^2 +4t +1$.
Its dimension is $72$, its top degree is $9$, an integral is given by
$121321324$.

\pf
As explained in Remark \ref{rm:wia}, the tetrahedron crossed set coincides
with the affine crossed set $(\mathbb{F}_4,w)$, where $w^2+w+1=w^2-w+1=0$.
The cases (\ref{rea1}) with $n=3$ and (\ref{rea4}) of Lemma \ref{lm:rea} apply
immediately and we see that the elements in \eqref{eq:rttet} are relations in
$\toba(V)$. Let $J$ be the ideal generated by these elements.
It can be seen that $J$ is $\gax$-stable. Since by Lemma \ref{lm:rea} part
(\ref{rea4}) the element $321321+213213+132132$ is $\gax$-homogeneous and
$Q_2$ is compatible with the braiding, then $J$ is compatible with the
braiding. Moreover, by Lemma \ref{lm:rea} part (\ref{rea4}) again, it is
a coideal.
Now, it is straightforward to see that
$$\partial_1\partial_2\partial_3\partial_4\partial_2
	\partial_4\partial_3\partial_4(121321324)=2\in\kk.$$
We now use   Theorem \ref{strate} part (\ref{str2}).
\epf

To realize this example, one can take the affine group
$\mathbb{F}_4\rtimes\mathbb{F}_4^\times\simeq\mathbb{A}_4$ and its direct
product with $C_2$. That is, we take $G=\mathbb{A}_4\times C_2$. Denote
by $t$ the generator of $C_2$ and let $g=(1\,2\,3)\times t\in G$.
Take $\chi\in\hat G$, $\chi(\sigma\times t^i)=(-1)^i$. Then
$V=M(g,\chi)\in\ydskg$ is isomorphic to $(\kk X,c^\qg)$. We get
a pointed Hopf algebra $\toba(V)\#\kk G$ whose dimension is
$72\times 24=2^63^3$. We get a family of link-indecomposable pointed
Hopf algebras replacing $C_2$ by $C_m$ ($m$ even), $g=(1\,2\,3)\times t$
($t$ a generator of $C_m$) and $\chi(\sigma\times t^i)=(-1)^i$.

\subsubsection*{\bf Nichols algebra related to the affine crossed set $(\Z_5,\trid^2)$}

Let $X= \{0,1,2,3,4\}$ be the affine crossed set $(\Z_5,\trid^2)$ and
consider the braided vector space $(V, c) = (\kk X, c^\qg)$
associated to the cocycle $\qg \equiv -1$. That is,
$c(i \otimes j) = - (2j-i)\otimes i$.

\begin{thm}
The Nichols algebra $\toba(V)$ can be presented by generators
$\{0,1,2,3,4\}$ with defining relations
\begin{equation}\label{eq:rta52}%relaciones toba afin 5-2
\begin{split}
& 0^2, \quad  1^2, \quad  2^2, \quad  3^2, \quad  4^2, \\
& 32+20+13+01, \quad 40+21+14+02, \quad 41+34+10+03, \\
& 42+30+23+04, \quad 43+31+24+12, \\
& 1010+0101.
\end{split}
\end{equation}
To obtain a basis, choose one element per row below, juxtaposing
them from top to bottom ($e$ is the unit element):
\begin{align*}
&(e, 0) \\
&(e, 1, 10, 101)\\
&(e, 2, 21, 212, 20, 201, 2012, 2010, 20102, 201020)\\
&(e, 3, 31, 312, 30, 303, 3031, 30312)\\
&(e, 4)\\
\end{align*}

Its Hilbert polynomial is then
\begin{align*}
P(t) &= (1 + t)^2  (1 + t + t^2  + t^3 )
	(1 + t + 2 t^2 + 2 t^3 + 2 t^4 + t^5 + t^6)(1+t+2t^2+2t^3+t^4+t^5) \\
&= t^{16}   + 5 t^{15}  + 15 t^{14}   + 35 t^{13}   + 66 t^{12}
	+ 105 t^{11} + 145 t^{10}   + 175 t^{9}  + 186 t^{8}  + 175 t^{7} \\
&\hspace{1cm} + 145 t^{6} + 105 t^{5}  + 66 t^{4}  + 35 t^{3}  + 15 t^{2}  + 5 t + 1.
\end{align*}
Its dimension is $1280$. Its top degree is $16$.
An integral is given by $0101201020303124$.
\end{thm}
\pf
The relations are given by Lemma \ref{lm:rea} parts (\ref{rea1}) and (\ref{rea3}).
By the same result, if $J$ is the ideal generated by \eqref{eq:rta52} then
it is a homogeneous coideal. It is not difficult to see that it is also
$\gax$-stable, whence it is compatible with the braiding. Using Gr\"obner
bases, it can be seen that relations \eqref{eq:rta52} yield the stated
dimensions in each degree. Using Theorem \ref{strate} part (\ref{str2}),
it is sufficient to see that $0101201020303124$ does not vanish in $\toba(V)$
in order to prove the Theorem. It is straightforward then to compute that
$$
\partial_1\partial_0\partial_4\partial_1\partial_4\partial_2
	\partial_3\partial_4\partial_2\partial_4\partial_3\partial_2
	\partial_3\partial_4\partial_3\partial_4(0101201020303124)=1\in\kk.
$$
\epf

To realize this example, one can take the group $G=\Z/5\rtimes(\Z/5)^\times$,
$\chi((a,2^j))=(-1)^j$, $g=(0,2)$. Then $V=M(g,\chi)\in\ydskg$ is
isomorphic to $(\kk X,c^\qg)$. We get a pointed Hopf algebra $\toba(V)\#\kk G$
whose dimension is $1280\times 20=2^{10}5^2$.
We get a family of link-indecomposable pointed Hopf algebras replacing
$(\Z/5)^\times$ by $C_{4m}$, where a generator $t$ of $C_{4m}$ acts as $2$,
i.e., $tit^{-1}=2i$ for $i\in\Z/5$. Then take $g=0\times t$,
$\chi(i\times t^j)=(-1)^j$, $V=M(g,\chi)\in\ydskg$.
The algebra $\toba(V)\#\kk G$ has dimension $2^{10}5^2m$.

%%%%%%%%%%%%%%%%%%%%%%%%%%%%%%%%%%%%%%%%%%%%%%%%%%%%%%%%%%%%%%%
%%%%%%%%%%%%%%%% SubSeccion Freeness %%%%%%%%%%%%%%%%%%%%%%%%%%
%%%%%%%%%%%%%%%%%%%%%%%%%%%%%%%%%%%%%%%%%%%%%%%%%%%%%%%%%%%%%%%
\subsection{A freeness result for extensions of crossed sets}

The concept of extension of crossed sets is not only useful in classification
problems of them. It turns out to be useful as well when one wants to
compute Nichols algebras, as the following Proposition asserts.

Let $X,Y$ be quandles and let $X\to Y$ be a surjective quandle homomorphism.
We assume that $X$ is indecomposable, hence $Y$ is also indecomposable
and $X\simeq Y\times_\al S$ for some dynamical $2$-cocycle
$\al$ and some set $S$ (see Definition \ref{df:dyncoc}).
Let $\qg:Y\times Y\to\rdu$ be a $2$-cocycle. Notice that
$\qg_{ii}=\qg_{jj}$ for all $i,j\in Y.$ Let $\widetilde \qg$ be the
pull-back of $\qg$ along $\pi$, that is $\widetilde \qg_{xy} :=
\qg_{\pi(x),\pi(y)}$. Let $(V,c)=(\kk X,c^{\widetilde \qg})$,
$(V',c')=(\kk Y, c^\qg)$. Let $P_V(t)$ be the Hilbert series of
$\toba(V)$ and $P_{V'}(t)$ the Hilbert series of $\toba(V')$.

\begin{prop}\label{pr:afr}
\begin{enumerate}
\item\label{afr1}
If the order of $\qg_{ii}$ is $>3$ and $\card S\ge 2$ then
$\dim\toba(V)=\infty$. If the order of $\qg_{ii}$ is $3$ and $\card
S\ge 3$ then $\dim\toba(V)=\infty$.
\item\label{afr2}
Let $P_S(t)$ be the Hilbert polynomial of $\toba(W)$, where
$(W,c_W)=(\kk S,c^{\qg_{ii}})$ (i.e., the cocycle is the constant
$\qg_{ii})$). Then $P_S\,|\,P_V$.
\item\label{afr3}
If $\alpha$ is a constant cocycle, $P_{V'}\,|\,P_V$.
\end{enumerate}
\end{prop}
\pf (\ref{afr1}) follows easily from \cite[Lemma 3.1]{gr}. (\ref{afr2}) follows at
once from \cite[Theorem 3.8,1]{grfrns}. (\ref{afr3}) can be proved using
a remark  right after the proof of  \cite[Theorem 3.2]{MiS}.
Actually, this remark is a generalization of Theorem 3.2 in {\it
loc. cit.}, which in turn is a generalization of \cite[Theorem
3.8,1]{grfrns}. The remark goes as follows: let $(R,c),\,(R',c')$
be braided Hopf algebras with maps $R'\fllad{i}R\fllad{\phi}R'$ of
algebras and coalgebras such that $\phi i=\id$, and such that
\begin{equation}\label{eq:csl}
(i\otimes\id)c'(\phi\otimes\id)=(\id\otimes\phi)c(\id\otimes i),
\quad c(i\phi\otimes\id)=(\id\otimes i\phi)c.
\end{equation}
Let $\rho:R\to
R\otimes R',\ \rho=(\id\otimes\phi)\Delta_R$, and let
$K=R^{\text{co} R'}=\{r\in R\,|\,\rho(r)=r\otimes 1\}$. Then the
conditions on $i$ and $\phi$ are sufficient to prove that
$\mu:K\otimes R'\to R,\ \mu=m_R(\id\otimes i)$ is an isomorphism.

Thus, we can find $i:\toba(V')\to\toba(V)$ and
$\phi:\toba(V)\to\toba(V')$ satisfying the previous conditions. By
the definitions of Nichols algebras, to give an algebra and
coalgebra map it is enough to give the maps at degree $1$ and
verify that they commute with the braidings. That is,
$V'\fllad{i}V\fllad{\phi}V'$ such that
$c(i\otimes i)=(i\otimes i)c'$ and similarly with $\phi$. We take
$i(y) = \frac 1{|S|}\sum_{\pi(x)=y}x$, 
$\phi(x) = \pi(x)$.
It is immediate to see that $i$ and $\phi$ commute with the
braidings, using that $\pi$ is a map of crossed sets and
$\widetilde \qg=\pi^{-1}(\qg)$.
The conditions in \eqref{eq:csl} are also easy to verify; for
the second one it is used that $\alpha$ is a constant cocycle.
\epf

%%%%%%%%%%%%%%%%%%%%%   Acknowledgments   %%%%%%%%%%%%%%
\subsection*{Acknowledgments}
We thank R. Guralnick for kindly communicating us Theorem \ref{clasigur2}.
We are grateful to P.~Etingof and J.~Alev for interesting discussions, as well
as to J.S.~Carter, A.~Harris and M.~Saito for pointing out several misprints
and statements which needed some clarification in earlier versions of the manuscript.
Our deep gratitude is due to F.~Fantino for carefully reading the paper
and sharing with us his comments.
We thank S.~Natale for encouraging us with Section \ref{sn:6}.

Results of this paper were announced at the XIV Coloquio latinoamericano
de \'Algebra, La Falda, august 2001, by the second-named author.
Part of the work of the first author was done during a visit 
to the University of Reims (October 2001 - January 2002); 
he is very grateful to J. Alev for his kind hospitality.

\end{document}